\documentclass[10pt,reqno]{amsart}
\usepackage{amsmath}
\usepackage{amsthm}
\usepackage{amsfonts}
\usepackage{amssymb}
\usepackage{graphicx}
\usepackage{bm}
\usepackage{multirow}
\usepackage{nomencl}
\newtheorem{theorem}{Theorem}[section]

\newtheorem{corollary}[theorem]{Corollary}

\newtheorem{definition}[theorem]{Definition}
\newtheorem{example}[theorem]{Example}
\newtheorem{lemma}[theorem]{Lemma}
\newtheorem{notation}[theorem]{Notation}
\newtheorem{proposition}[theorem]{Proposition}
\newtheorem{remark}[theorem]{Remark}
\renewenvironment{proof}{Proof.-}{$\square$}
\begin{document}

\title{Projective and Direct limits of Banach $G$ and tensor structures}
\author{P. Cabau}
\author{F. Pelletier}
%\address{Unit\'{e} mixte de recherche 5127 CNRS\\
%Universit\'{e} de Savoie Mont-Blanc\\
%Laboratoire de Math\'{e}matiques (LaMa)\\
%Campus scientifique\\
%73370 Le Bourget-du-Lac\\
%France}
%\email{patrickcabau@gmail.com, fernand.pelletier@univ-smb.fr}
\maketitle

\begin{abstract}
We endow projective (resp. direct) limits of Banach tensor structures with
Fr\'{e}chet (resp. convenient) structures and study adapted connections to
$G$-structures in both frameworks. This situation is illustrated by a lot of examples.
\end{abstract}


\textbf{M.S.C. 2010.} 53C10, 58B25, 18A30, 22E65.

\tableofcontents

\section{Introduction}

\label{_Introduction}

The concept of $G$-structure provides a unified framework for a lot of
interesting geometric structures. This notion is defined in finite dimension
as a reduction of the frame bundle. Different results are obtained in the
finite case in \cite{Kob} and in \cite{Mol} ($G$-structures equivalent to a
model, characteristic class of a $G$-structure, intrinsic geometry of
$G$-structures, ...).  

\smallskip

In order to extend this notion to the Banach
framework, following Bourbaki (\cite{Bou}), the frame bundle $\ell\left(
TM\right)  $ is defined as an open submanifold of the linear map bundle
$L\left(  M\times\mathbb{M},TM\right)  $ where the manifold $M$ is modelled on
the Banach space $\mathbb{M}$. In \cite{Klo}, Klotz
proofs that the automorphism group of some manifolds can be turned into a
Banach-Lie group acting smoothly on $M$.

\smallskip

The notion of tensor structure which corresponds to an intersection of
$G$-structures, where the different $G$ are isotropy groups for tensors, is
relevant in a lot of domains in Differential Geometry: Krein metric ( cf.
\cite{Bog}), almost tangent and almost cotangent structures (cf. \cite{ClGo}),
symplectic structures (cf.\cite{Vai}, \cite{Wei}), complex structures (cf.\cite{Die}, inner products and decomposable complex structures (cf. \cite{ChMa}).

\smallskip

In this paper, we are interested in the study of projective (resp. direct)
limits of Banach frame bundles sequences $\left(  \ell\left(  E\right)
\right)  _{n\in\mathbb{N}^{\ast}}$, associated tensor structures and adapted
connections to $G$-structures.

\smallskip

The main problem to endow the projective (resp. direct) limit $M\mathbb{=}%
\underleftarrow{\lim}M_{n}$ (resp. $M=\underrightarrow{\lim}M_{n}$) of Banach
manifolds $M_{n}$ modelled on the Banach spaces $\mathbb{M}_{n}$ with a
structure of Fr\'{e}chet (resp. convenient) manifold is to find a chart around
any point of $M$.

\smallskip

It is worth noting that the convenient setting of Kriegl and Michor is
particularly adapted to the framework of direct limit of supplemented Banach
manifolds because the direct limit topology $\tau_{\mathrm{DL}}$ is
$c^{\infty}$-complete (cf. \cite{CaPe}). It appears as a generalisation of
calculus \textit{\`{a} la Leslie} on Banach spaces used in the projective
limit framework.

\smallskip

Because the general linear group $\operatorname{GL}\left(  \mathbb{M}\right)
$ does not admit a reasonable structure in both cases, it has to be replaced by the following groups:

\begin{enumerate}
\item[--] 
The Fr\'{e}chet topological group $\mathcal{H}_{0}\left(
\mathbb{M}\right)  $ (cf. \cite{DGV}), recalled in $\S$ \ref{__TheFrechetTopologicalGroupH0(F)}, in the case of projective sequences;
\item[--] 
The Fr\'{e}chet topological group $\operatorname{G}(\mathbb{E})$ as
introduced in \cite{CaPe}, recalled in $\S$ \ref{__TheFrechetTopologicalGroupG(E)}, in the case of ascending sequences.
\end{enumerate}

So it is possible to endow the projective (resp. direct) limit of a sequence
of generalized linear frame bundles $\left(  \mathbf{\ell}\left(
E_{n}\right)  \right)  _{n\in\mathbb{N}^{\ast}}$, as defined in \cite{DGV},
(resp. linear frame bundles $\left(  \ell\left(  E_{n}\right)  \right)
_{n\in\mathbb{N}^{\ast}}$) with a structure of Fr\'{e}chet (resp. convenient)
principal bundle.

We then obtain similar results for some projective (resp. direct) limits of
sequences of Banach $G$-structures and tensor structures.

\smallskip

This paper is organized as follows.

In section \ref{_BanachStructures}, we recall or introduce some Banach
structures. 
In section \ref{_ExamplesOfTensorStructuresOnABanchBundle}, a lot of examples of tensor structures on a Banach vector bundle is given.
In section \ref{_ExamplesOfIntegrableTensorStructuresOnABanachManifold}, we give examples of integrable tensor structures on a Banach manifolds. 
Section \ref{_ProjectiveLimitsOfTensorStructures} (resp. \ref{_DirectLimitsOfTensorStructures}) is devoted to the study of projective (resp. direct) limit of $G$-structures  and tensor structures.
In section \ref{_ProjectiveLimitsOfAdaptedConnectionsToGStructures} (resp.
\ref{_DirectLimitsOfAdaptedConnectionsToGStructures}), we obtain results about
projective (resp. direct) limits of adapted connections to $G$-structures.\\

In order to have a self-contained paper, we recall the convenient framework in Appendix \ref{_ConvenenientFramework} and different linear operators on Banach spaces in Appendix \ref{_LinearOperators}.

\section{Banach structures}
\label{_BanachStructures}

We give a brief account of different Banach structures which will be used in
this paper. The main references of this section are \cite{Bou}, \cite{DGV},
\cite{Lan}, \cite{Pay} and \cite{Pel1}.

\subsection{Submanifolds}
\label{__Submanifolds}
%\begin{definition}
%\label{D_BanachSubManifolds}A subset $L$ of a Banach manifold $M$ is called a
%submanifold, if for $x\in L$ there exists a chart $\left(  U,\phi\right)  $ of
%$M$ where $x\in U$ such that $\phi\left(  U\cap L\right)  =\phi\left(
%U\right)  \cap\mathbb{L}$, where $\mathbb{L}$ is a closed subspace of the
%model space $\mathbb{M}$.
%\end{definition}
%Then clearly $L$ is itself a Banach manifold with $\left(  U\cap
%L,\phi_{|U\cap L}\right)  $ as charts, where $\left(  U,\phi\right)  $ runs
%through all the submanifold charts from above.
%\begin{definition}
%\label{D_BanachSplittingSubmanifolds}A submanifold $L$ of the Banach manifold
%$M$ is called a splitting submanifold if there exists a cover of $L$ by
%submanifold charts $\left(  U,\phi\right)  $ such that $\mathbb{L}$ is
%supplemented in the Banach space $\mathbb{M}$.
%\end{definition}
Among the different notions of submanifolds, in this paper, we are interested
in the weak submanifolds as used in \cite{Pel1} in the framework of
integrability of distributions.

\begin{definition}
\label{D_WeakSubmanifold} 
Let $M$ be a Banach manifold modelled on the Banach
space $\mathbb{M}$.\\ 
A weak submanifold of $M$ is a pair $(S,\varphi)$ where $S$ is a connected Banach manifold modelled on a Banach space $\mathbb{S}$ and $\varphi:S\longrightarrow M$ is a smooth map such that:

\begin{description}
\item[\textbf{(WSM 1)}] There exists a continuous injective linear map
$i:\mathbb{S}\longrightarrow\mathbb{M}$;

\item[\textbf{(WSM 2)}] $\varphi$ is an injective smooth map and the tangent
map $T_{x}\varphi:T_{x}N\longrightarrow T_{\varphi(x)}M$ is an injective
continuous linear map with closed range for all $x\in N$.
\end{description}
\end{definition}

Note that for a weak submanifold $\varphi:S\longrightarrow M$, on the subset
$\varphi(S)$ of $M$, we have two topologies:

--- the induced topology from $M$;

--- the topology for which $\varphi$ is a homeomorphism from $S$ to
$\varphi(S)$. \\

With this last topology, via $\varphi$, we get a Banach manifold structure on
$\varphi(S)$  modelled on $\mathbb{S}$. Moreover,
the inclusion from $\varphi(S)$ into $M$ is continuous as a map from the
manifold $\varphi(S)$ to $M$. In particular, if $U$ is an open set of $M$,
then $\varphi(S)\cap U$ is an open set for the topology of the manifold on
$\varphi(S)$.

\subsection{$G$-structures and tensor structures on a Banach space}
\label{__GStructuresTensorStructuresOnABanachSpace}

According to \cite{Bel2}, Corollary 3.7, we have the following result.

\begin{theorem}
\label{T_ExistenceWeakLieSubgroup} Assume that $H$ is a Banach-Lie group with the Lie algebra $\mathfrak{h}$, $G$ is a closed subgroup of $H$ and denote
\[
\mathfrak{g}=\{X\in\mathfrak{h}\;|\;\forall t\in\mathbb{R},\ \exp_{H}(tX)\in G\}
\]
Then $\mathfrak{g}$ is a closed subalgebra of $\mathfrak{h}$ and there exist on
$G$ a uniquely determined topology $\tau$ and a manifold structure
making $G$ into a Banach-Lie group such that $\mathfrak{g}$ is the Lie algebra of $G$, the inclusion of $G$ into $H$ is smooth and the following diagram:

\[%
\begin{array}
[c]{ccc}%
\mathfrak{g} & \underrightarrow{\;\;\;} & \mathfrak{h}\\
{\;\exp_{G}\;} \downarrow &  & \downarrow{\;\exp_{H}\;}\\
G & \underrightarrow{\;\;\;} & H
\end{array}
\]
is commutative, where the horizontal arrows stand for inclusion maps.
\end{theorem}

\begin{definition}
\label{D_WeakLieSubgroup} 
A closed topological subgroup $G$ of $H$ which satisfies the assumption of Theorem \ref{T_ExistenceWeakLieSubgroup} is called a weak Lie subgroup of $H$.
\end{definition}

We fix a Banach space $\mathbb{E}_{0}$. For any Banach space $\mathbb{E}$, the Banach
space of linear bounded maps from $\mathbb{E}_{0}$ to $\mathbb{E}$ is denoted
by $\mathcal{L}(\mathbb{E}_{0},\mathbb{E})$ and, if it is non empty, the open
subset of $\mathcal{L}(\mathbb{E}_{0},\mathbb{E})$ of linear isomorphisms from
$\mathbb{E}_{0}$ onto $\mathbb{E}$ will denoted by $\mathcal{L}is(\mathbb{E}%
_{0},\mathbb{E})$. \\
\textit{We always assume that }$\mathcal{L}is(\mathbb{E}_{0},\mathbb{E)}\neq\emptyset$. \\ 
We then have a natural right transitive action
of $\operatorname{GL}(\mathbb{E}_{0})$ on $\mathcal{L}is(\mathbb{E}_{0},\mathbb{E})$ given by
$(\phi,g)\mapsto\phi\circ g$.\\
According to \cite{PiTa}, we define the notion of $G$-structure on the Banach space $\mathbb{E}$.

\begin{definition}
\label{D_GStructureOnBanachSpace} 
Let $G$ be a weak Banach Lie subgroup of
$\operatorname{GL}(\mathbb{E}_{0})$. \\
A $G$-structure on $\mathbb{E}$ is a subset
$\mathcal{S}$ of $\mathcal{L}is(\mathbb{E}_{0},\mathbb{E})$ such that:

\begin{description}
\item[\textbf{(GStrB 1)}] $\forall\left(  \phi,\psi\right)  \in\mathcal{S}%
^{2},\ \phi\circ\psi^{-1}\in G$;

\item[\textbf{(GStrB 2)}] $\forall\left(  \phi,g\right)  \in\mathcal{S}\times
G,\ \phi\circ g$ $\in\mathcal{S}$.
\end{description}
\end{definition}

The set ${L}_{r}^{s}(\mathbb{E}_{0})$ of tensors of type $(r,s)$ on
$\mathbb{E}_{0}$ is the Banach space of $(r+s)$-multilinear maps from
$(\mathbb{E}_{0})^{r}\times(\mathbb{E}_{0}^{\ast})^{s}$ into $\mathbb{R}$.

Each $g\in \operatorname{GL}(\mathbb{E}_{0})$ induces an automorphism $g_{r}^{s}$ of ${L}%
_{r}^{s}(\mathbb{E}_{0})$ which gives rise to a right action of $\operatorname{GL}(E_{0})$ on
${L}_{r}^{s}(\mathbb{E}_{0})$ which is $(\mathbb{T},g)\mapsto g_{r}%
^{s}(\mathbb{T})$ where, for any $\left(  u_{1},\dots,u_{r},\alpha_{1}%
,\dots,\alpha_{s}\right)  \in(\mathbb{E}_{0})^{r}\times(\mathbb{E}_{0}^{\ast
})^{s}$,
\[
g_{r}^{s}(\mathbb{T})\left(  u_{1},\dots,u_{r},\alpha_{1},\dots,\alpha
_{s}\right)  =\mathbb{T}\left(  g^{-1}.u_{1},\dots,g^{-1}.u_{r},\alpha
_{1}\circ g,\dots,\alpha_{s}\circ g\right).
\]

\begin{definition}
\label{D_IsotropyGroupOfTensor}
The isotropy group of a tensor $\mathbb{T}_{0}$ $\in{L}_{r}^{s}(\mathbb{E}_{0})$ is the set
\[
G(\mathbb{T}_{0})=\{g\in \operatorname{GL}(\mathbb{E}_{0})\;:\;g_{r}^{s}(\mathbb{T}%
_{0})=\mathbb{T}_{0}\}.
\]
\end{definition}

\begin{definition}
\label{D_TensorStructureBanachSpace}
Let $\mathbb{T}$ be a $k$-uple $\left(\mathbb{T}_{1},\dots,\mathbb{T}_{k}\right)  $ of tensors on $\mathbb{E}_{0}$.\\
A tensor structure on $\mathbb{E}_{0}$ of type $\mathbb{T}$ or a $\mathbb{T}%
$-structure on $\mathbb{E}_{0}$ is a $G=\bigcap\limits_{i=1}^{k}%
G(\mathbb{T}_{i})$-structure on $\mathbb{E}_{0}$.
\end{definition}

It is clear that such a group $G=\bigcap\limits_{i=1}^{k}G(\mathbb{T}_{i})$ is
a closed topological subgroup of the Banach Lie subgroup of $\operatorname{GL}(\mathbb{E}%
_{0})$. Therefore, from Theorem \ref{T_ExistenceWeakLieSubgroup}, it is a weak
Lie subgroup of $\operatorname{GL}(\mathbb{E}_{0})$.

\subsection{The frame bundle of a Banach vector bundle}
\label{__FrameBundleOfABanachVectorBundle}

Let $\left(  E,\pi_{E},M\right)  $ be a vector bundle of fibre type the Banach
space $\mathbb{E}$, with total space $E$, projection $\pi$ and base $M$
(modelled on the Banach space $\mathbb{M}$). Because $\mathbb{E}$ has not
necessarily a Schauder basis, it is not possible to define the frame bundle of
this vector bundle as it is done in finite dimension.

An extension of the notion of frame bundle in finite dimension to the Banach
framework can be found in \cite{Bou}, $\S$ 7.10.1 and in \cite{DGV}, $\S$ 1.6.5.

The set of linear bicontinuous isomorphisms from $\mathbb{E}$ to $E_{x}$
(where $E_{x}$ is the fibre over $x\in M$) is denoted by $\mathcal{L}is\left(
\mathbb{E},E_{x}\right)  $.

The set
\[
P\left(  E\right)  =\left\{  \left(  x,f\right)  :x\in M,\ f\in\mathcal{L}%
is\left(  \mathbb{E},E_{x}\right)  \right\}
\]
is an open submanifold of the linear map bundle $L\left(  M\times
\mathbb{E},E\right)  =\bigcup\limits_{x\in M}\mathcal{L}\left(  \mathbb{E}%
,E_{x}\right)  $. The Lie group $\operatorname{GL}\left(  \mathbb{E}\right)
$ of the continuous linear automorphisms of $\mathbb{E}$ acts on the right of $P\left(
E\right)  $ as follows:
\[%
\begin{array}
[c]{cccc}%
\widehat{R}: & P\left(  E\right)  \times\operatorname{GL}\left(
\mathbb{E}\right)  & \longrightarrow & P\left(  E\right) \\
& \left(  \left(  x,f\right)  ,g\right)  & \longmapsto & \left(  x,f\right)
.g=\left(  x,f\circ g\right)
\end{array}
.
\]
\begin{definition}
\label{D_FrameBundleOfABanachVectorBundle}
The quadruple $\ell\left(  E\right)  =\left(  P\left(  E\right)
,\pi,M,\operatorname{GL}\left(  \mathbb{E}\right)  \right)  $, where 
\[
\begin{array}
[c]{cccc}%
\pi: & P\left(  E\right)  & \longrightarrow & M\\
& \left(  x,f\right)  & \mapsto & x
\end{array}
\] 
is the projection on the base, is a principal bundle, called the frame bundle of $E$.
\end{definition}

We can write
\[
P\left(  E\right)  =\bigcup\limits_{x\in M}\mathcal{L}is\left(  \mathbb{E}%
,E_{x}\right)  .
\]

Let us describe the local structure of $P\left(  E\right)  $.\\
Let $\lbrace\left(  U_{\alpha},\tau_{\alpha}\right)\rbrace_{\alpha \in A}$ be a local trivialization of $E$
with $\tau_{\alpha}:\pi_{E}^{-1}\left(  U_{\alpha}\right)  \longrightarrow
U_{\alpha}\times\mathbb{E}$. This trivialization gives rise to a local section
$s_{\alpha}:U_{\alpha}\longrightarrow P\left(  E\right)  $ of $P\left(
E\right)  $ as follows:%
\[
\forall x\in U_{\alpha},\ s_{\alpha}\left(  x\right)  =\left(  x,\left(
\overline{\tau}_{\alpha,x}\right)  ^{-1}\right)
\]
where $\overline{\tau}_{\alpha,x}\in\mathcal{L}is\left(  E_{x},\mathbb{E}%
\right)  $ is defined by $\overline{\tau}_{\alpha,x}=\operatorname*{pr}%
_{2}\circ\tau_{\alpha}|_{E_{x}}$.\\
We then get a local trivialization of $P\left(  E\right)  $:
\[%
\begin{array}
[c]{cccc}%
\Psi_{\alpha}: & \pi^{-1}\left(  U_{\alpha}\right)  & \longrightarrow &
U_{\alpha}\times\operatorname{GL}\left(  \mathbb{E}\right) \\
& \left(  x,f\right)  & \mapsto & \left(  x,\overline{\tau}_{\alpha,x}\circ
f\right)
\end{array}
.
\]
In particular, we have
\[
\Psi_{\alpha}\left(  s_{\alpha}\left(  x\right)  \right)  =\left(
x,\overline{\tau}_{\alpha,x}\circ\left(  \overline{\tau}_{\alpha,x}\right)
^{-1}\right)  =\left(  x,\operatorname{Id}_{\mathbb{E}}\right)  .
\]
$\Psi_{\alpha}$ gives rise to $\overline{\Psi}_{\alpha,x}$ defined by
\[
\overline{\Psi}_{\alpha,x}\left(  x,f\right)  =\overline{\tau}_{\alpha,x}\circ
f.
\]
Moreover, for $g\in\operatorname{GL}\left(  \mathbb{E}\right) $, we have
\[
\widehat{R}_{g}\left(  s_{\alpha}\left(  x\right)  \right)  =\left(  x,\left(
\overline{\tau}_{\alpha,x}\right)  ^{-1}\circ g\right)  .
\]

Because the local structure of $P\left(  E\right)  $ is derived from the local
structure of the vector bundle $E$, we get%
\[
\overline{\Psi}_{\alpha,x}\circ\left(  \overline{\Psi}_{\beta,x}\right)
^{-1}\left(  \overline{\tau}_{\beta,x}\circ f\right)  =\overline{\Psi}%
_{\alpha,x}\left(  x,f\right)  =\overline{\tau}_{\alpha,x}\circ f.
\]

We then have the following result.

\begin{proposition}
\label{P_CoincidenceTransitionFunctions-E-PE} The transition functions
\[%
\begin{array}
[c]{cccc}%
\mathrm{T}_{\alpha\beta}: & U_{\alpha}\cap U_{\beta} & \longrightarrow &
\operatorname{GL}\left(  \mathbb{E}\right) \\
& x & \longmapsto & \overline{\tau}_{\alpha,x}\circ\left(  \overline{\tau
}_{\beta,x}\right)  ^{-1}%
\end{array}
\]
of $E$ coincide with the transition functions of $P\left(  E\right)  $.
\end{proposition}

The transition functions form a cocycle; that is%
\[
\forall x\in U_{\alpha}\cap U_{\beta}\cap U_{\gamma},\ \mathrm{T}%
_{\alpha\gamma}\left(  x\right)  =\mathrm{T}_{\alpha\beta}\left(  x\right)
\circ\mathrm{T}_{\beta\gamma}\left(  x\right)
\]

\begin{definition}
\label{D_TangentFrameBundle} The tangent frame bundle $\ell\left(  TM\right)
$ of a Banach manifold $M$ is the frame bundle of $TM$.
\end{definition}

\begin{proposition}
\label{P_BanachBundlesMorphism}Let $\pi_{E_{1}}:E_{1}\longrightarrow M$ and
$\pi_{E_{2}}:E_{2}\longrightarrow M$ be two Banach vector bundles with the
same fibre type $\mathbb{E}$ and let $\Phi:E_{1}\longrightarrow E_{2}$ be a
bundle morphism above $\operatorname{Id}_{M}$. Then $\Phi$ induces a unique
bundle morphism $\ell\left(  \Phi\right)  :$ $\ell\left(  E_{1}\right)
\longrightarrow\ell\left(  E_{2}\right)  $ that is injective (resp.
surjective) if and only if $\Phi$ is injective (resp. surjective).
\end{proposition}

\subsection{$G$-structures and tensor structures on a Banach vector bundle}
\label{__GStructureAndTensorStructuresOnABanachVectorBundle}

A reduction of the frame bundle $\ell\left(  E\right)  $ of a Banach fibre
bundle $\pi_{E}:E\longrightarrow M$ of fibre type $\mathbb{E}$ corresponds to
the data of a weak Banach Lie subgroup $G$ of $\operatorname{GL}\left(
\mathbb{E}\right)  $ and a topological principal subbundle $\left(  F,\pi
_{|F},M,G\right)  $ of $\ell\left(  E\right)  $ such that $\left(  F,\pi
_{|F},M,G\right)  $ has its own smooth principal structure, and the inclusion
is smooth. In fact, such a reduction can be obtained in the following way.

Assume that there exists a bundle atlas $\left\{  \left(  U_{\alpha},\tau
_{a}\right)  \right\}  _{\alpha\in A}$ whose transition functions
\[
\mathrm{T}_{\alpha\beta}:x\longmapsto\overline{\tau}_{\alpha,x}\circ\left(
\overline{\tau}_{\beta,x}\right)  ^{-1}%
\]
belong to a weak Banach-Lie subgroup $G$ of $\operatorname{GL}\left(
\mathbb{E}\right)  $.

To any local trivialization $\tau_{\alpha}$ of $E$, is associated the
following map:
\[
\begin{array}
[c]{cccc}%
\phi_{\alpha}: & U_{\alpha}\times G & \longrightarrow & P\left(  E\right) \\
& \left(  x,g\right)  & \mapsto & s_{\alpha}\left(  x\right)  .g
\end{array}
\]
where
\[
s_{\alpha}\left(  x\right)  .g=\left(  x,\left(  \overline{\tau}_{\alpha
,x}\right)  ^{-1}\right)  .g=\left(  x,\left(  \overline{\tau}_{\alpha
,x}\right)  ^{-1}\circ g\right)  .
\]

If $U_{\alpha}\cap U_{\beta}\neq\phi$, we have%
\[
\forall x\in U_{\alpha}\cap U_{\beta},\ s_{\beta}\left(  x\right)  =s_{\alpha
}\left(  x\right)  .\mathrm{T}_{\alpha\beta}\left(  x\right).
\]

Let $F=\bigcup\limits_{\alpha\in A}V_{\alpha}$ be the subset of $P\left(
E\right) $ where $V_{\alpha}=$ $\phi_{\alpha}\left(  U_{\alpha}\times
G\right) $. 
Because a principal bundle is determined by its cocycles (cf., for
example, \cite{DGV}, $\S$ 1.6.3), the quadruple $\left(  F,\pi_{|F},M,G\right)  $ can
be endowed with a structure of Banach principal bundle and is a topological
principal subbundle of $\ell\left(  E\right) $. 
Moreover, we have the following lemma.

\begin{lemma}
\label{L_ClosedTopologicalSubbundle} 
$F$ is closed in $\ell\left(  E\right)  $
and $F$ is a weak submanifold of $\ell\left(  E\right)  $.
\end{lemma}

\begin{proof}
Consider a sequence $(x_{n},g_{n})\in F$ which converges to some $(x,g)\in
\ell\left(  E\right)  $. Let $\phi_{\alpha}:U_{\alpha}\times
\operatorname{GL}(\mathbb{E})\longrightarrow P\left(  E\right) $ be the mapping associated to a local trivialization of $P(E)$ where $(x,g) \in P(E)$.
Now, for $n$ large enough, $\phi_{\alpha}^{-1}(x_{n},g_{n})$ belongs to
$U_{\alpha}\times G\subset U_{\alpha}\times \operatorname{GL}(\mathbb{E})$. But $U_{\alpha
}\times G$ is closed in $U_{\alpha}\times \operatorname{GL}(\mathbb{E})$. This implies that
$(x,g)$ belongs to $F$. Since the inclusion of $G$ into $\operatorname{GL}(\mathbb{E})$ is
smooth, this implies that the inclusion of $F$ into $\ell\left(  E\right)  $
is also smooth.
\end{proof}

\begin{definition}
\label{D_GStructure}
The weak subbundle $\left(  F,\pi_{|F},M,G\right)  $ of the frame bundle $\ell\left(  E\right)  =\left(  P\left(  E\right),\pi,M,\operatorname{GL}\left(  \mathbb{E}\right)  \right)  $ is called a
$G$-structure on $E$.\newline 
When $E=TM$, a $G$-stucture on $TM$ is called a $G$-structure on $M$.
\end{definition}

Intuitively, a $G$-structure $F$ on a Banach bundle $\pi_{E}:E\longrightarrow
M$ may be considered as a family $\left\{  F_{x}\right\}  _{x\in M}$ which varies smoothly with $x$ in $M$, in the sense that $F=\bigcup\limits_{x\in M}F_{x}$ has a principal bundle structure over $M$
with structural group $G$.

\begin{definition}
\label{D_GStructuresPreservingMorphism}
Let $\pi_{E_{1}}:E_{1}\longrightarrow M$ and $\pi_{E_{2}}:E_{2}\longrightarrow M$ be two Banach vector bundles with the same fibre $\mathbb{E}$ provided with the $G$-structures $\left(
F_{1},\pi_{_{E_{1}}|F_{1}},M,G\right)  $ and $\left(  F_{2},\pi_{_{E_{2}%
}|F_{2}},M,G\right)  $respectively. A morphism $\Phi:E_{1}\longrightarrow
E_{2}$ is said to be $G$-structure preserving if for any $x\in M$,
$\Phi\left(  \left(  F_{1}\right)  _{x}\right)  =\left(  F_{2}\right)  _{x}$.
\end{definition}

\begin{definition}
\label{D_TensorStructureTypeTOnBanachVectorBundle}Let $\mathbb{T}$ be a
$k$-uple $\left(  \mathbb{T}_{1},\dots,\mathbb{T}_{k}\right)  $ of tensors on
the fibre type $\mathbb{E}$ of a Banach vector bundle $\pi_{E}%
:E\longrightarrow M$. A tensor structure on $E$ of type $\mathbb{T}$ is a
$G\left(  \mathbb{T}\right)  =\bigcap\limits_{i=1}^{k}G\left(  \mathbb{T}%
_{i}\right)  $ structure on $E$ where $G\left(  \mathbb{T}_{i}\right)  $ is
the isotropy group of $\mathbb{T}_{i}$ in $\operatorname{GL}\left(
\mathbb{E}\right)  $.
\end{definition}

\begin{definition}
\label{D_TensorBanachBundle}
The Banach fibre bundle $\pi_{\widehat{L_{r}^{s}}}:\widehat{L_{r}^{s}}\left(  E\right)  \longrightarrow M$ whose typical fibre is $L_{r}^{s}\left(  \mathbb{E}\right) $  is the tensor Banach bundle
of type $\left(  r,s\right)  $.
\end{definition}

\begin{definition}
\label{D_TensorLocallyModelled}
Let $\pi_{\widehat{L_{r}^{s}}}:\widehat
{L_{r}^{s}}\left(  E\right)  \longrightarrow M$ be the tensor Banach bundle of
type $(r,s)$.

\begin{enumerate}
\item A smooth (local) section $\mathcal{T}$ of this bundle is called a
(local) tensor of type (r,s) on $E$.

\item Let $\Phi:E\longrightarrow E^{\prime}$ be an isomorphism from a Banach
vector bundle $\pi_{E}:E\longrightarrow M$ to another Banach vector bundle
$\pi_{E^{\prime}}:E^{\prime}\longrightarrow M^{\prime}$ over a map
$\phi:M\longrightarrow M^{\prime}$. Given a (local) tensor $\mathcal{T}$ of
type (r,s) on $E^{\prime}$, defined on an open set $U^{\prime}$ of $M^{\prime
},$ the pullback of $\mathcal{T}$ is the (local) tensor $\Phi^{\ast
}\mathcal{T}$ of the same type defined (on $\phi^{-1}\left(  U^{\prime
}\right)  $) by

$%
\begin{array}
[c]{c}%
\left(  \Phi^{\ast}\mathcal{T}\right)  _{\phi^{-1}\left(  x^{\prime}\right)
}\left(  \Phi^{-1}\left(  u_{1}^{\prime}\right)  ,\dots,\Phi^{-1}\left(
u_{r}^{\prime}\right)  ,\Phi^{\ast}\left(  \alpha_{1}^{\prime}\right)
,\dots,\Phi^{\ast}\left(  \alpha_{s}^{\prime}\right)  \right) \\
=\mathcal{T}_{x^{\prime}}\left(  u_{1}^{\prime},\dots,u_{r}^{\prime}%
,\alpha_{1}^{\prime},\dots,\alpha_{s}^{\prime}\right)
\end{array}
$

\item A local tensor $\mathcal{T}$ is called locally modelled on
$\mathbb{T}\in L_{r}^{s}\left(  \mathbb{E}\right)  $ if there exists a
trivialization $\tau:E_{|U}\longrightarrow U\times\mathbb{E}$ such that%
\[
\mathcal{T}_{|U}\left(  x\right)  =\overline{\tau}_{x}^{\ast}\left(
\mathbb{T}\right)
\]

\end{enumerate}
\end{definition}

\begin{proposition}
\label{P_ConditionTensorStructure}A tensor $\mathcal{T}$ of type $\left(
r,s\right)  $ on a Banach bundle $\pi_{E}:E\longrightarrow M$, where $M$ is a
connected manifold, defines a tensor structure on $E$ if and only if there
exists a tensor $\mathbb{T}\in L_{r}^{s}\left(  \mathbb{E}\right)  $ and a
bundle atlas $\left\{  \left(  U_{\alpha},\tau_{a}\right)  \right\}
_{\alpha\in A}$ such that $\mathcal{T}_{\alpha}=\mathcal{T}_{|U_{a}}$ is
localy modelled on $\mathbb{T}$ and the transition functions $T_{\alpha\beta
}\left(  x\right)  $ belong to the isotropy group $G\left(  \mathbb{T}\right)
$ for all $x\in U_{\alpha}\cap U_{\beta}$ (where $U_{\alpha}\cap U_{\beta
}\not =\emptyset$) and all $\left(  \alpha,\beta\right)  \in A^{2}$.
\end{proposition}

\begin{proof}
Fix some tensor $\mathbb{T}$ on $\mathbb{E}$ and denote the isotropy group of
$\mathbb{T}$ by $G$. Assume that we have a $G$-structure on $E$ and let
$\{\left(  U_{\alpha},\tau_{\alpha}\right)  \}_{\alpha\in A}$ be a bundle
atlas such that each transition function $T_{\alpha\beta}:U_{\alpha}\cap
U_{\beta}\rightarrow \operatorname{GL}(\mathbb{E})$ takes values in $G$. Fix some point
$x_{0}\in M$ and denote an open set of the previous atlas which contains
$x_{0}$ by $U_{\alpha}$. Since we have a trivialization $\tau_{\alpha
}:E_{U_{\alpha}}\rightarrow U_{\alpha}\times\mathbb{E}$, we consider on
$U_{\alpha}$ the section of $\widehat{L_{r}^{s}}(E)$ defined by $\mathcal{T}%
(x)=\bar{\tau}_{\alpha,x}^{\ast}\mathbb{T}$. If there exists $U_{\beta}$ that
contains $x_{0}$, then, on $U_{\alpha}\cap U_{\beta}$, the transition function
$T_{\alpha\beta}$ takes values in $G$; it follows that the restriction of
$\mathcal{T}$ on $U_{\alpha}\cap U_{\beta}$ is well defined. Therefore, there
exists an open set $U$ in $M$ on which is defined a smooth section
$\mathcal{T}$ of $\widehat{L_{r}^{s}}(E)$ such that $\mathcal{T}_{|U_{\alpha
}\cap U}=\bar{\tau}_{\alpha,x}^{\ast}\mathcal{T}$, for all $\alpha$ such that
$U_{\alpha}\cap U\not =\emptyset$. If $x$ belongs to the closure of $U$, let
$U_{\beta}$ be an open set of the previous atlas which contains $x$. Then
using the previous argument, we can extend $\mathcal{T}$ to $U\cup U_{\beta}$.
Since $M$ is connected, we can defined such a tensor $\mathcal{T}$ on $M$.
\newline 
Conversely, if $\mathcal{T}$ is a tensor such that there exists a
bundle atlas $\left\{  (U_{\alpha},\tau_{\alpha})\right\}  _{\alpha\in A}$
such that $\mathcal{T}_{|U_{\alpha}}$ is locally modeled on $\mathbb{T}$, this
clearly implies that each transition function $T_{\alpha\beta}$ must take
values in $G$.
\end{proof}

\subsection{Integrable tensor structures on a Banach manifold}
\label{__IntegrableTensorStructuresOnABanachManifold}

\begin{definition}
\label{D_IntegrableTensorStructureOnABanachManifold}
A tensor $\mathcal{T}$ on a Banach manifold $M$ modelled on the Banach space $\mathbb{M}$ is called an integrable tensor structure if there exists a tensor $\mathbb{T}\in L_{r}^{s}\left(  \mathbb{M}\right)  $ and a bundle atlas $\left\{  \left(  U_{\alpha},\tau_{a}\right)  \right\}  _{\alpha\in A}$
such that $\mathcal{T}_{\alpha}=\mathcal{T}_{|U_{a}}$ is locally modelled on
$\mathbb{T}$ and the transition functions $T_{\alpha\beta}\left(  x\right)  $
belong to the isotropy group $G\left(  \mathbb{T}\right)  $ for all $x\in
U_{\alpha}\cap U_{\beta}$ (where $U_{\alpha}\cap U_{\beta}\not =\emptyset$)
and all $\left(  \alpha,\beta\right)  \in A^{2}$.
\end{definition}

\section{Examples of tensor structures on a Banach bundle}
\label{_ExamplesOfTensorStructuresOnABanchBundle} 
Let $(E,\pi_{E},M)$ be a vector bundle whose fibre is the Banach space $\mathbb{E}$ and whose base is modelled on the Banach space $\mathbb{M}$.

\subsection{Krein metrics}
\label{KreinMetrics}

\begin{definition}
\label{D_PseudoRiemannianMetric} 
A pseudo-Riemannian metric on $(E,\pi_{E},M)$ is a smooth field of weak non degenerate symmetric forms $g$ on $E$. Moreover
$g$ is called
\begin{enumerate}
\item[--] 
a weak Riemannian metric if each $g_{x}$ is a pre-Hilbert inner product on the
fibre $E_{x}$
\item[--] 
a Krein metric if there exists a decomposition $E=E^{+}\oplus E^{-}$ in
a Whitney sum of Banach bundles such that for each fibre $E_{x}$, $g_{x}$ is a
Krein inner product associated to the decomposition $E_{x}=E^{+}_{x}\oplus
E^{-}_{x}$
\item[--]
a neutral metric if it is a Krein metric such that there exists a decomposition  $E=E^{+}\oplus E^{-}$ 
where $E^{+}$ and $E^{-}$ are isomorphic sub-bundles of $E$
\item[--]
a Krein indefinite metric if there exists a decomposition $E=E_{1}\oplus
E_{2}$ in a Whitney sum of Banach bundles such that for each fiber $E_{x}$,
$g_{x}$ is a Krein indefinite inner product associated to the decomposition
$E_{x}=(E_{1})_{x}\oplus(E_{2})_{x}$.
\end{enumerate}
\end{definition}

Given a Krein metric $g$ on $E$, according to Proposition \ref{P_ExtensionKreinInnerProduct}, we get the following result.

\begin{theorem}
\label{T_GStructureKreinMetric}${}$
\begin{enumerate}
\item 
Let $g$ be a Krein metric on a Banach bundle $(E,\pi_{E},M)$. Consider a
decomposition $E=E^{+}\oplus E^{-}$ in a Whitney sum such that the restriction
$g^{+}$ (resp. $-g^{-}$) of $g$ (resp. $-g$) to $E^{+}$ (resp. $E^{-}$) is a
(weak) Riemannian metric. Let $\widehat{E}_{x}$ be the Hilbert space which is the
completion of the pre-Hilbert space $(E_{x},g_{x})$. 
Assume that $\widehat{E}=\bigcup\limits_{x\in M}\widehat{E}_{x}$ is Banach bundle over $M$ such that the inclusion of $E$ into $\widehat{E}$ is a bundle morphism. 
Then $g$ can be extended to a strong Krein metric $\widehat{g}$ on $\widehat{E}$ and we have a
decomposition $\widehat{E}=\widehat{E}^{+}\oplus\widehat{E}^{-}$ such that the
restriction $\widehat{g}^{+}$ (resp. $\widehat{g}^{-}$) of $\widehat{g}$
(resp. $-\widehat{g}$) to $\widehat{E}^{+}$ (resp. $\widehat{E}^{-}$) such
that $\widehat{g}^{+}$ (resp. $-\widehat{g}^{-}$) is a strong Riemannian
metric. In fact $\widehat{E}^{+}$ (resp. $\widehat{E}^{-}$) is the closure of
$E^{+}$ (resp. $E^{-}$) in $\widehat{E}$.\newline 
Given any point $x_{0}$ in
$M$, we identify the typical fibre of $\widehat{E}$ (resp. $E$) with
$\widehat{E}_{x_{0}}$ (resp. $E_{x_{0}}$). Then $\widehat{\gamma}=\widehat
{g}^{+}-\widehat{g}^{-}$ defines a $\widehat{\gamma}_{x_{0}}$-structure and a
$\widehat{g}_{x_{0}}$-structure on $\widehat{E}$. Moreover, $g$ defines a
$g_{x_{0}}$-structure on $E$

\item 
Let $(\widehat{E},\pi_{\widehat{E}},M)$ be a Banach bundle whose typical
fibre is reflexive and let $\widehat{g}$ be a strong Krein metric on
$\widehat{E}$. Then if $\widehat{g}_{x_{0}}$ is the Krein inner product on
$\widehat{E}_{x_{0}}$, then $\widehat{g}$ is a $\widehat{g}_{x_{0}}$-structure
on $\widehat{E}$. \newline 
Let $(E,\pi_{E},M)$ be any Banach bundle such that we have an injective morphism of Banach bundle $\iota: E\longrightarrow\widehat{E}$. Then the restriction $g=\iota^{*}\widehat{g}$ of $\widehat{g}$
is a Krein metric on $E$ and $g$ induces a $g_{x_{0}}$-structure on $E$.
\end{enumerate}
\end{theorem}

\begin{proof}
1. Consider a decomposition $E=E^{+}\oplus E^{-}$ in a Whitney sum such that
the restriction $g^{+}$ (resp. $-g^{-}$) of $g$ (resp. $-g$) to $E^{+}$ (resp.
$E^{-}$) is a (weak) Riemannian metric. \newline

Assume that $\widehat{E}=\bigcup\limits_{x\in M}\widehat{E}_{x}$ is a Banach bundle
over $M$, such that the inclusion $\iota$ of $E$ in $\widehat{E}$ is a
bundle morphism. \newline

This implies that we have a bundle atlas $\{( U_{\alpha},\tau_{\alpha
})\}_{\alpha\in A}$ (resp. $\{(U_{\alpha},\sigma_{\alpha})\}_{\alpha\in
A}$) for $E$ (resp. for $\widehat{E}$) such that, for each $\alpha\in A$, we have
$\sigma_{\alpha}\circ\iota=\tau_{\alpha}$, and, on $U_{\alpha}\cap U_{\beta}$,
the transition function $T_{\alpha\beta}$ is the restriction to $E$ of the
transition function $S_{\alpha\beta}$ associated $\{(U_{\alpha},\tau_{\alpha
})\}_{\alpha\in A}$ to $\mathbb{E}$.\\

Assume that $g$ is a weak Riemannian metric. Then $\widehat{g}$ is a strong
Riemannian metric on $\widehat{E}$.\\
 
Now given any $x \in M$, each fibre $\widehat{E}_{x}$ can be provided with the Hilbert inner product $\widehat{g}_{x}$. Therefore from \cite{Lan}, Theorem 3.1, there exists a bundle atlas
$\{(U^{\prime}_{\alpha},\sigma^{\prime}_{\alpha})\}_{\alpha\in A}$ such that the
transition functions $S^{\prime}_{\alpha\beta}$ take values in the group of
isometries of the typical fibre. 
Now, since $\iota$ is a bundle morphism,
$\{(U^{\prime}_{\alpha}, \tau^{\prime}_{\alpha}=\sigma^{\prime}_{\alpha}\circ\iota)\}_{\alpha \in A}$ is a bundle atlas for $E$ and the transition function $T^{\prime}_{\alpha\beta}$ for this atlas is nothing but the restriction to $\mathbb{E}$ of $S^{\prime}_{\alpha\beta}$ and so is an isometry
of $\mathbb{E}$ provided with the pre-Hilbert inner product induced on $E$
from $\widehat{E}$. Now, for a given $x_{0}\in M$, the typical fibre of
$\widehat{E}$ and $E$ can be identified with $\widehat{E}_{x_{0}}$ and
$E_{x_{0}}$ respectively which ends the proof when $E^{-}$ is reduced to
$\{0\}$. \\
Now, we can apply the previous result to $E^{+}$ and $E^{-}$ and for
the weak Riemannian metric $g^{+}$ and $-g^{-}$ respectively. The result is
then a consequence of Proposition \ref{P_IsometryKreinProduct}.\newline

2. According to the assumption of Point 2, let $\mathbb{E}$ and $\widehat
{\mathbb{E}}$ be the typical fibres of $E$ an $\widehat{E}$, respectively.
Since $\iota:E\longrightarrow\widehat{E}$ is a bundle morphism, as for the
proof of Point 1, we have a bundle atlas $\{(U_{\alpha},\tau_{\alpha})
\}_{\alpha\in A}$ for $E$ and $\{(U_{\alpha},\sigma_{\alpha})\}_{\alpha\in A}$
for $\widehat{E}$ such that for all $\alpha\in A$ we have $\sigma_{\alpha
}\circ\iota=\tau_{\alpha}$ and, on $U_{\alpha}\cap U_{\beta}$ the transition
functions $T_{\alpha\beta}$ are the restrictions to $\mathbb{E}$ of the
transition functions $S_{\alpha\beta}$ associated $\{(U_{\alpha},\tau_{\alpha
})\}_{\alpha\in A}$. Let $\bar{\mathbb{E}}$ be the closure of
$\mathbb{E}$ in $\widehat{E}$. For each $\alpha\in A$ we set $\bar
{E}_{U_{\alpha}}=\bigcup\limits_{x\in U_{\alpha}} \bar{E_{x}}$ where $\bar{E_{x}}$ is
the closure of $E_{x}$ in $\widehat{E}_{x}$. Since $\sigma_{\alpha}\circ
\iota(x,u)=\tau_{\alpha}(x,u)$, then we can extend the trivialization
$\tau_{\alpha}:E_{U_{\alpha}}\longrightarrow U_{\alpha}\times\mathbb{E}$ to a
trivialization $\bar{\tau}_{\alpha}:\bar{E}_{U_{\alpha}}\longrightarrow
U_{\alpha}\times \bar{\mathbb{E}}$ and if $\bar{\iota}:\bar{E}_{U_{\alpha}}\longrightarrow
\widehat{E}_{U_{\alpha}}$ the canonical extension of $\iota: E_{U_{\alpha}%
}\longrightarrow\widehat{E}_{U_{\alpha}}$ then we also have $\sigma_{\alpha
}\circ\bar{\iota} =\bar{\tau}_{\alpha}$. In the same way, the transition
function $\bar{T}_{\alpha,\beta}$ associated to the bundle atlas $\{(U_{\alpha
},\bar{\tau}_{\alpha})\}_{\alpha\in A}$ is the restriction of $S_{\alpha\beta}$
to $\bar{E}$ which is an Hilbert space. The result is a consequence of Point 1.
\end{proof}

\subsection{Almost tangent, decomposable  complex  and para-complex structures}
\label{__AlmostTangentStructureAndAlmostDecomposableComplexStructure}

Let $(E,\pi_{E},M)$ be a Banach bundle such that $M$ is connected. As in
finite dimension we introduce the notions of almost tangent and almost complex structures.

\begin{definition}
\label{D_AlmostTangentStructureOnABanachBundle} 
An endomorphism $J$ of $E$ is called an almost tangent structure on $E$ if $\mathrm{im}J=\ker J$ and $\ker J$ is a supplemented sub-bundle of $E$.
\end{definition}

\begin{definition}
\label{D_AlmostComplexStructureOnABanachBundle} 
An endomorphism $\mathcal{I}$ (resp. $\mathcal{J}$ is called an almost complex (resp. paracompact) structure on $E$ if $\mathcal{I}^{2}=-\operatorname{Id}_{E}$ (resp. $\mathcal{J}^{2}=\operatorname{Id}_{E}$).
\end{definition}

According to $\S$ \ref{__TangentStructures} and $\S$ \ref{___ComplexStructures}, if
$J$ is an almost tangent structure, there exists  a decomposition $E=\ker J\oplus K$ in a
Whitney sum, the restriction $J_{K}$ of $J$ to $K$ is a bundle isomorphism
from ${K}$ to $\ker J$. Moreover, we can associate to $J$ an almost complex structure
$\mathcal{I}$ (resp. $\mathcal{J}$ on $E$ given by $\mathcal{I}(u)=-J_{K} u$ (resp. $\mathcal{J}(u)=J_{K} u$) for $u\in K$ and
$\mathcal{I}(u)={J}_{K}^{-1}(u)$ (resp. $\mathcal{I}(u)={J}_{K}^{-1}(u)$) for $u\in\ker J$.

\begin{definition}
\label{D_DecomposableParaComplexStructure} An almost  complex structure  is called
decomposable if there exists a decomposition in a Whitney sum $E=E_{1}\oplus
E_{2}$ and an isomorphism $I:E_{2}\longrightarrow E_{1}$ such that $\mathcal{I}$ can be written 
$$\begin{pmatrix} 0&-I\\
I^{-1}&0\\
\end{pmatrix}.$$
\end{definition}

Given a Withney decomposition $E=E_{1}\oplus
E_{2}$ associated to a decomposable almost complex structure $\mathcal{I}$ let  $\mathcal{S}$ be a the isomorphism of $E$ defined by the matrix
$$\begin{pmatrix} -\operatorname{Id}_{E_1}&0\\
0&\operatorname{Id}_{E_2}\\
\end{pmatrix}.$$
According to $\S$ \ref{___ParaComplexStructures}, $\mathcal{J}=\mathcal{S}\mathcal{I}$ is an almost para-complex structure. Conversely, if $\mathcal{J}$ is an almost para-complex structure there exists a decomposition in a Whitney sum $E=E_{1}\oplus
E_{2}$ and an isomorphism $I:E_{2}\longrightarrow E_{1}$ such that $\mathcal{I}$ can be written 
$$\begin{pmatrix} 0&-I\\
I^{-1}&0\\
\end{pmatrix}.$$
and in this way, $\mathcal{I}=\mathcal{S}\mathcal{J}$ is decomposable almost complex structure.\\

From Propositions \ref{P_IsomorphismDecompositionTangentStructure} and
\ref{P_IsomorphismDecomposableComplexStructure}, we obtain the following result.

\begin{theorem}
\label{T_GStructureAlmostTangentAlmostDecomposableComplex}
Let $(E,\pi_{E},M)$ be a Banach vector bundle.
\begin{enumerate}
\item[1.] 
Let $J$ be an almost tangent structure on $(E,\pi_{E},M)$.\\
For a fixed $x_{0}$ in $M$, we identify the fibre $E_{x_{0}}$ with the typical fibre of $E$. If $J_{x_{0}}$ is the induced tangent structure
on $E_{x_{0}}$, then $J$ defines a $J_{x_{0}}$-structure on $E$.

\item[2.] 
Let $\mathcal{I}$ be an almost decomposable complex structure on $(E,\pi_{E},M)$.\\
For a fixed $x_{0}$ in $M$, we identify the fibre $E_{x_{0}}$ with the typical fibre of $E$. If $\mathcal{I}_{x_{0}}$ is the induced complex structure on $E_{x_{0}}$, then $\mathcal{I}$ defines a
$\mathcal{I}_{x_{0}}$-structure on $E$.

\item[3.]
Let $\mathcal{J}$ be an almost para-complex structure on $(E,\pi_{E},M)$.\\
For a fixed $x_{0}$ in $M$, we identify the fibre $E_{x_{0}}$ with the typical fibre of $E$. If $\mathcal{J}_{x_{0}}$ is
the induced para-complex structure on $E_{x_{0}}$, then $\mathcal{J}$ defines a
$\mathcal{J}_{x_{0}}$-structure on $E$.

\end{enumerate}
\end{theorem}

\begin{proof}
Consider a decomposition $E=\ker J\oplus K$ in Whitney sum. 
Let $\{(U_{\alpha},\tau_{\alpha})\}_{\alpha\in A}$ be a bundle atlas of $E$ which
induces an atlas on each sub-bundle $\ker J$ and $K$. For any $x\in
U_{\alpha}$, $\tau_{\alpha}(x):E_{x} \longrightarrow E_{x_{0}}$ is a an
isomorphism and so $\tau_{\alpha}(x)^{*} J_{x_{0}}$ is a tangent structure on
$E_{x}$. Moreover, if $\tau^{1}_{\alpha}$ and $\tau^{2}_{\alpha}$ are the
restriction of $\tau_{\alpha}$ to $\ker J{| U_{\alpha}}$ and $K_{| U_{\alpha}%
}$ respectively, then $\tau_{\alpha}^{1}(x)$ and $\tau_{\alpha}^{2}(x)$ is an
isomorphism of $\ker J_{x}$ and $K_{x}$ onto $\ker J_{x_{0}}$ and $K_{x_{0}}$
respectively. Therefore $(\tau_{\alpha}(x)^{-1})^{*}(J_{x})$ is a (linear)
tangent structure on $E_{x_{0}}$ whose kernel is also $\ker J_{x_{0}}$ with
$K_{x_{0}}=(\tau^{2}_{\alpha}(x))^{-1}(K_{x})$ and $(\tau_{\alpha}%
(x)^{-1})^{*}J_{| K_{x}}$ is an isomorphism from $K_{x_{0}}$ to $\ker
J_{x_{0}}$. From Proposition \ref{P_IsomorphismDecompositionTangentStructure},
$T_{x}: E_{x_{0}}\longrightarrow E_{x_{0}}$ defined by $T_{x}(u,v)= (u,
(\tau_{\alpha}(x)^{-1})^{*}J_{| K_{x}}\circ J_{x_{0}}^{-1}$ is an automorphism
of $E_{x_{0}}$ such that $T_{x}^{*} (J_{x_{0}})= (\tau_{\alpha}(x)^{-1})J_{x}$
and so $(T_{x}\circ\tau_{\alpha}(x) )^{*}(J_{x_{0}})= J_{x}$. Since
$\tau_{\alpha}$ is smooth on $U$, this implies that $\tau^{\prime}_{\alpha
}(x)=T_{x}\circ\tau_{\alpha}(x)$ defines a trivialization $\tau^{\prime
}_{\alpha}:E_{U_{\alpha}} \longrightarrow U_{\alpha}\times E_{x_{0}}$, so $J_{|
U_{\alpha}}$ is locally modelled on $J_{x_{0}}$. Now it is clear that each
transition map $T_{\alpha\beta}$ belongs to the isotropy group of $J_{x_{0}}%
$.\\
The proofs of Point 2  and Point 3 use, step by step, the same type of arguments as in the
previous one and is left to the reader.
\end{proof}

\subsection{Compatible almost tangent and almost cotangent structures}
\label{_CompatibleAlmostTangentAndAlmostCotangentStructures}

Let $(E,\pi_{E},M)$ be a Banach bundle such that $M$ is connected.

\begin{definition}
\label{D_AlmostCotangentStructureOnABanachBundle}
A weak non degenerated $2$-form $\Omega$ on $E$ is called an almost cotangent structure on $E$ if there exists a a decomposition of $E$ in a Whitney sum $L\oplus K$ of Banach sub-bundles of $E$ such that each fibre $L_{x}$ is a weak Lagrangian subspace
of $\Omega_{x}$ in $E_{x}$ for all $x\in M$. In this case $L$ is called a weak
Lagrangian bundle.
\end{definition}

\begin{definition}
\label{D_CompatibleAlmostTangentAlmostCotangentStructuresOnABanachBundle} An
almost tangent structure $J$ on $E$ is compatible with an almost cotangent
structure $\Omega$ if $\Omega_{x}$ is compatible with $J_{x}$ for all $x\in M$.
\end{definition}

Assume that $E$ is a Whitney sum $E_{1}\oplus E_{2}$ of two sub-bundles of $E$
and there exists a bundle isomorphism $J:E_{2}\longrightarrow E_{1}$, then as
we have already seen, we can associate to $J$ an almost complex structure
$\mathcal{I}$ on $E$. On the other hand, if $g$ is a pseudo-Riemannian metric
on $E_{1}$, we can extend $g$ to a canonical pseudo-Riemannian metric $\bar
{g}$ on $E$ such that $E_{1}$ and $E_{2}$ are $\bar{g}$ orthogonal and the
restriction of $\bar{g}$ on $E_{2}$ is $\bar{g}(u,v)=g(Ju,Jv)$. Now, according
to $\S$ \ref{___TangentStructureCompatibleCotangentStructure}, by application of Proposition
\ref{P_CompatibilityBetweenTangentCotangentStructuresAndIndefineInnerProduct}
on each fibre, we get the following result.

\begin{proposition}\label{P_CompatibilityAlmostTangentAlmostCotangentStructuresPseudoRiemannianMetricOnABanachBundle}${}$
\begin{enumerate}
\item[1.] 
Assume that there exists an almost tangent structure $J$ on $E$ and
let $E=\ker J\oplus K$ be an associated decomposition in a Whitney sum. If
there exists a pseudo-Riemannian metric $g$ on $\ker J$, then there exists a
cotangent structure $\Omega$ on $E$ compatible with $J$ such that $\ker J$ is
a weak Lagrangian bundle of $E$ where
\[
\forall(u,v) \in E^{2}, \Omega( u,v)=\bar{g}(\mathcal{I}u,
v)-\bar{g}(u,\mathcal{I}v)
\label{eq_OmegaIg}
\]

if $\bar{g}$ is the canonical extension of $g$ and $\mathcal{I}$ is the almost
complex structure on $E$ defined by $J$.

\item[2.] 
Assume that there exists an almost cotangent structure $\Omega$ on
$E$ and let $E={L}\oplus{K}$ be an associated Whitney decomposition. Then
there exists a tangent structure $J$ on $E$ such that $\ker J={L}$

\item[3.] 
Assume that there exists a tangent structure $J$ on $E$ which is
compatible with a cotangent structure $\Omega$. Then there exists a
decomposition $E=\ker J\oplus{K}$ where $\ker J$ is a weak Lagrangian bundle.
Moreover $g(u,v)=\Omega(Ju,v)$ is a pseudo-Riemannian metric on $\ker J$ and
$\Omega$ satisfies the relation (\ref{eq_OmegaIg}).
\end{enumerate}
\end{proposition}

In the context of the previous Proposition, if $g$ is a weak Riemannian
metric, so is its extension $\bar{g}$. In this case, as we have already seen
in $\S$ \ref{KreinMetrics}, each fibre $E_{x}$ can be continuously and
densely embedded in a Hilbert space which will be denoted $\widehat
{E}_{x}$. According to Theorem \ref{T_GStructureKreinMetric}, with these
notations, we get the following theorem.

\begin{theorem}\label{T_GStructureTangentCotangentStructureWeakRiemannianMetric}${}$
\begin{enumerate}
\item[1.] 
Consider a tangent structure $J$ on $E$, a weak Riemannian metric
$g$ on $\ker J$ and $\Omega$ the cotangent structure compatible with $J$ as
defined in (\ref{eq_OmegaIg}). Assume that $\widehat{E}=\bigcup\limits_{x\in
M}\widehat{E}_{x}$ is a Banach bundle over $M$ such that the inclusion of $E$
in $\widehat{E}$ is a bundle morphism. For any $x_{0}$ in $M$, the triple
$(J,g,\Omega)$ defines a $(J_{x_{0}}, \bar{g}_{x_{0}},\Omega_{x_{0}}%
)$-structure on $E$ where $\bar{g}$ is the natural extension of $g$ to $E$.
%Moreover $\O$ is a Darboux form.

\item[2.] 
Let $(\widehat{E},\pi_{\widehat{E}},M)$ be a Hilbert bundle over a
connected Banach manifold $M$. Consider an almost tangent structure
$\widehat{J}$ on $\widehat{E}$, a Riemannian metric $\widehat{g}$ on
$\ker\widehat{J}$ and an almost cotangent structure $\widehat{\Omega}$
compatible with $\widehat{J}$ which satisfies (\ref{eq_OmegaIg}) on
$\widehat{E}$. Then for any Banach bundle $(E,\pi_{E},M)$ such that there
exists an injective bundle morphism $\iota:E\longrightarrow\widehat{E}$ with
dense range, the restriction $J$ of $\widehat{J}$ to $E$ is an almost tangent
structure on $E$ and the restriction $\Omega$ of $\widehat{\Omega}$ to $E$ is
an almost cotangent structure compatible with $J$ and $\Omega$ satisfies the
relation (\ref{eq_OmegaIg}) on $E$. In particular, for any $x_{0}$ in $M$, the
triple $({J},{g},{\Omega})$) defines a $(J_{x_{0}}, {g}_{x_{0}},\Omega_{x_{0}%
})$-structure on $E$.
\end{enumerate}
\end{theorem}

\begin{remark}
\label{R_TangentStructureCompatibleWithCotangentStructure} According to
Proposition
\ref{P_CompatibilityAlmostTangentAlmostCotangentStructuresPseudoRiemannianMetricOnABanachBundle},
Point 3, if there exists a tangent structure $J$ on $E$ which is compatible with
a cotangent structure $\Omega$ such that $g(u,v)=\Omega(Ju,v)$ is a Riemannian
metric on $\ker J$, then $\Omega$ satisfies the relation (\ref{eq_OmegaIg}).
Therefore we have a corresponding version of Theorem \ref{eq_OmegaIg} with the
previous assumption.
\end{remark}

\begin{proof}
1.
%Let  $\widehat{g}$ be the canonical extension of $g$ to $E$.
According to our assumption, by density of $\iota(E)$ (identified with $E$) in
$\widehat{E}$ and using compatible bundle atlases for $\widehat{E}$ and $E$
(\textit{cf.} proof of Theorem \ref{T_GStructureKreinMetric}) it is easy to
see that we can extend $J$ and ${g}$, to an almost tangent structure
$\widehat{J}$ and a strong Riemannian metric $\widehat{g}$, on $\widehat{E}$.
By the way, if $\mathcal{I}$ is the almost complex structure on $E$ associated
to $J$, we can also extend $\mathcal{I}$ to an almost complex structure
$\widehat{\mathcal{I}}$ on $\widehat{E}$ which is exactly the almost complex
structure associated to $\widehat{J}$. 
Thus this implies, by relation (\ref{eq_OmegaIg}) that we can extend $\Omega $ into $2$-form $\widehat{\Omega}$ on $\widehat{E}$ which will satisfy also relation (\ref{eq_OmegaIg}) on
$\widehat{E}$, using $\widehat{g}$ and $\widehat{\mathcal{I}}$. Then
$\widehat{\Omega}$ is a strong non-degenerate $2$-form on $\widehat{E}$. Note
that $\ker\widehat{J}$ is the closure in $\widehat{E}$ of $\ker J$ an since
$J: K\longrightarrow\ker J$ is an isomorphism, this implies that $\widehat
{K}=\widehat{J}(\ker\widehat{J})$ is a Banach sub-bundle of $\widehat{E}$ such
that $\widehat{E}=\ker\widehat{J}\oplus\widehat{K}$. Moreover, each fibre
$\ker\widehat{J}_{x}$ of $\ker J$ is exactly the Hilbert space which is the
closure of $\ker J_{x}$ provided with the pre-Hilbert product $g_{x}$. 
Now, from the construction of the extension $\bar{g}$ of $g$ to $E$, the
restriction of $J$ to $K$ is an isometry and so the restriction of
$\widehat{J}$ to $\widehat{K}$ is also an isometry. This implies that
$\mathcal{I}$ (resp. $\widehat{\mathcal{I}}$) is an isometry of $\bar{g}$
(resp. $\widehat{g}$). 
Fix a point $x_{0}$ in $M$ and identify the fibres of $E$
and $\widehat{E}$ with $E_{x_{0}}$ and $\widehat{E}_{x_{0}}$. From Theorem
\ref{T_GStructureKreinMetric} the exists an atlas $\{(U_{\alpha},\tau_{\alpha
})\}_{\alpha\in A}$ of $E$ and $\{(U_{\alpha},\sigma_{\alpha})\}_{\alpha\in A}$
such that $\tau_{\alpha}=\sigma\alpha\circ\iota$ and the respective transition
functions $T_{\alpha\beta}$ and $S_{\alpha\beta}$ are isometries of $E$ relatively to $\bar
{g}_{x_{0}}$ and $\widehat{g}_{x_{0}}$. Since $\mathcal{I}$ and $\widehat
{\mathcal{I}}$ are isometries for $\bar{g}$ and $\widehat{g}$, this implies
that $\mathcal{I}$ is also $\mathcal{I}_{x_{0}}$-structure on $E$ and
$\widehat{\mathcal{I}}$ is a $\widehat{\mathcal{I}}_{x_{0}}$-structure on
$\widehat{E}$. But since $J$ and $\widehat{J}$ are canonically defined from
$\mathcal{I}$ and $\widehat{\mathcal{I}}$ respectively, we obtain a similar
result for $J$ and $\widehat{J}$. Finally, according to relation
(\ref{eq_OmegaIg}), we also obtain a similar result for $\Omega$ and
$\widehat{\Omega}$ which ends the proof.\newline

2. Under the assumptions of Point 2, we can apply all the results of Point 1
to $(\widehat{E},\pi_{\widehat{E}},M)$. In particular, for a decomposition
$\widehat{E}=\ker\widehat{J}\oplus\widehat{K}$, $\widehat{J}$ is an isometry
from $\widehat{K}$ to $\ker\widehat{J}$ and the almost complex structure $\widehat
{\mathcal{I}}$ associated to $\widehat{J}$ is an isometry of $\widehat{g}$ .
We can identify $E$ with $\iota(E)$ in $\widehat{E}$. We set $K=E\cap
\widehat{K}$. Then the range of the restriction of $\widehat{J}$ to $K$ is
included in $\ker\widehat{J}$. Using compatible bundle atlases for $E$ and
$\widehat{E}$ we can show that if $J$ is the restriction of $\widehat{J}$ then
$\ker J= \ker\widehat{J}\cap E$ and $E=\ker J\oplus K$ so that $J$ is an
almost tangent structure on $E$. Now $\widehat{g}$ induces a weak Riemannian
metric $\bar{g}$ on $E$ and $J$ is again an isometry from $K$ to $\ker J$.
This implies that the almost complex structure $\mathcal{I}$ associated to $J$
is nothing but else than the restriction of $\widehat{\mathcal{I}}$ to $E$. Of
course we obtain that the restriction $\Omega $ of $\widehat{\Omega}$ satisfies
the relation (\ref{eq_OmegaIg}) relative to $\bar{g}$ and $\mathcal{I}$. So
the assumptions of Point 1 are satisfied which ends the proof.
\end{proof}

\subsection{Compatible weak symplectic form, weak Riemannian metric and almost complex structures}
\label{__CompatibleWeakSymplecticFormWeakRiemannianMetricAlmostComplexStructuresOnABanachBundle}%

If $g$ (resp. $\Omega$) is a weak neutral metric (resp. a weak symplectic
form) on a Banach bundle $(E,\pi_{E},M)$, as in the linear context, we denote
by $g^{\flat}$ (resp. $\Omega^{\flat}$ the associated morphism from $E$ to
$E^{*}$ where $(E^{*},\pi_{E^{*}},M)$ is the dual bundle of $(E,\pi_{E},M)$. Moreover if $\mathcal{I}$ is an almost complex structure on $E$,  following $\S$ \ref{__CompatibilityBetwenDifferentStructures}, we introduce the following notions.

\begin{definition}\label{D_CompatibilityWeakSymplecticFormAlmostComplexStructureWeakRiemannianMetricOnABanachBundle}${}$
\begin{enumerate}
\item[1.] 
We say that a weak symplectic form $\Omega$ and an almost complex
structure $\mathcal{I}$ on $E$ are compatible if $(u,v)\mapsto\Omega
(u,\mathcal{I}v)$ is a weak Riemannian metric on $E$ and $\Omega
(\mathcal{I}u,\mathcal{I}v)=\Omega(u,v)$ for all $u$ and $v$ in $E$.

\item[2.] 
We say that a weak Riemannian metric $g$ and an almost complex
structure $\mathcal{I}$ on $E$ are compatible if $g(\mathcal{I}%
u,\mathcal{I}v)=g(u,v)$ for all $u$ and $v$ in $E$.

\item[3.] 
We say that a weak Riemannian metric $g$ and a weak symplectic
structure $\Omega$ on $E$ are compatible, if $\mathcal{I}=(g^{\flat}%
)^{-1}\circ\Omega^{\flat}$ is well defined and is a complex structure on
$E$.
\item[4.] A weak symplectic form $\Omega$ on $E$ will be called a Darboux form if there exists a decomposition $E=E_1\oplus E_2$ such that for each $x\in M$ each fiber $(E_1)_x$ and $(E_2)_x$ is Lagrangian. 
\end{enumerate}
\end{definition}

Now by application of Proposition \ref{P_CompatibilitySymplecticStructureInnerProductDecomposableComplexStructures}
and Corollary \ref{C_DarbouxFormPreHilbertProductComplexStructure} we obtain:

\begin{theorem}
\label{T_CompatibilityDarbouxFormWeakRiemannianMetricComplexStructureOnABanachBundle}
Consider a Darboux form $\Omega$, a weak Riemannian $g$ and a decomposable complex
structure $\mathcal{I}$ on a Banach space $\mathbb{E}$. Assume that any pair
among such a triple exists on $E$ and is compatible. Then the third one also
exists and is compatible with any element of the given pair. Denote
$\widehat{E}_{x}$ the Hilbert space defined by $g_{x}$ for each $x\in M$. If
$\widehat{E}=\bigcup\limits_{x\in M}\widehat{E}_{x}$ is a Banach bundle over $M$ and
there exists an injective morphism $\iota: E\longrightarrow\widehat{E}$ then
we can extend $\Omega$, $g$ and $\mathcal{I}$ to a strong symplectic form
$\widehat{\Omega}$, a strong Riemannian metric $\widehat{g}$, an almost
complex structure $\widehat{\mathcal{I}}$ on $\widehat{E}$ respectively.
Moreover, for any $x_{0}\in M$, if we identify the typical fibre of $E$ and
of $\widehat{E}$ with $E_{x_{0}}$ and $\widehat{E}_{x_{0}}$ respectively, then
the triple $(\Omega,g,\mathcal{I})$ (resp. $(\widehat{\Omega},\widehat
{g},\widehat{\mathcal{I}})$) defines a $(\Omega_{x_{0}},g_{x_{0}}%
,\mathcal{I}_{x_{0}})$-structure (resp. $(\widehat{\Omega}_{x_{0}}\widehat
{g}_{x_{0}},\widehat{\mathcal{I}}_{x_{0}})$-structure) on $E$ (resp. $\widehat{E})$).\\
Conversely, assume that we have such a triple $(\widehat{\Omega},\widehat
{g},\widehat{\mathcal{I}})$ which is compatible on a Hilbert bundle
$(\widehat{E},\pi_{\widehat{E}},M)$. For any Banach bundle which can be
continuously and densely embedded in $\widehat{E}$, then this triple induces,
by restriction to $E$, a triple $(\Omega,g,\mathcal{I})$ which is compatible.
Therefore, we obtain a $(\Omega_{x_{0}},g_{x_{0}},\mathcal{I}_{x_{0}}%
)$-structure (resp. $(\widehat{\Omega}_{x_{0}},\widehat{g}_{x_{0}}%
,\widehat{\mathcal{I}}_{x_{0}})$-structure) on $E$ (resp. $\widehat{E})$).
\end{theorem}

\begin{proof} The first property is a direct application of Proposition
\ref{P_CompatibilitySymplecticStructureInnerProductDecomposableComplexStructures}
and Corollary \ref{C_DarbouxFormPreHilbertProductComplexStructure}, the other ones are obtained as in the proof of the corresponding parts of Theorem
\ref{T_GStructureTangentCotangentStructureWeakRiemannianMetric}.
\end{proof}

\begin{definition}
\label{D_KahlerBanachBundle} 
A Banach bundle $\pi_E:E\longrightarrow M$  has a weak (resp. strong) almost  K\"ahler structure if there exists on $E$ a weak (resp. strong) Darboux form $\Omega$, a weak (resp strong) Riemanian metric $g$ and a decomposable almost complex structure $\mathcal{I}$ such that there exists a compatible  pair among these three data.\\
\end{definition}
From Remark \ref{R_DecompositionInLagrangianOrthogonalspaces}, to a weak or strong almost  K\"ahler structure on $E$  is associated a decomposition $E=E_1\oplus E_2$ of isomorphic sub-bundle which are lagrangian and orthogonal.
Note that  Theorem \ref{T_CompatibilityDarbouxFormWeakRiemannianMetricComplexStructureOnABanachBundle} can be seen as sufficient conditions under which a weak K\"ahler structure on a Banach bundle $E$  is a tensor structure on $E$.  \\
 
\subsection{Compatible weak symplectic forms, weak neutral metrics and almost para-complex structures }
\label{__CompatibleWeakSymplecticFormsWeakNeutralMetricsAlmostParaComplexStructuresOnABanachBundle}%

If $g$ (resp. $\Omega$) is a weak Riemannian metric (resp. weak symplectic
form) on a Banach bundle $(E,\pi_{E},M)$, as in the linear context, we denote
by $g^{\flat}$ (resp. $\Omega^{\flat}$ the associated morphism from $E$ to
$E^{*}$ where $(E^{*},\pi_{E^{*}},M)$ is the dual bundle of $(E,\pi_{E},M)$. Moreover, if $\mathcal{J}$ is an almost para-complex structure on $E$,  following $\S$ \ref{__CompatibilityBetwenDifferentStructures}, we introduce the following notions.

\begin{definition}
\label{D_CompatibilityWeakSymplecticFormAlmostParaComplexStructureWeakNeutralMetricOnABanachBundle}${}$
\begin{enumerate}
\item[1.] 
We say that a weak symplectic form $\Omega$ and an almost para-complex
structure $\mathcal{J}$ on $E$ are compatible if $(u,v)\mapsto\Omega
(u,\mathcal{J}v)$ is a weak neutral metric on $E$ and $\Omega
(\mathcal{J}u,\mathcal{J}v)=-\Omega(u,v)$ for all $u$ and $v$ in $E$.

\item[2.] 
We say that a weak neutral metric $g$ and an almost para-complex
structure $\mathcal{J}$ on $E$ are compatible if $g(\mathcal{J}%
u,\mathcal{J}v)=-g(u,v)$ for all $u$ and $v$ in $E$.

\item[3.] 
We say that a weak neutral metric $g$ and a weak symplectic
structure $\Omega$ on $E$ are compatible, if $\mathcal{J}=(g^{\flat}%
)^{-1}\circ\Omega^{\flat}$ is well defined and is an almost  para-complex  structure on
$E$.
%\item[4.] A weak symplectic form $\Omega$ on $E$ will be called a Darboux form if there exists a decomposition $E=E_1\oplus E_2$ such that for each $x\in M$ each fiber $(E_1)_x$ and $(E_2)_x$ is Lagrangian. 
\end{enumerate}
\end{definition}

According to Remark \ref{R_AlmostDecomposableComplexAlmostParaComplex} and Theorem \ref{T_CompatibilityDarbouxFormWeakRiemannianMetricComplexStructureOnABanachBundle} we have 

\begin{theorem}
\label{T_CompatibilityDarbouxFormsWeakNeutralMetricParaComplexStructuresOnABanachBundle}
Consider a Darboux form $\Omega$, a weak neutral metric  $g$ and a para-complex
structure $\mathcal{J}$ on a Banach space $\mathbb{E}$. Assume that any pair
among such a triple exists on $E$ and is compatible. Then the third one also
exists and is compatible with any element of the given pair. Let $h$ be a Riemaniann metric canonically associated to some decomposition $E=E^+\oplus E^-$ relative to $g$ and let
$\widehat{E}_{x}$ be the Hilbert space defined by the inner product associated to $h_x$ for each $x\in M$. If
$\widehat{E}=\bigcup\limits_{x\in M}\widehat{E}_{x}$ is a Banach bundle over $M$ and
there exists an injective morphism $\iota: E\longrightarrow\widehat{E}$ then
we can extend $\Omega$, $g$ and $\mathcal{J}$ to a strong symplectic form
$\widehat{\Omega}$, a strong neutral metric $\widehat{g}$ and an almost
para-complex structure $\widehat{\mathcal{J}}$ on $\widehat{E}$ respectively.
Moreover, for any $x_{0}\in M$, if we identify the typical fibre of $E$ and
of $\widehat{E}$ with $E_{x_{0}}$ and $\widehat{E}_{x_{0}}$ respectively, then
the triple $(\Omega,g,\mathcal{J})$ (resp. $(\widehat{\Omega},\widehat
{g},\widehat{\mathcal{J}})$) defines a $(\Omega_{x_{0}},g_{x_{0}}%
,\mathcal{J}_{x_{0}})$-structure (resp. $(\widehat{\Omega}_{x_{0}}\widehat
{g}_{x_{0}},\widehat{\mathcal{J}}_{x_{0}})$-structure) on $E$ (resp. $\widehat{E})$.
Conversely, assume that we have such triple $(\widehat{\Omega},\widehat
{g},\widehat{\mathcal{J}})$ which are compatible on a Hilbert bundle
$(\widehat{E},\pi_{\widehat{E}},M)$. For any Banach bundle which can be
continuously and densely embedded in $\widehat{E}$, then this triple induces,
by restriction to $E$, a triple $(\Omega,g,\mathcal{J})$ which is compatible.
Therefore, we obtain a $(\Omega_{x_{0}},g_{x_{0}},\mathcal{J}_{x_{0}}%
)$-structure (resp. $(\widehat{\Omega}_{x_{0}},\widehat{g}_{x_{0}}%
,\widehat{\mathcal{J}}_{x_{0}})$-structure) on $E$ (resp. $\widehat{E})$.
\end{theorem}

As in finite dimension (cf. \cite{Lib} for instance) we introduce the para-K\"ahler structures.
 
\begin{definition}
\label{D_ParaKahlerBanachBundle} 
A Banach bundle $\pi_E:E\longrightarrow M$ 
has a weak (resp. strong)  almost para-K\"ahler structure if there exists on $E$ a weak (resp. strong) Darboux form $\Omega$, a weak (resp strong) neutral metric $g$, a  almost para-complex 
structure $\mathcal{J}$ such that  a compatible  pair among these three data.\\
\end{definition}

As for weak or strong  almost  para-K\"ahler bundles, from Remark \ref{R_DecompositionInLagrangianOrthogonalspacesParaComplex}, to a weak or strong para-K\"ahler structure on $E$  is associated a decomposition $E=E_1\oplus E_2$ of isomorphic sub-bundle which are lagrangian  and the restriction of ${g}$ to ${E}_1$ (resp. ${E}_2$) is positive definite (resp. negative definite).\\
Note also  that  Theorem \ref{T_CompatibilityDarbouxFormsWeakNeutralMetricParaComplexStructuresOnABanachBundle} can be seen as sufficient conditions under which a weak almost para-K\"ahler structure on a Banach bundle $E$  is a tensor structure on $E$.  \\

\section{Examples of integrable tensor structures on a Banach manifold}
\label{_ExamplesOfIntegrableTensorStructuresOnABanachManifold}

\subsection{The Darboux Theorem on a Banach manifold}
\label{__DarbouxTheoremOnABanachManifold}

In finite dimension, from the Darboux theorem, a symplectic form on a manifold
defines an integrable tensor structure on $M$. The extension of such a result
to the Banach framework is given in \cite{Bam} for weak symplectic Banach
manifolds.\newline

\begin{definition}
A weak symplectic form on a Banach manifold $M$ modelled on a Banach space
$\mathbb{M}$ is a closed $2$-form $\Omega$ on $M$ which is non-degenerate.
\end{definition}

If $\Omega^{\flat}:TM\rightarrow T^{*}M$ is the associated morphism, then
$\Omega^{\flat}$ is an injective bundle morphism. The symplectic form $\Omega$
is weak if $\Omega^{\flat}$ is not surjective. Assume that $\mathbb{M}$ is
reflexive. We denote by $\widehat{T_{x}M}$ the Banach space which is the
completion of $T_{x}M$ provided with the norm $||\;||_{\Omega_{x}}$. Recall
that $\widehat{T_{x}M}$ does not depend on the choice of the norm on $T_{x}M$
(cf. Remark \ref{R_IndependanceOfTheNorms}). Then $\Omega$ can be
extended to a continuous bilinear map $\widehat{\Omega}_{x}$ on $T_{x}%
M\times\widehat{T_{x}M}$ and $\Omega^{\flat}$ becomes an isomorphism from
$T_{x}M$ to $(\widehat{T_{x}M})^{*}$. We set
\[
\widehat{TM}=\displaystyle\bigcup_{ x\in M}\widehat{T_{x}M} \;\text{ and }\;
(\widehat{TM})^{*}=\displaystyle\bigcup_{ x\in M}(\widehat{T_{x}M})^{*}.
\]

We have the following Darboux Theorem (\cite{Bam}) (see also \cite{Pel2}):

\begin{theorem}
[Local Darboux theorem] 
\label{T_LocalDarboux} 
Let $\Omega$ be a weak symplectic form on a Banach manifold modelled on a reflexive Banach space
$\mathbb{M}$. Assume that we have the following assumptions:

\begin{enumerate}
\item[(i)] 
There exists a neighbourhood $U$ of $x_{0}\in M$ such that
$\widehat{TM}_{| U}$ is a trivial Banach bundle whose typical fiber is the
Banach space $(\widehat{T_{x_{0}}M}, ||\;||_{\Omega_{x_{0}}})$;
%\item[(ii)] the map $||\;||_\omega$ defined by $||(x,u)||_\omega:=||u||_{\omega_x}$ is a Finsler metric on  $\widehat{TM}_{| U}$;

\item[(ii)]
$\Omega$ can be extended to a smooth field of continuous
bilinear forms on \newline$TM_{| U}\times\widehat{TM}_{| U}$.
\end{enumerate}

Then there exists a chart $(V, F)$ around $x_{0}$ such that $F^{*}\Omega
_{0}=\Omega$ where $\Omega_{0}$ is the constant form on $F(V)$ defined by
$(F^{-1})^{*}\Omega_{x_{0}}$.\newline
\end{theorem}

\begin{definition}
\label{D_DarbouxChart} 
The chart $(V,F)$ in Theorem \ref{T_LocalDarboux} will be called a Darboux chart around $x_{0}$.
\end{definition}

\begin{remark}
\label{R_CompareBambusi} 
The assumptions of Darboux theorem in \cite{Bam} (Theorem 2.1) are formulated in a different way. In fact, the assumption on
all those norms $||\;||_{\Omega_{x}}$ in this Theorem 2.1 on the typical fiber
$\widehat{\mathbb{M}}$ is a consequence of assumptions (i) and (ii) of Theorem
\ref{T_LocalDarboux} after shrinking $U$ if necessary (cf. \cite{Pel2}).
\end{remark}

\begin{remark}
\label{R_DarbouxTheoremStrongSymplecticForm} 
If $\Omega$ is a strong symplectic form on $M$, then $\mathbb{M}$ is reflexive and $\Omega^{\flat}$ is
a bundle isomorphism from $TM$ to $T^{*}M$. In particular, the norm
$||\;||_{\Omega_{x}}$ is equivalent to any norm $||\;||$ on $T_{x}M$ which
defines its Banach structure and so all the assumptions (i) and (ii) of
Theorem \ref{T_LocalDarboux} are always locally satisfied. Thus Theorem
\ref{T_LocalDarboux} recovers the Darboux Theorem which is proved in
\cite{Mar} or \cite{Wei}.
\end{remark}

Weinstein gives an example of a weak symplectic form $\Omega$ on a
neighbourhood of $0$ of a Hilbert space $\mathbb{H}$ for which the Darboux
Theorem is not true. The essential reason is that the operator $\Omega^{\flat
}$ is an isomorphism from $T_{x}U$ onto $T_{x}^{*}U$ on $U\setminus
\left\lbrace 0\right\rbrace $, but $\Omega^{\flat}_{0} $ is not surjective. 

%In this way, we have the following Corollary:
%
%\begin{corollary}
%\label{C_EquivalenceDarboux} 
%Let $\Omega$ be a weak symplectic form on a Banach manifold modelled on a reflexive Banach space $\mathbb{M}$. Then there exists a Darboux chart $(V, F)$ around $x_{0}$ if and only
%
%\begin{enumerate}
%\item[(i)] 
%There exists a chart $(U,\phi)$ around $x_{0}$ such that
%$\widehat{TM}_{| U}$ is a trivial Banach bundle whose typical fiber
%is the completion $\widehat{\mathbb{M}}$ of the normed space $(\mathbb{M},||\;||_{(\phi^{-1})^{*}\Omega_{x_{0}}})$ and $T\phi$ can be extended to a trivialization $\widehat{T\phi}: \widehat{TM}_{| U}$ on $U\times \widehat{\mathbb{M}}$;
%
%\item[(ii)] 
%$\Omega$ can be extended to a smooth field of continuous bilinear form on $TM_{| U}\times\widehat{TM}_{| U}$.
%\end{enumerate}
%\end{corollary}

Finally from Theorem \ref{T_LocalDarboux}, we also obtain the following global
version of a Darboux Theorem:

\begin{theorem}
[Global Darboux Theorem] 
\label{T_GlobalDarboux} 
Let $\Omega$ be a weak symplectic form on a Banach manifold modelled on a reflexive Banach space
$\mathbb{M}$. Assume that we have the following assumptions:

\begin{enumerate}
\item[(i)] $\widehat{TM}\rightarrow M$ is a Banach bundle whose typical fiber
is $\widehat{\mathbb{M}}$;
%and the inclusion $j:TM\longrightarrow \widehat{TM}$ is a  bundle morphism;
%\item[(ii)] the map $||\;||_\omega$ defined by $||(x,u)||_\omega:=||u||_{\omega_x}$ is a Finsler metric on  $\widehat{TM}$;

\item[(ii)] $\Omega$ can be extended to a smooth field of continuous bilinear
forms on \newline$TM\times\widehat{TM}_{| U}$.
\end{enumerate}

Then, for any $x_{0}\in M$, there exists a Darboux chart $(V, F)$ around
$x_{0}$.\newline 
In particular $\Omega$ defines an integrable tensor structure on $M$ if and only if the assumptions (i) and (ii) are satisfied.
\end{theorem}

Note that, from Remark \ref{R_DarbouxTheoremStrongSymplecticForm}, a strong
symplectic form on a Banach manifold is always an integrable tensor structure.

\subsection{Flat pseudo-Riemannian metrics on a Banach manifold}
\label{__FlatPseudeRiemannianMetricsOnABanachManifold}

In finite dimension, a pseudo-riemannian metric on a manifold $M$ defines an
integrable tensor structure on $M$. We give a generalization of this result to
the Banach framework.

\begin{definition}
\label{D_PseudoRiemannianMetricKreinMetricOnABanchManifold}
A pseudo-Riemannian metric (resp. a Krein metric) $g$ on a Banach manifold $M$
is a pseudo-Riemannian metric (resp. a Krein metric) on the tangent bundle $TM$.
\end{definition}

When $g$ is a Krein metric, we have a decomposition $TM=TM^{+}\oplus TM^{-}$
in a Whitney sum such that the restriction $g^{+}$ (resp. $-g^{-}$) of $g$
(resp. $-g$). to $TM^{+}$ (resp. $TM^{-}$) is a (weak) Riemannian metric. To
$g$ is associated a canonical Riemannian metric $\gamma=g^{+}-g^{-}$.
According to Appendix B, § \ref{_LinearOperators}, on each fibre $T_{x}M$, we have
a norm $||\;||_{g_{x}}$ which is associated to the inner product ${\gamma}%
_{x}$; We denote by $\widehat{T_{x}M}$ Banach Hilbert space associated to the
normed space ($T_{x}M,||\;||_{g_{x}})$. \newline

By application of Theorem \ref{T_GStructureKreinMetric} to $TM$ we have the following theorem.

\begin{theorem}${}$
\label{T_ExtensionToAStrongKreinMetric}
\begin{enumerate}
\item[1.] Let $g$ be a Krein metric on a Banach $M$. Consider a decomposition
$TM=TM^{+}\oplus TM^{-}$ in a Whitney sum such that the restriction $g^{+}$
(resp. $-g^{-}$) of $g$ (resp. $-g$) to $TM^{+}$ (resp. $TM^{-}$) is a (weak)
Riemannian metric. Moreover, assume that $\widehat{TM}=\bigcup\limits_{x\in
M}\widehat{T_{x}M}$ is a Banach bundle over $M$ such that the inclusion of $TM$
in $\widehat{TM}$ is a bundle morphism. 
Then $g$ can be extended to a strong Krein metric $\widehat{g}$ on the bundle $\widehat{TM}$ and we have a decomposition $\widehat{TM}=\widehat{TM}^{+}\oplus\widehat{TM}^{-}$ such that the
restriction $\widehat{g}^{+}$ (resp. $\widehat{g}^{-}$) of $\widehat{g}$
(resp. $-\widehat{g}$) to $\widehat{TM}^{+}$ (resp. $\widehat{TM}^{-}$) is a
strong Riemannian metric. In fact $\widehat{TM}^{+}$ (resp. $\widehat{TM}^{-}%
$) is the closure of $TM^{+}$ (resp. $TM^{+}$) in $\widehat{TM}$.
%Moreover, the Levi-Civita connection $\nabla$ of $g$ is well defined.

\item[2.] 
Let $(\widehat{TM},\widehat{\pi},M)$ be a Banach bundle whose
typical fibre is a reflexive Banach space and let $\widehat{g}$  be a strong Krein metric on $\widehat{TM}$.
%We denote by $\widehat{\g}$ the strong Riemannian metric associate to $\widehat{g}$.
%on $\widehat{E}_{x_0}$ then $\widehat{g}$ is a $\widehat{g}_{x_0}$-structure  on $\widehat{E}$. \\
Assume that there exists an injective morphism of Banach bundles $\iota:
TM\longrightarrow\widehat{TM}$ whose range is dense. Then the restriction
$g=\iota^{*}\widehat{g}$ of $\widehat{g}$ is a Krein metric on $M$.
%and the Levi-Civita connection of $g$ is well defined
\end{enumerate}
\end{theorem}

Now, as for a strong Riemannian metric on a Banach manifold, to a strong
pseudo-Riemannian metric $g$ is associated a Levi-Civita connection which is a
Koszul connection $\nabla$ characterized by

\label{eq_LeviCivita}
$
\begin{array}{cc}
2g(\nabla_{X}Y,Z)= &  X(g(Y,Z))+Y(g(Z,X))-Z(g(X,Y)) \\ 
 & +g([X,Y],Z)-g([Y,Z],X)-g([X,Z],Y)
\end{array} 
$

for all (local) vector fields $X,Y,Z$ on $M$. \\
When $g$ is a weak pseudo-Riemannian metric, a Levi-Civita connection need not exist,
but if it exists, this connection is unique. \\
For a weak Riemannian metric, there exist many examples of weak Riemannian metrics for which its Levi-Civita connection is defined. When the model $\mathbb{M}$ of $M$ is a
\textit{reflexive} Banach space we will give sufficient conditions for
a weak pseudo-Riemannian metric under which the Levi-Civita connection
exists. Before, we need to introduce some preliminaries.\newline

According to § \ref{__OperatorsFromBanachSpaceToItsDual}, for each $x\in M$,
from a given norm $||\;||$ on $T_{x}M$, we can define a norm $||\;||_{g_{x}}$
on $T_{x}M$ and the completion $T_{x}M_{g}$ of the normed space $(T_{x}%
M,||\;||_{g_{x}})$ does not depend on the choice of the norm $||\;||$ (see
Remark \ref{R_IndependanceOfTheNorms}). 
In this way, $g_{x}^{\flat}$ can be extended to an isomorphism from $T_{x}M_{g}$ to $T_{x}^{*}M$ and
$g_{x}^{\flat}$ becomes an isomorphism from $T_{x}M$ to $(T_{x}M_{g})^{*}$.
Now we set:
\[
TM_{g}=\bigcup_{x\in M}T_{x}M_{g}\;\;\text{ and }\;\; T^{*}M_{g}=\bigcup_{x\in
M}(T_{x}M_{g})^{*}.
\]
Note that $TM$ is a subset of $TM_{g}$ and each fibre $T_{x}M$ is dense in
$T_{x}M_{g}$ and $g_{x}^{\flat}$ can be extended to an isometry from
$T_{x}M_{g}$ to $T_{x}^{*}M$. With these notations, we have:

\begin{proposition}
\label{P_LeviCivitaConnectionWellDefined} 
Let $g$ be a pseudo-Riemannian metric on a Banach manifold modelled on a reflexive Banach space. Assume that
the following assumptions are satisfied:

\begin{enumerate}
\item[(i)] 
$TM_{g}$ has a structure of Banach bundle over $M$;

\item[(ii)] 
The injective morphism $g^{\flat}:TM\longrightarrow T^{*}M$ can be
extended to a bundle morphism from $TM_{g}$ to $T^{*}M$.
\end{enumerate}

Then the Levi-Civita connection of $g$ is well defined.
\end{proposition}

\begin{proof}
Fix some $x_{0}\in M$ and, for the sake of simplicity, denote by $\mathbb{E}$
the Banach space $T_{x_{0}}M_{g}$. From assumption (i), it follows that around
$x_{0}$, there exists a trivialization $\Psi:(TM_{g})_{| U}\longrightarrow
U\times\mathbb{M}_{g}$ where $\mathbb{M}_{g}$ denotes the typical fibre of
$TM_{g}$. After shrinking $U$, if necessary, we may assume that $U$ is a
chart domain associated to some $\phi:U\longrightarrow\mathbb{M}$ and so
$T\pi:TM_{| U}\longrightarrow\phi(U)\times\mathbb{M}$ is a trivialization. We
may assume that $\Psi: (TM_{g})_{| U}\longrightarrow\phi(U)\times\mathbb{M}$ is a trivialization and so, without loss of generality, we can also assume that $U$ is an open subset of $\mathbb{M}$ and $TM_{|U}=U\times\mathbb{M}$. Now since $\mathbb{M}$ is dense in $\mathbb{M}_{g}$, the
inclusion of the trivial bundle $U\times\mathbb{M}$ into $U\times\mathbb{E}$ is
an injective bundle morphism.\newline

For fixed local vector fields $X$ an $Y$ defined on $U$, the map which, to any
local vector field $Z$ defined on $U$, associates the second member of
(\ref{eq_LeviCivita}) is a well defined local $1$-form on $U$ denoted
$\alpha_{X,Y}$. \\ 
Note that from Proposition \ref{P_ExtensionLinearMappingLinearSpaceIntoDual}, 3, since $\mathbb{M}$ is reflexive, $g_{x}^{\flat}$ can be extended to an isomorphism from $T_{x}M$ to
$(T_{x}M_{g})^{*}$ and also gives rise to an isomorphism from $T_{x}M$ to
$(T_{x}M_{g})^{*}$. Thus, from assumption (ii), the extension of $g^{\flat}$
gives rise to a bundle isomorphism from $TM$ to $T^{*}M_{g}$ again denoted
$g^{\flat}$. But it is clear that any $1$-form $\alpha$ on $U\times\mathbb{M}$
can be extended to a smooth $1$-form (again denoted $\alpha$) on $(TM_{g})_{|
U}=U\times\mathbb{M}$. This implies that $(g^{\flat})^{-1}(\alpha)$ is a
smooth vector fields on $U$, which ends the proof.

\end{proof}

\begin{theorem}
\label{T_gIntegrableTensorStructureNullCurvature} 
Let $g$ be a pseudo-Riemannian metric on a Banach manifold. Assume that the Levi-Civita
connection $\nabla$ of $g$ is defined. Then $g$ defines an integrable tensor
structure if and only if the curvature of $\nabla$ vanishes.
\end{theorem}

This result is well known but we have no precise and accessible reference to a
complete proof of such a result. We give a sketch of the proof.

\begin{proof}
On the principal frame bundle $\ell(TM)$, the Levi-Civita gives rise to a
connection form $\omega$ with values in the Lie algebra $gl(\mathbb{M})$ of
the Banach Lie group $\operatorname{GL}(\mathbb{M})$ (cf.\cite{KrMi} 8). Let
$\Omega$ be the curvature of $\omega$. If $X$ and $Y$ are local vector fields
on $M$, let $X^{h}$ and $Y^{h}$ their horizontal lifts in $\ell(TM)$ and we
have $\Omega(X^{h},Y^{h})=-\omega([X^{h},Y^{h}])$. Therefore, the horizontal
bundle is integrable if the curvature vanishes. 
Moreover, in this case, over any simply connected open set $U$ of $M$, the bundle $\ell(TM)_{|U}$ is
trivial (consequence of \cite{KrMi}, Theorem 39.2 for instance). Thus, around
any point $x_{0}\in M$, there exists a chart domain $(U,\phi)$ such that
$\ell(TM)_{|U}\equiv\phi(U)\times\operatorname{GL}(\mathbb{M})$. Without loss
of generality, we may assume that $U=\mathbb{M}$, $x_{0}=0\in\mathbb{M}$ and
so $\ell(TM)_{|U}= \mathbb{M}\times\operatorname{GL}(\mathbb{M})$. The
horizontal leaves are obtained from $\mathbb{M}\times \operatorname{Id}_{\mathbb{M}}$ by
right translation in $\operatorname{GL}(\mathbb{M})$. Consider an horizontal
section $\psi: \mathbb{M}\longrightarrow\operatorname{GL}(\mathbb{M})$ such
that $\psi(0)=\operatorname{Id}_{\mathbb{M}}$ then $\psi(x)=\operatorname{Id}_{\mathbb{M}}$ for all
$x\in\mathbb{M}$. Since the parallel transport from some point $x$ to some
point $y$ does not depend on the curve which joins $x$ to $y$, this implies
that $\psi(x)$ is the parallel transport from $T_{0}\mathbb{M}$ to
$T_{x}\mathbb{M}$. For the same reason, if $\nabla$ is the Levi-Civita
connection, then for any constant vector field $X$ and $Y$ on $M$, we have
$\nabla_{Y}X=0$. As $\nabla$ is compatible with $g$, this implies that, for
any $X\in T_{x}\mathbb{M}$, we have $X\{g_{x}(Y,Z)\}=0$, for all
$x\in\mathbb{M}$, all $X\in T_{x}\mathbb{M}$ and any local constant vector
field $Y, Z$ around $x$. Consider the diffeomorphism $\Psi(x)=x+\operatorname{Id}_{\mathbb{M}}$ of $\mathbb{M}$. The previous property implies that $\Psi^{*}g_{0}%
=g$.\newline 
Let $(U_{1},\phi_{1})$ and $(U_{2},\phi_{2})$ be two charts such
that $U_{1}\cap U_{2}\not =\emptyset$ and if $x_{0}\in U_{1}\cap U_{2}$ then,
for $i\in\{1,2\}$, $\phi_{i}^{*}g=g_{x_{0}}$ on $\phi_{i}(U_{i}) \subset
\mathbb{M}\equiv T_{x_{0}}M$ where $g_{x_{0}}$ is the constant
pseudo-Riemannian metric on $\mathbb{M}$ defined by $g_{x_{0}}$. Then it is easy to
see that the the transition function $T_{12}$ on $U_{1}\cap U_{2}$ associated to
$T\phi_{1}$ and $T\phi_{2}$ belongs to the isotropy group of $g_{x_{0}}$.

Conversely, if $g$ is an integrable tensor structure, locally there exists a
chart $(U,\phi)$ such that $\phi^{*}g=g_{x_{0}}$ on $\phi
(U)\subset\mathbb{M}\equiv T_{x_{0}}M$ where $g_{x_{0}}$ is the constant
pseudo-Riemannian defined on $\mathbb{M}\equiv T_{x_{0}}M$ by $g_{x_{0}}$.
Then it is clear that the curvature of the Levi-Civita connection of the constant
metric $g_{x_{0}}$ vanishes. Now since $\phi$ is an isometry from $(U,g_{|U})$
onto $(\phi(U),(g_{x_{0}})_{| \phi(U)})$, the curvature of the Levi-Cevita connection of
$g$ must also vanish on $U$, which ends the proof.
\end{proof}

\begin{corollary}
\label{C_CharacterizationIntegrabilityKreinMetricOnABanachManifold}
A Krein metric $g$ on a Banach manifold $M$ is an integrable tensor structure
if and only if the curvature of the Levi-Civita connection of $g$ vanishes.
\end{corollary}

\subsection{Integrability of almost tangent,  para-complex and decomposable complex structures}
\label{__IntegrabilityAlmostTangentStructuresOnABanachManifold}

\begin{definition}
\label{D_AlmostTangentComplexStructureOnABanachManifold} 
An almost tangent structure $J$ (resp. para-complex structure $\mathcal{J}$, resp. decomposable complex) structure on a Banach manifold $M$ is an almost tangent (resp. para-complex structure $\mathcal{J}$, resp. decomposable complex) structure
structure on the tangent bundle $TM$ of $M$.
\end{definition}

If $J$ is an almost tangent structure, we have a decomposition  $TM=\ker J\oplus L$,  such that $K$ and $L$ are isomorphic sub-bundles of $TM$. Moreover, there exists an isomorphism $J_L:L\rightarrow \ker J$ such that ${J}$ can be written as a matrix of type  
$\begin{pmatrix}
0 &J_L\\
0&0\\
\end{pmatrix}$.\\

If $\mathcal{J}$  is an almost para-complex  structure, there exists a Whitney decomposition $TM=E^+\oplus E^-$ where $E^+$ and $E^-$ are the eigen-bundles associated to the eigenvalues $+1$ and $-1$ of $\mathcal{J}$ respectively.

%For a global proof of the integrability of such structure  we consider a bundle morphism $\mathcal{J}$ of $TM$ such that $\mathcal{J}^2=\lambda Id_{TM}$ with $\lambda=0$ or $\lambda=-1$. In this way, we have a decomposition  $TM=E_1\oplus E_2$ , where $E_1=\ker \mathcal{J}$ for $\lambda=0$, such that $E_1$ and $E_2$ are isomorphic sub-bundle of $TM$. According to this decomposition there exists an isomorphism $I:E_2\longrightarrow E_1$ such that $\mathcal{J}$ can be written as a matrix of type  $\begin{pmatrix} 0 &-A\\
%-\lambda A^{-1}&0\\
%\end{pmatrix}$.

%For an almost tangent structure $J$ on the manifold $M$, $\ker J$ is a
%supplemented sub-bundle of $TM$ and if ${L}$ is a supplemented sub-bundle of
%$\ker J$, and the restriction $J_{L}$ of $J$ to $L$ is an isomorphism onto
%$\ker J$.
Now for a classical criterion of the integrability of such  structures we need the notion of Nijenhuis tensor.

\begin{definition}
\label{D_NijenhuisTensor} 
If $A$ is an endomorphism of $TM$, the Nijenhuis tensor of $A$ is defined, for all (local) vector fields $X$ and $Y$ on $M$, by:
\[
N_{A}(X,Y)= [AX,AY]-A[AX,Y]-A[X,AY]+A^{2}[X,Y].
\]
\label{eq_NijenhuisTensor}
\end{definition}

\begin{theorem}
\label{T_CaracterizationOfIntegrableAlmostTangentParaComplexDecomposableComplexStructure} 
An almost tangent   (resp. para-complex) structure $J$ (resp. $\mathcal{J}$) structure on $M$  is integrable if and only its Nijenhuis tensor is  null. 
\end{theorem}

\begin{proof}
For any almost tangent structure $J$, since $J^{2}=0$, we have:
\[
N_{J}(X,Y)=[JX,JY]-J[JX,Y]-J[X,JY].
\]
Moreover, if $X$ or $Y$ is a section of $\ker J$, we always have $N_{J}(X,Y)=0$.  Therefore we only
 have to consider $N_{J}$ in restriction to $L$. In particular, we can note that if  $N_{J}\equiv0$,
then we have
\[
[JX,JY]=J[JX,Y]+J[X,JY].
\]
Therefore, when we restrict this relation to a section of $L$, since $J_{L}$ is an isomorphism, this implies that $\ker J$ is an involutive supplemented sub-bundle of $TM$. \\
 
%Note that $N_J\equiv 0$ is equivalent  to $\ker J$ is involutive.
Fix some $x_{0}\in M$ and consider a chart
$(U,\phi)$ around $x_{0}$ such that $(\ker J)_{| U}$ and $K_{| U}$ are
trivial. Therefore, $T\phi(TM)=\phi(U)\times\mathbb{M}$, $T\phi(\ker
J)=\phi(U)\times\mathbb{K}$ and $T\phi(L)=\phi(U)\times\mathbb{L}$.

If $\mathbb{K}\subset\mathbb{M}$ and $\mathbb{L}\subset\mathbb{M}$ are the
typical fibres of $\ker J$ and $L$ respectively, we have $\mathbb{M}%
=\mathbb{K}\oplus\mathbb{L}$. Therefore, without loss of generality, we may
assume that $U$ is an open subset of $\mathbb{M}\equiv \mathbb{K}\times \mathbb{L}$ and $TM=U\times\mathbb{M}$,
$\ker J=U\times\mathbb{K}$ and $L=U\times\mathbb{L}$. Thus  $x\mapsto(J_{L})_{x}$ is
a smooth field which takes values in the set $\operatorname{Iso}(\mathbb{L}, \mathbb{K})$ of
isomorphisms from $\mathbb{L}$ to $\mathbb{K}$ and $x\mapsto J_{x}$ is a smooth
field from $M$ to ${L}(\mathbb{M})$ such that $\ker J_{x}=\mathbb{K}$.

As in \cite{Bel1},  Proposition 2.1 we obtain in our context:

\begin{lemma}
\label{L_ExpressionOfNJ} 
With the previous notations, we have
\[
N_{J}(X,Y)=J^{\prime}(Y,JX)-J^{\prime}(JX,Y)+J^{\prime}(X,JY)-J^{\prime
}(X,JY)
\]
where $J^{\prime}$ stands for the differential of $J$ as a map from $M$ to 
$L(\mathbb{M})$.
\end{lemma}

At first, assume that $J$ is an integrable almost tangent structure. This means that,
for any $x\in M$, the value $J_{x}$ in ${L}(\mathbb{M})$ (resp. $(J_{L})_{x}$
in $\operatorname{Iso}(\mathbb{L}, \mathbb{K})$) is a constant. So from Lemma
\ref{L_ExpressionOfNJ}. it follows that $N_{J}=0$.\newline

Conversely,  assume that $N_{J}=0$ and as we have already seen $\ker J$ is involutive. If $X$ and $Y$ are sections of $L$, in the expression of $N_{J}$ given in Lemma \ref{L_ExpressionOfNJ}, we can replace $J$ by $J_{L}$ and this implies that  $J_L$ satisfies the assumption of IV., Theorem 1.2 in \cite{Lan}.   By the same arguments as in the proof of the Frobenius Theorem in \cite{Lan}, we  produce  a local diffeomorphism $\Psi$ from a neighbourhood
$U_{0}\times V_{0}\in \mathbb{K}\times \mathbb{L}$ of $x_{0}$ such that $\Psi^{*}J_L$ is a smooth field from
$U_{0} \times V_{0}$ to $\operatorname{Iso}(\mathbb{K},\mathbb{L})$ which is constant. Thus
the same is true for $J$, which ends the proof in this case.\\

Let us consider now a para-complex structure $\mathcal{J}$. Recall that we have  a canonical Whitney decomposition $TM=E^+\oplus E^-$ where $E^+$ (resp. $E^-$) is the eigen-bundle associated to the eigenvalue $+1$ (resp. $-1$) of $\mathcal{J}$. Note that if $P^{\pm}=\frac{1}{2}( \operatorname{Id}\pm \mathcal{J})$, then $E^{\pm}=\textrm{im } P^{\pm}$, $\ker P^{\pm}=E^{\mp}$  and $TM=E^+\oplus E^-$. Now, we have   $N_\mathcal{J}= 0$ if and only if  $N_{P^+}=N_{P^-}=0$  and  $E^+$ and $E^-$ are integrable sub-bundles of $TM$.
Note that  since $\mathcal{J}^2=Id$, Lemma \ref{L_ExpressionOfNJ}  is also valid for $\mathcal{J}$; so if $\mathcal{J}$ is integrable, this implies that $N_\mathcal{J}= 0$.\\

Conversely, assume that $N_\mathcal{J}= 0$.  Fix some $x_0\in M$ and denote by $\mathbb{E}^{\pm}$ the typical  fiber of  $E^{\pm}$. If  the sub-bundles $E^{\pm}$ are integrable, from Frobenius Theorem, there exists a chart $(U^{\pm},\phi_{\pm})$  around $x_0$ where $\phi_{\pm}$ is a diffeomorphism from $U^{\pm}$ onto an open neighbourhood $ V_1^{\pm}\times V_1^{\pm}$ in $\mathbb{E}^{+}\times\mathbb{E}^{-}$ such that   $\phi_{+}^*P^{+}$ and  $\phi_{-}^*P^{-}$ are the fields of constant matrices
  $\begin{pmatrix}
0&\operatorname{Id}_{\mathbb{E}^+}\\
0&0\\
\end{pmatrix}$ and $\begin{pmatrix}
0&0\\
0&\operatorname{Id}_{\mathbb{E}^-}\\
\end{pmatrix}$.
Now the transition map $\phi_{-}\circ\phi_{+}^{-1}$ is necessarily of type $(\bar{x},\bar{y})\mapsto (\alpha(\bar{x}),\beta(\bar{y}))$. So the restriction of $\phi_{+}$ to $U^+\cap U^-$ is such that $\phi_{+}^*P^{+}$ and  $\phi_{+}^*P^{-}$ are also the previous matrices, which ends the proof since $\mathcal{J}=P^++P^-$.
\end{proof}

Unfortunately,  the problem of integrability of an almost complex structure  $\mathcal{I}$ on a Banach  manifold $M$ is not equivalent to the relation $N_{\mathcal{I}}\equiv 0$. The reader can find in \cite{Pat}  an example of an almost complex structure  $\mathcal{I}$  on a smooth Banach manifold $M$ such that $N_{\mathcal{I}}\equiv 0$ for which there exists {\it no holomorphic chart} in which $\mathcal{I}$ is isomorphic to a linear complex structure.  
However, if $M$ is analytic and if  $N_{\mathcal{I}}\equiv 0$ there exists an holomorphic chart in which $\mathcal{I}$ is isomorphic to a linear complex structure (Banach version of  Newlander-Nirenberg theorem see \cite{Bel1}). In particular, there exists a structure of holomorphic manifold on $M$.  

\begin{definition}
\label{D_FormalntagrebilityComplexStructure} 
An almost complex structure $\mathcal{I}$ on a Banach manifold $M$  is called formally integrable if $N_{\mathcal{I}}\equiv 0$. \\
$\mathcal{I}$ is called integrable if there exists a structure of complex manifold on $M$
\end{definition}

\subsection{Flat K\"ahler and para-K\"ahler Banach manifolds}
 Following  \cite{Tum1}, we introduce the following notions.
\begin{definition}
\label{D_KalherBanachStructure}
A K\"ahler (resp. formal-K\"ahler) Banach structure on a Banach manifold $M$ is an almost K\"ahler structure $(\Omega,g,\mathcal{I})$ on $TM$ such that the almost complex structure  $\mathcal{I}$  is  integrable (resp. formally integrable).
\end{definition}
   %In this case, as we have seen in subsection \ref{__CompatibleWeakSymplecticFormWeakRiemannianMetricAlmostComplexStructureOnABanachBundle}, we have a decomposition $TM=E_1\oplus E_2$ where $E_1$ and $E_2$ are isomorphic sub-bundle, Lagrangian and orthogonal sub-bundle of $TM$.  
 Note that  if $M$ is analytic  and $(\Omega,g,\mathcal{I})$ on $TM$ is an almost  K\"ahler structure $(\Omega,g,\mathcal{I})$ on $TM$ if an only if $N_\mathcal{I}\equiv 0$  then we have a K\"ahler structure on $M$. The reader will find many examples of weak and strong formal-K\"ahler Banach manifolds in \cite {Tum2}.\\

 In the same way, a  {\it   weak (resp. strong)  para-K\"ahler  Banach structure} on a Banach manifold $M$ is an almost  para-K\"ahler structure $(\Omega,g,\mathcal{J})$ on $TM$   (cf. Definition \ref{D_KahlerBanachBundle}) {\it such that the almost para-complex structure  $\mathcal{J}$  is  integrable}. 
 
\smallskip 
 
\textbf{In what follows, we use the terminology K\"ahler or para-K\"ahler manifold  instead of weak  K\"ahler or weak  para-K\"ahler manifold.}\\
 
 According to $\S$  
\ref{__CompatibleWeakSymplecticFormWeakRiemannianMetricAlmostComplexStructuresOnABanachBundle}, if we have a K\"ahler or para-K\"ahler manifold structure on $M$, we have a decomposition $TM=E_1\oplus E_2$ where $E_1$ and $E_2$ are isomorphic sub-bundles, Lagrangian and orthogonal sub-bundles of $TM$.
\begin{definition}
\label{D_FlatFormalBanachStructure}
A  formal-K\"ahler (resp. para-K\"ahler) Banach structure  $(\Omega,g,\mathcal{I})$ (resp. $(\Omega,g,\mathcal{J})$) is said to be flat if the Levi-Civita connection $\nabla$ of $g$ exists and is flat.
\end{definition}
 
\begin{lemma}
Let $(\omega,g,\mathcal{I})$ (resp. $(\Omega,g,\mathcal{J})$)  be  a K\"ahler (resp. para-K\"ahler) Banach structure such that the Levi-Civita connection $\nabla$ of $g$ exists. Then we have 
 $\nabla \mathcal{I}\equiv 0$ (resp.  $\nabla \mathcal{J}\equiv 0$)
\end{lemma}
\begin{proof}
cf. proof of \cite{Tum2}, Proposition 91. 
\end{proof}
 
\begin{theorem}
\label{T_FlatKahlerBanachManifold} 
If a  formal-K\"ahler (resp. para-K\"ahler) Banach structure  $(\Omega,g,\mathcal{I})$ 
 (resp. $(\Omega,g,\mathcal{J})$) on $M$ is flat, for each $x_0\in M$, there exists a chart $(U,\phi)$ and a linear Darboux form $\Omega_0$, a linear decomposable complex 
  structure $\mathcal{I}_0$, an inner product $g_0$ (resp.  a linear Darboux form $\Omega_0$, a linear para-complex structure $\mathcal{J}_0$ and a neutral inner product $g_0$) on
  $\mathbb{M}$ such that $\phi^*\Omega_0=\Omega$, $\phi^*\mathcal{I}_0=\mathcal{I}$ and $\phi^*g_0=g$ (resp. $\phi^*\Omega_0=\Omega$, $\phi^*\mathcal{J}_0=\mathcal{J}$ 
   and $\phi^*g_0=g$).  
\end{theorem}

\begin{proof}%At first note that the parallel transport along any curve $\gamma=[0,1]\longrightarrow M$  is an isometry from $T_{\gamma(0)}M$ to $T_{\gamma(1)}M$ which induces an isometry from $(E_i)_{\gamma(0)}$ to $(E_i)_{\gamma(1)}$ for $i=1,2$. 
We only consider the case of an almost  para-K\"ahler manifold, the case of an almost K\"ahler manifold is similar.\\
 
Consider the decomposition $TM=E_1\oplus E_2$ as recalled previously.  Since $\nabla g=0$, and $E_1$ and $E_2$ are orthogonal, the connection $\nabla$ induces a connection $\nabla^i$ on $E_i$ which preserves the restriction $g_i$ of $g$ to $E_i$  for $i=1,2$ and $\nabla=\nabla^1+\nabla^2$. Thus, for $i=1,2$,  if $X$ and $Y$ are sections of $E_i$ we have 
$$[X,Y]=\nabla_XY-\nabla_YX=\nabla^i_XY-\nabla^i_YX$$ and so $E_i$ is integrable. Moreover the parallel transport from $T_{\gamma(0)}M$ into $T_{\gamma(1)}M$ along an curve $\gamma:[0,1]\longrightarrow M$ induces an isometry of the fibre $(E_i)_{\gamma(0)}$ into $(E_i)_{\gamma(1)}$ for $i=1,2$.
 Fix some $x_0\in M$. Since $\nabla$ is flat,  according to the proof Theorem \ref{T_gIntegrableTensorStructureNullCurvature}, we may assume that $M=\mathbb{M}$ and 
$x_0=0\in \mathbb{M}$. Moreover, if $\mathbb{E}_1$ (resp. and $ \mathbb{E}_2$) is the typical fiber of the integrable sub-bundle $E_1$ (resp. $E_2$), we may assume that 
$\mathbb{M}=\mathbb{E}_1\oplus \mathbb{E}_2$.  Recall that under the flatness assumption, the diffeomorphism $\Psi(x)= x+Id_\mathbb{M}$ is such that $T_0\Psi$ is the parallel 
transport from $T_0\mathbb{M}$ to $T_x\mathbb{M}$ and $\Psi^0g_0=g$. But since $E_1$ and $E_2$ are invariant by parallel transport, and and $\mathcal{J}$ is invariant by 
parallel transport we must have $\mathcal{J}\circ T\Psi=T\Psi\circ \mathcal{J}$ and so $\Psi^*\mathcal{J}_0=\mathcal{J}$. Now as $\Omega(u,v)=g(u,\mathcal{J}v)$ we also get $\Psi^*\Omega_0=\Omega$.
\end{proof}

 \begin{remark}
 \label{R_CryterionIntegrabilityDecomposableComplex}
Under the assumption of  flatness of $\nabla$ in Theorem   \ref{T_FlatKahlerBanachManifold} we obtain the integrability  of $\mathcal{I}$.   But as we have already seen, in general the condition of nullity of $N_\mathcal{I}$ is not sufficient to ensure the integrability of $\mathcal{I}$.
\end{remark}

\section{Projective limits of tensor structures}
\label{_ProjectiveLimitsOfTensorStructures}

The reference for this section is the book \cite{DGV}.

\subsection{Projective limits}
\label{__ProjectiveLimits}

\begin{definition}
\label{D_ProjectiveSystem}
Let $\left(  I,\leq\right)  $ be a directed set and
$\mathbb{A}$ a category.\\
$\mathcal{S}=\left\{ \left(  X_{i},\mu_{i}^{j}\right) \right\}_{\left(  i,j\right)  \in I^{2},\ i\leq j}$
 is called a projective system if
\begin{enumerate}
\item[--]
For all $i \in I$, $X_i$ is an object of the category;
\item[--] 
For all $(i,j)\in I^2:j \geq i$, $\mu_{i}^{j}:X_{j}\longrightarrow X_{i}$ is a morphism (bonding map)
\end{enumerate} 
such that:

\begin{description}
\item[\textbf{(PS 1)}] 
$\forall i \in I, \; \mu_{i}^{i}=\operatorname{Id}_{X_{i}}$;

\item[(PS 2)] 
$\forall (i,j,k) \in I^3:i\leq j\leq k, \; \mu_{i}^{j}\circ\mu_{j}^{k}=\mu_{i}^{k}$.
\end{description}
\end{definition}

\begin{definition}
\label{D_ProjectiveLimit}
Let $\mathcal{S}=\left\{ \left(  X_{i},\mu_{i}^{j}\right) \right\}_{\left(  i,j\right)  \in I^{2},\ i\leq j}$ be a projective system.\\
The set
\[
X=\left\{  \left(  x_{i}\right)  \in\prod\limits_{i\in I}X_{i}:\forall\left(
i,j\right)  \in I^{2},\mu_{i}^{j}\left(  x_{j}\right)  =x_{i}\right\}
\]
is called the projective limit of the system $\mathcal{S}$ and is denoted by
$\underleftarrow{\lim}X_{i}$.
\end{definition}

When $I=\mathbb{N}$ with the usual order relation, countable projective
systems are called \textit{projective sequences}.

\begin{remark}
\label{R_SurjectiveInjectiveBondingMaps}
The bonding maps can be surjective, e.g. projective limits of $k$-jets bundles (cf. \cite{Sau}) or injective, e.g. ILB structures as introduced in \cite{Omo}.
\end{remark}

\subsection{Fr\'{e}chet spaces}
\label{__FrechetSpaces}

In Analysis and Mathematical Physics, Banach representations break down. By
weakening the topological requirement, replacing the norm by a sequence of
semi-norms, we get the notion of Fr\'{e}chet space.

\begin{definition}
\label{D_FrechetSpace}A Fr\'{e}chet space is a topological vector space which
is locally convex, Hausdorff, metrizable and complete.
\end{definition}

\begin{example}
\label{Ex_CountableProductOfR}
The space $\mathbb{R}^{\mathbb{N}}=\prod\limits_{n\in\mathbb{N}}\mathbb{R}^{n}$, endowed with the cartesian topology, is a Fr\'{e}chet space with corresponding sequence of semi-norms
$\left(  p_{n}\right)  $ defined by
\[
p_{n}\left(  x_{0},x_{1},\dots\right)  =\sum\limits_{i=0}^{n}x_{i}%
\]
where the distance is given by
\[
d\left(  x,y\right)  =\sum\limits_{n=0}^{+\infty}\dfrac{\left\vert x_{n}%
-y_{n}\right\vert }{2^{n}\left(  1+\left\vert x_{n}-y_{n}\right\vert \right)
}%
\]
(\cite{DGV}, 2.1, Examples 2.1.5., 2.).
\end{example}

Some important properties of Fr\'{e}chet spaces are analogous to Banach ones:

\begin{enumerate}
\item[--] Every closed subspace of a Fr\'{e}chet space is also a Fr\'{e}chet space;

\item[--] The open mapping theorem as well as the Hahn-Banach hold true in
Fr\'{e}chet spaces.
\end{enumerate}

\smallskip

But a number of critical deficiencies emerge in this framework:

\begin{enumerate}
\item[--] The space of continuous linear maps $\mathcal{L}\left(
\mathbb{F}_{1},\mathbb{F}_{2}\right)  $ between two Fr\'{e}chet spaces
$\mathbb{F}_{1}$ and $\mathbb{F}_{2}$ is not necessarly a Fr\'{e}chet space;

\item[--] The inverse function theorem is not valid in general;

\item[--] There is no general solvability theory for differential equations in
Fr\'{e}chet spaces.
\end{enumerate}

However, one can consider, under sufficient conditions, geometrical objects,
structures and properties in the Fr\'{e}chet context as projective limits of
sequences of their Banach components. So different pathological entities in
the Fr\'{e}chet framework can be replaced by approximations compatible with
the limit process, e.g. ILB-Lie groups (\cite{Omo}) or projective limits of
Banach Lie groups (\cite{Gal1}), manifolds (\cite{AbMa}), bundles
(\cite{Gal2}, \cite{AgSu}), algebroids (\cite{Cab}), connections and
differential equations (\cite{ADGS}).

\subsection{Projective limits of Banach manifolds}
\label{__ProjectiveLimitsOf BanachManifolds}

\begin{definition}
\label{D_ProjectiveLimitOfBanachManifolds}Let $\left(  M_{i},\mu_{i}%
^{j}\right)  _{\left(  i,j\right)  \in\mathbb{N}^{\ast}\times\mathbb{N}^{\ast
},\ i\leq j\ }$be a projective sequence of Banach $C^{\infty}$-manifolds,
where $\mu_{i}^{j}:M_{j}\longrightarrow M_{i}$ are smoth maps, $M_{i}$ being
modelled on the Banach spaces $\mathbb{M}_{i}$. \newline The space
$\underleftarrow{\lim}M_{i}$ is called a projective limit of Banach manifolds,
provided that the sequence of models $\left(  \mathbb{M}_{n}\right)
_{n\in\mathbb{N}^{\ast}}$ forms a projective sequence with connecting
morphisms $\overline{\mu_{i}^{j}}:\mathbb{M}_{j}\longrightarrow\mathbb{M}_{i}$
whose projective limit is the Fr\'{e}chet space $\mathbb{M=}\underleftarrow
{\lim}\mathbb{M}_{n}$ and has the projective limit chart property at any point:

\begin{description}
\item[\textbf{(PLCP)}] for all $x=\left(  x_{n}\right)  \in M=\underleftarrow
{\lim}M_{n}$, there exists a projective system of local charts $\left(
U_{n},\phi_{n}\right)  _{n\in\mathbb{N}^{\ast}}$ such that $x_{n}\in U_{n}$
where $\phi_{i}\circ\mu_{i}^{j}=\overline{\mu_{i}^{j}}\circ\phi_{j}$ and where
$U=\underleftarrow{\lim}U_{n}$ is open in $M$.
\end{description}
\end{definition}

In this situation, the projective limit $M=\underleftarrow{\lim}M_{n}$ has a
structure of Fr\'{e}chet manifold modelled on the Fr\'{e}chet space
$\mathbb{M}$ where the differentiable structure is defined via the charts
$\left(  U,\phi\right)  $ where $\phi=\underleftarrow{\lim}\phi_{n}%
:U\rightarrow\left(  \phi_{n}\left(  U_{n}\right)  \right)  .$\newline$\phi$
is a homeomorphism (projective limit of homeomorphisms) and the charts
changings
\[
\left(  \phi^{\alpha}\circ\left(  \phi^{\beta}\right)  ^{-1}\right)
_{|\phi^{\beta}\left(  U^{\alpha}\cap U^{\beta}\right)  }=\underleftarrow
{\lim}\left(  \left(  \phi_{n}^{\alpha}\circ\left(  \phi_{n}^{\beta}\right)
^{-1}\right)  _{|\phi_{n}^{\beta}\left(  U_{n}^{\alpha}\cap U_{n}^{\beta
}\right)  }\right)
\]
between open sets of Fr\'{e}chet spaces are smooth.

\begin{example}
\label{Ex_InverseLimitManifolds}In \cite{Omo}, H. Omori introduced the notion
of inverse manifolds as the intersection of a descending sequence of Banach
manifolds $M_{1}\supset M_{2}\supset\cdots$ for which local charts can be
defined on $M=\bigcap\limits_{n\in\mathbb{N}^{\ast}}M_{n}$.
\end{example}

\subsection{Projective limits of Banach-Lie groups}
\label{__ProjectiveLimitsOfBanachLieGroups}

The group of diffeomorphisms of a compact manifold studied in \cite{Les} and
\cite{Omo} appears as a projective limit of Lie groups.

Let $\left(  G_{i},\gamma_{i}^{j}\right)  _{_{\left(  i,j\right)
\in\mathbb{N}^{\ast}\times\mathbb{N}^{\ast},\ i\leq j}}$be a projective sequence
of groups where $\gamma_{i}^{j}:G_{j}\longrightarrow G_{i}$ are group
homomorphisms. The projective limit $G=\underleftarrow{\lim}G_{n}$ can be
endowed with a structure of group in an obvious way.

\begin{definition}
\label{D_ProjectiveSequenceOfBanachLieGroups}
A projective sequence of Banach-Lie groups is a projective sequence of groups
satisfying the conditions of Definition \ref{D_ProjectiveLimitOfBanachManifolds} where $\gamma_{i}^{j}$ are morphisms
of Banach-Lie groups.
\end{definition}

\begin{proposition}
\label{P_StructureOfFrechetLieGroupOnProjectiveLimitOfBanachLieGroups}
Let
$\left(  G_{i},\gamma_{i}^{j}\right)  _{\left(  i,j\right)  \in\mathbb{N}
^{\ast}\times\mathbb{N}^{\ast},\ i\leq j}$ be a projective sequence of
Banach-Lie groups where $G_{i}$ is modelled on the Banach space $\mathbb{G}%
_{i}$. \newline 
Then the projective limit of Banach manifolds $G=\underleftarrow{\lim}G_{i}$ can be endowed with a structure of Lie group modelled on the Fr\'{e}chet space $\mathbb{G}=\underleftarrow{\lim}
\mathbb{G}_{i}$.
\end{proposition}

In this case $G$ is called a projective limit of Banach-Lie groups.

We then have the fundamental properties for any projective limit of Banach-Lie
group (\cite{Gal1}, Theorem 1.3).

\begin{proposition}
\label{P_PropertiesOfProjectiveLimitsOfBanachLieGroups} Any projective limit
of Banach Lie groups $G=\underleftarrow{\lim}G_{i}$ has the following properties:

\begin{description}
\item[\textbf{(PLBLG 1)}] 
$G$ is a Fr\'{e}chet Lie group;

\item[\textbf{(PLBLG 2)}] 
$\forall x=\left(  x_{i}\right)  \in G,\ T_{x}%
G=\underleftarrow{\lim}T_{x_{i}}G_{i}$;

\item[\textbf{(PLBLG 3)}] 
$TG=\underleftarrow{\lim}TG_{i}$;

\item[\textbf{(PLBLG 4)}] 
If $X$ is a left invariant vector field of $G$ (i.e. $X$ belongs to the Lie algebra $\mathfrak{g}$), there exists a unique integral curve of $X$ through $e$.
\end{description}
\end{proposition}

As a consequence of (PLBLG 4), $G$ admits an exponential map $\exp_{G}%
:T_{e}G\longrightarrow G$ such that $\exp_{G}=\underleftarrow{\lim}\exp G_{i}%
$. \newline Unfortunately, $\exp_{G}$ is not necessarily a local
diffeomorphism at $0\in T_{e}G$. This property becomes true for commutative
Fr\'{e}chet-Lie groups in the sense of \cite{Gal1}.

\begin{example}
\label{Ex_GroupSmoothMapsfromCompactManifolfsToFiniteDimensionalLieGroup}
The group of smooth maps $C^{\infty}\left(  M,G\right)  =\bigcap\limits_{n\in
\mathbb{N}^{\ast}}C^{n}\left(  M,G\right)  $ from a compact manifold $M$ to a
finite dimensional Lie group is a Fr\'{e}chet projective limit of Banach Lie
groups (\cite{DGV}, Examples 3.1.6.,4.) under the pointwise operations; the
Lie groups $C^{n}\left(  M,G\right)  $ are sometimes referred to as group of
currents and, in the case $M=\mathbb{S}^{1}$, as loop groups (\cite{OnVi}%
).\newline The differential structure is given by the charts $\left\{  \left(
C^{n}\left(  M,N\right)  ,\psi_{n}\right)  \right\}  _{n\in\mathbb{N}^{\ast}}$
where $N$ is an open neighbourhood of $e\in G$ such that the exponential map
of $G$ $\exp_{G}:\exp_{G}^{-1}\left(  N\right)  \longrightarrow N$ is a
diffeomorphism and $\psi_{n}\left(  f\right)  =\exp_{G}^{-1}\circ f$
(\cite{Gal1}, 1. Examples (iv)).
\end{example}

\subsection{The Fr\'{e}chet topological group $\mathcal{H}^{0}\left(\mathbb{F}\right) $}
\label{__TheFrechetTopologicalGroupH0(F)}

Let $\mathbb{F}_{1}$ and $\mathbb{F}_{2}$ be Fr\'{e}chet spaces. The space of
continuous linear mappings $\mathcal{L}\left(  \mathbb{F}_{1},\mathbb{F}%
_{2}\right)  $ between these spaces drops out the Fr\'{e}chet category.
However, using the realization of the Fr\'{e}chet space $\mathbb{F}_{1}$
(resp. $\mathbb{F}_{2}$) as projective limit of Banach spaces, say
$\mathbb{F}_{1}=\underleftarrow{\lim}\mathbb{E}_{1}^{n}$ (resp. $\mathbb{F}%
_{2}=\underleftarrow{\lim}\mathbb{E}_{2}^{n}$), where $\left(  \overline{\mu
}_{1,i}^{j}\right)  _{i\leq j}$ (resp. $\left(  \overline{\mu}_{2,i}%
^{j}\right)  _{i\leq j}$) are the bonding maps$,$ $\mathcal{L}\left(
\mathbb{F}_{1},\mathbb{F}_{2}\right)  $ can be replaced by a new space within
the Fr\'{e}chet framework, as defined in \cite{Gal2}.

For each $n\in\mathbb{N}^{\ast},$ we define the set%
\[
\mathcal{H}^{n}\left(  \mathbb{F}_{1},\mathbb{F}_{2}\right)  =\left\{  \left(
f_{i}\right)  _{1\leq i\leq n}\in\prod\limits_{i=1}^{n}\mathcal{L}\left(
\mathbb{E}_{1}^{i},\mathbb{E}_{2}^{i}\right)  :\forall j\in\left\{
i,\dots,n\right\}  :f_{i}\circ\overline{\mu}_{1,i}^{j}=\overline{\mu}%
_{2,i}^{j}\circ f_{j}\right\}  .
\]

$\mathcal{H}^{n}\left(  \mathbb{F}_{1},\mathbb{F}_{2}\right)  $ is a Banach
space as a closed space of the Banach space $\prod\limits_{i=1}^{n}%
\mathcal{L}\left(  \mathbb{E}_{1}^{i},\mathbb{E}_{2}^{i}\right)  $.

\begin{proposition}
\label{P_StructureFrechetSpaceHF1F2}
The space $\mathcal{L}\left(\mathbb{F}_{1},\mathbb{F}_{2}\right)  $ can be represented as the Fr\'{e}chet space
\[
\mathcal{H}\left(  \mathbb{F}_{1},\mathbb{F}_{2}\right)  =\left\{  \left(
f_{n}\right)  \in\prod\limits_{n\in\mathbb{N}^{\ast}}\mathcal{L}\left(
\mathbb{E}_{1}^{n},\mathbb{E}_{2}^{n}\right)  :\underleftarrow{\lim}%
f_{n}\text{ exists}\right\}
\]
isomorphic to the Fr\'{e}chet space $\underleftarrow{\lim}\mathcal{H}^{n}\left(
\mathbb{F}_{1},\mathbb{F}_{2}\right)  $.
\end{proposition}

Now we consider
\[
\mathcal{H}_{0}^{n}\left(  \mathbb{F}_{1},\mathbb{F}_{2}\right)
=\mathcal{H}^{n}\left(  \mathbb{F}_{1},\mathbb{F}_{2}\right)  \bigcap
\prod\limits_{i=1}^{n}\mathcal{L}is\left(  \mathbb{E}_{1}^{i},\mathbb{E}%
_{2}^{i}\right)  .
\]

For a Fr\'{e}chet space $\mathbb{F}$, we denote $\mathcal{H}_{0}^{n}\left(
\mathbb{F},\mathbb{F}\right)  $ (resp. $\mathcal{H}_{0}\left(  \mathbb{F}%
,\mathbb{F}\right)  $) by $\mathcal{H}_{0}^{n}\left(  \mathbb{F}\right)  $
(resp. $\mathcal{H}_{0}\left(  \mathbb{F}\right)  $). \newline 
We then have the following result (\cite{DGV}, Proposition 5.1.1)

\begin{proposition}
\label{P_FrechetTopologicalGroupH0F}
Every $\mathcal{H}_{0}^{n}\left(\mathbb{F}\right)  $ is a Banach Lie group modelled on $\mathcal{H}^{n}\left(\mathbb{F}\right)  $.\newline 
Moreover the projective limit $\underleftarrow{\lim}\mathcal{H}_{0}^{n}\left(  \mathbb{F}\right)  $ exists and coincides, up to an isomorphism of topological groups, with $\mathcal{H}_{0}\left(
\mathbb{F}\right)  $. \newline 
Thus $\mathcal{H}_{0}\left(  \mathbb{F}\right)$ is a Fr\'{e}chet topological group.
\end{proposition}

In this situation, the maps $\left(  h_{0}\right)  _{i}^{j}:\mathcal{H}%
_{0}^{j}\left(  \mathbb{F}\right)  \longrightarrow\mathcal{H}_{0}^{i}\left(
\mathbb{F}\right)  $ are morphisms of topological groups satisfying, for every
$i\leq j\leq k$, the relations $\left(  h_{0}\right)  _{i}^{k}=\left(
h_{0}\right)  _{i}^{j}\circ\left(  h_{0}\right)  _{j}^{k}$. Thus $\left(
\mathcal{H}_{0}^{i}\left(  \mathbb{F}\right)  ,\left(  h_{0}\right)  _{i}%
^{j}\right)  $ is a projective system of Banach-Lie groups, but $\mathcal{H}%
_{0}\left(  \mathbb{F}\right)  =\underleftarrow{\lim}\mathcal{H}_{0}%
^{n}\left(  \mathbb{F}\right)  $ is not necessarily a Fr\'{e}chet-Lie group
because the projective limit of the open sets $\mathcal{H}_{0}^{n}\left(
\mathbb{F}\right)  $ is not necessarily open.

\subsection{Some  Fr\'{e}chet topological subgroups of $\mathcal{H}^{0}\left(\mathbb{E}\right)  $}
\label{__FrechetTopologicalSubgroupsOfH0E}

Let $\left(  \mathbb{E}%
_{i},\overline{\lambda_{i}^{j}}\right)  _{\left(  i,j\right)  \in
\mathbb{N}^{\ast}\times\mathbb{N}^{\ast},\ i\leq j\ }$ be  a projective
sequence of Banach spaces. Consider a  sequence $(G_n)$ of  Banach Lie groups such that, for all  $n\in \mathbb{N}^*$,  $G_n$ is a weak Banach Lie subgroup of $\operatorname{GL}(\mathbb{E}_{n})$ and denote by  $ \mathcal{L}_{G_n}(\mathbb{E}_n)\subset \mathcal{L}(\mathbb{E}_n)$  the Lie  algebra of $G_n$. We set  

$$\mathcal{G}^{n}\left(  \mathbb{E}\right)  =
\left\{ (g_i)_{1\leq i\leq n}\in\prod\limits_{j=1}^{n}\mathcal{L}_{G_j}(\mathbb{E}_j)\;:\; \forall j\in \{1,\dots,n\}\;:\; g_i\circ\overline{\lambda_{i}^{j}}=\overline{\lambda_{i}^{j}}\circ g_j\right\}$$
\noindent Note that $\mathcal{G}^{n}\left(  \mathbb{E}\right) $ is a closed vectorial subspace of  $\mathcal{H}^{n}\left(  \mathbb{E}\right) $ and so has a natural structure of Banach space.\\

For $1\leq i\leq j$, let us consider the natural projection
\[%
\begin{array}
[c]{cccc}%
\gamma_{i}^{j}: & \mathcal{G}^{j}\left(  \mathbb{E}\right)  &
\longrightarrow & \mathcal{G}^{j}\left(  \mathbb{E}\right) \\
& \left(  g_{1},\dots,g_{j}\right)  &  & \left(  g_{1},\dots,g_{i}\right)
\end{array}
.
\]
\noindent  We consider:\[
\mathcal{G}_{0}^{n}\left(  \mathbb{E}\right)
=\mathcal{G}^{n}\left(  \mathbb{E}\right)  \bigcap
\prod\limits_{i=1}^{n}\mathcal{L}is\left(  \mathbb{E}_i\right)  .
\]
By the same arguments used in the proof of \cite{DGV}, Proposition 5.1.1, we obtain
 
\begin{proposition}
\label{T_WeakFrechetSubgroup}
Every $\mathcal{G}_{0}^{n}\left(\mathbb{F}\right)  $ is a Banach Lie group modelled on $\mathcal{G}^{n}\left(\mathbb{F}\right)  $.\\ 
Moreover the projective limit $\underleftarrow{\lim}\mathcal{G}_{0}^{n}\left(  \mathbb{E}\right)  $ exists and coincides, up to an isomorphism of topological groups, with $\mathcal{G}_{0}\left(
\mathbb{E}\right)  $. \newline 
Thus $\mathcal{G}_{0}\left(  \mathbb{E}\right)$ is a Fr\'{e}chet topological group and is a closed topological subgroup of $\mathcal{H}_{0}\left(  \mathbb{E}\right)$.
\end{proposition}

The group  $\mathcal{G}_{0}\left(  \mathbb{E}\right)  $  will play the role of structure group for the
principal bundle of projective limit of $G$-structures (cf. $\S$ \ref{_ProjectiveLimitsOfTensorStructures}).\\

As in the previous section all projections $\gamma_i^j$ induce a projection   
\[
\left(  \gamma_{0}\right)  _{i}^{j}:\mathcal{G}%
_{0}^{j}\left(  \mathbb{F}\right)  \longrightarrow\mathcal{G}_{0}^{i}\left(
\mathbb{E}\right)
\] 
which are morphisms of topological groups satisfying, for all $i\leq j\leq k$, the relations
\[ 
\left(  \gamma_{0} \right)  _{i}^{k}=
\left( \gamma_{0} \right)  _{i}^{j}\circ\left(  \gamma_{0} \right)_{j}^{k}.
\]
Thus 
$\left(\mathcal{G}_{0}^{i}\left(  \mathbb{E}\right)  ,\left(  \gamma_{0}\right)  _{i}^{j}\right) $
 is a projective system of Banach-Lie groups, but 
$
\mathcal{G}_{0} \left( \mathbb{E} \right)  =\underleftarrow{\lim}\mathcal{G}_{0}^{n}\left( \mathbb{E} \right) 
$ 
is not necessarily a Fr\'{e}chet-Lie group. \\

{\it From now on, without ambiguity, these projections $\left(  \gamma_{0}\right)_{i}^{j}$ will simply denoted $\gamma_i^j$ for short.}

\subsection{Projective limits of Banach vector bundles}
\label{__ProjectiveLimitOfBanachVectorBundles}

\begin{definition}
\label{D_ProjectiveSequenceOfBanachVectorBundles} 
$\mathcal{E}=\left\{\left(  E_{n},\pi_{n},M_{n}\right)  \right\}  _{n\in\mathbb{N}^{\ast}}$ is
called a projective sequence of Banach vector bundles if the following
conditions are satisfied:

\begin{description}
\item[\textbf{(PSBVB 1)}] $\left(  M_{i},\mu_{i}^{j}\right)  _{\left(
i,j\right)  \in\mathbb{N}^{\ast}\times\mathbb{N}^{\ast},\ i\leq j\ }$ is a
projective sequence of Banach $C^{\infty}$-manifolds;

\item[\textbf{(PSBVB 2)}] $\left(  E_{i},\lambda_{i}^{j}\right)  _{\left(
i,j\right)  \in\mathbb{N}^{\ast}\times\mathbb{N}^{\ast},\ i\leq j\ }$ is a
sequence of Banach $C^{\infty}$-manifolds and $\left(  \mathbb{E}%
_{i},\overline{\lambda_{i}^{j}}\right)  _{\left(  i,j\right)  \in
\mathbb{N}^{\ast}\times\mathbb{N}^{\ast},\ i\leq j\ }$ is a projective
sequence of Banach spaces where $\mathbb{E}_{i}$ is the fibre of $E_{i}$;

\item[\textbf{(PSBVB 3)}] For all $i\leq j$, $\pi_{j}\circ\lambda_{i}^{j}%
=\mu_{i}^{j}\circ\pi_{i};$\newline

\item[\textbf{(PSBVB 4)}] $\left(  M_{n}\right)  _{n\in\mathbb{N}^{\ast}}$ has
the projective limit chart property \emph{(PLCL)} at any point relatively to
$(U=\underleftarrow{\lim}U_{n},\phi=\underleftarrow{\lim}\phi_{n})$~;

\item[\textbf{(PSBVB 5)}] For any $x\in M=\underleftarrow{\lim}M_{n}$, there
exists a local trivialization $\tau_{n}:\pi_{n}^{-1}\left(  U_{n}\right)
\longrightarrow U_{n}\times\mathbb{E}_{n}$ such that the following diagram is
commutative:%
\[%
\begin{array}
[c]{ccc}%
\left(  \pi_{i}\right)  ^{-1}\left(  U_{i}\right)  & \underleftarrow
{\lambda_{i}^{j}} & \left(  \pi_{j}\right)  ^{-1}\left(  U_{j}\right) \\
\tau_{i}\downarrow &  & \downarrow\tau_{j}\\
U_{i}\times\mathbb{E}_{i} & \underleftarrow{\mu_{i}^{j}\times\overline
{\lambda_{i}^{j}}} & U_{j}\times\mathbb{E}_{j}%
\end{array}
\]

\end{description}
\end{definition}

Adapting the result of \cite{DGV}, 5.2, we get the following theorem:

\begin{theorem}
\label{T_FrechetStructureOnProjectiveLimitOfBanachVectorBundles}Let $\left\{
\left(  E_{n},\pi_{n},M_{n}\right)  \right\}  _{n\in\mathbb{N}^{\ast}}$ be a
projective sequence of Banach vector bundles. The triple $\left(
\underleftarrow{\lim}E_{n},\underleftarrow{\lim}\pi_{n},\underleftarrow{\lim
}M_{n}\right)  $ is a Fr\'{e}chet vector bundle.\\
\end{theorem}

\begin{definition}\label{D_FrechetBundleAtlas}  A set $\left\{ U^\alpha= \underleftarrow{\lim}\left(  U_{n}^{\alpha}, \tau\alpha=\tau_{n}^{\alpha}\right)  \right\}
_{n\in\mathbb{N}^{\ast},\alpha\in A}$ such that  
the following diagrams are
commutative:%
\[%
\begin{array}
[c]{ccc}%
\left(  \pi_{i}\right)  ^{-1}\left(  U_{i}^{\alpha}\right)  & \underleftarrow
{\lambda_{i}^{j}} & \left(  \pi_{j}\right)  ^{-1}\left(  U_{j}^{\alpha}\right)
\\
\tau_{i}^{\alpha}\downarrow &  & \downarrow\tau_{j}^{\alpha}\\
U_{i}^{\alpha}\times\mathbb{E}_{i} & \underleftarrow{\mu_{i}^{j}%
\times\overline{\lambda_{i}^{j}}} & U_{j}^{\alpha}\times\mathbb{E}_{j}%
\end{array}
\]
  and $\left\{U_n^\alpha\right\}_{\alpha\in A}$ is a covering of $M_n$ for all  $n\in \mathbb{N}^*$ will be called a Fr\'echet bundle atlas.\\
\end{definition}

\subsection{Projective limits of principal bundles}
\label{__ProjectiveLimitsOfPrincipalBundles}

\begin{definition}
$\mathcal{P}=\left\{  \left(  P_{n},\pi_{n},M_{n},G_{n}\right)  \right\}
_{n\in\mathbb{N}^{\ast}}$ is called a projective sequence of Banach principal
bundles if the following conditions are satisfied:

\begin{description}
\item[\textbf{(PSBPB 1)}] $\left(  M_{i},\mu_{i}^{j}\right)  _{\left(
i,j\right)  \in\mathbb{N}^{\ast}\times\mathbb{N}^{\ast},\ i\leq j\ }$ is a
projective sequence of Banach $C^{\infty}$-manifolds;

\item[\textbf{(PSBPB 2)}] $\left(  P_{i},\lambda_{i}^{j}\right)  _{\left(
i,j\right)  \in\mathbb{N}^{\ast}\times\mathbb{N}^{\ast},\ i\leq j\ }$ is a
projective sequence of Banach $C^{\infty}$-manifolds;

\item[\textbf{(PSBPB 3)}] $\left(  G_{i},\gamma_{i}^{j}\right)  _{\left(
i,j\right)  \in\mathbb{N}^{\ast}\times\mathbb{N}^{\ast},\ i\leq j\ }$ is a
projective sequence of Banach Lie groups;

\item[\textbf{(PSBVB 4)}] For all $i\leq j$, $\left(  \mu_{i}^{j},\lambda
_{i}^{j},\gamma_{i}^{j}\right)  $ is a principal bundle morphism$;$\newline

\item[\textbf{(PSBVB 5)}] $\left(  M_{n}\right)  _{n\in\mathbb{N}^{\ast}}$ has
the projective limit chart property \emph{(PLCL)} at any point relatively to
$(U=\underleftarrow{\lim}U_{n},\phi=\underleftarrow{\lim}\phi_{n})$;

\item[\textbf{(PSBVB 6)}] For any $x\in M=\underleftarrow{\lim}M_{n}$, there
exists a local trivialization $\tau_{n}:\pi_{n}^{-1}\left(  U_{n}\right)
\longrightarrow U_{n}\times G_{n}$ such that the following diagram is
commutative:%
\[%
\begin{array}
[c]{ccc}%
\left(  \pi_{i}\right)  ^{-1}\left(  U_{i}\right)  & \underleftarrow
{\lambda_{i}^{j}} & \left(  \pi_{j}\right)  ^{-1}\left(  U_{j}\right) \\
\tau_{i}\downarrow &  & \downarrow\tau_{j}\\
U_{i}\times G_{i} & \underleftarrow{\mu_{i}^{j}\times\gamma_{i}^{j}} &
U_{j}\times G_{j}%
\end{array}
\]

\end{description}
\end{definition}

We then have the following result.

\begin{theorem}
\label{T_FrechetStructureOnProjectiveLimitOfBanachPrincipalBundles} Let
$\left\{  \left(  P_{n},\pi_{n},M_{n},G_{n}\right)  \right\}  _{n\in
\mathbb{N}^{\ast}}$ be a projective sequence of Banach principal bundles. The
quadruple $\left(  \underleftarrow{\lim}P_{n},\underleftarrow{\lim}\pi
_{n},\underleftarrow{\lim}M_{n},\underleftarrow{\lim}G_{n}\right)  $ is a
Fr\'{e}chet principal bundle.
\end{theorem}

\subsection{Projective limits of generalized frame bundles}
\label{__ProjectiveLimitsOfGeneralizedFrameBundles}

We extend the notion developed in \cite{DGV}, 6.5 to the framework of
generalized frame bundles defined over a projective sequence of Banach manifolds.

Let us consider a projective sequence of Banach vector bundles 
$\mathcal{E}=\left\{  \left(  E_{n},\pi_{n},M_{n}\right)  \right\}  _{n\in\mathbb{N}^{\ast}}$ (cf. Definition \ref{D_ProjectiveSequenceOfBanachVectorBundles})
where the manifold $M_{n}$ is modelled on the Banach space $\mathbb{M}_{n}$
and where $\mathbb{E}_{n}$ is the fibre type of $E_{n}$. The projective limit
of $\mathcal{E}$ is a Fr\'{e}chet vector bundle $(E,\pi, M)$ (cf. Theorem
\ref{T_FrechetStructureOnProjectiveLimitOfBanachVectorBundles}).

We then define the space
\[
\mathbf{P}\left(  E_{n}\right)  =\bigcup\limits_{x_{n}\in M_{n}}%
\mathcal{H}_{0}^{n}\left(  \mathbb{E}\text{,}\left(  E_{n}\right)  _{x_{n}%
}\right)
\]
where  $\mathcal{H}_{0}^{n}\left(  \mathbb{E}\text{,}\left(
E_{n}\right)  _{x_{n}}\right)  $ is
\[
\left\{  \left(  q_{1},\dots,q_{n}\right)  \in\prod\limits_{i=1}
^{n}\mathcal{L}is\left(  \mathbb{E}_{i},\left(  E_{i}\right)  _{\mu_i^n(x_n)}\right)
:\forall\left(  i,j\right)  \in\left\{  1,\dots,n\right\}  ^{2},i\leq
j,\ \lambda_{i}^{j}\circ q_{j}=q_{i}\circ\overline{\lambda_{i}^{j}}\right\}  .
\]

We denote by $\mathbf{p}_n:\mathbf{P}\left(  E_{n}\right)\longrightarrow M_n$ the natural projection. Then we have:

\begin{proposition}
\label{P_GeneralizedLinearFrameBundleOverMn}
The quadruple $\left(  \mathbf{P}\left(  E_{n}\right)  ,\mathbf{p}_{n},M_{n}
,\mathcal{H}_{0}^{n}\left(  \mathbb{E}\right)  \right)  $ is a Banach
principal bundle over $M_{n}$ with structural group $\mathcal{H}_{0}^{n}\left(  \mathbb{E}\right)  $.
\end{proposition}

\begin{proof} 
Consider a Fr\'echet bundle atlas $\left\{\left(U^\alpha= \underleftarrow{\lim}U^\alpha_i,\tau^\alpha= \underleftarrow{\lim}\tau^\alpha_i\right)\right\}_{\alpha\in A}$ for the Fr\'echet bundle $(E,\pi, M)$. For each $i\in \mathbb{N}^*$,  this bundle atlas gives rise to a a bundle atlas $\left(U^\alpha_i,\Phi^\alpha_i\right)_{\alpha\in A}$ for the principal bundle $\ell\left( E_i\right)$ and the transition maps $T^{\alpha\beta}_i:U_i^\alpha\cap U_i^\beta\longrightarrow \operatorname{GL}(\mathbb{E}_{i})$ are the same for both these bundles (cf. \ref{P_CoincidenceTransitionFunctions-E-PE}). \\
We consider the map 
\[
\begin{matrix}
  \mathbf{\Phi}^\alpha_n: & \left(\mathbf{p}_n\right)^{-1}(U_n^\alpha) & \longrightarrow & U_n^\alpha\times \prod\limits_{i=1}^{n}\operatorname{GL}(\mathbb{E}_{n}) \\ 
   & \left(x_n;  q_{1},\dots,q_{n}\right) & \mapsto & \left(x_n; (\bar{\tau}^\alpha_1)_{\mu_1^n(x_n)}\circ q_1,\dots,(\bar{\tau}^\alpha_1)_{\mu_n^n(x_n)}\circ q_n\right)
  \end{matrix}   
\]
where $(\tau^\alpha_i)_{x_i}=pr_2\circ (\tau^\alpha_i)_{| (E_n)_{x_i}}$ (for more details see $\S$ \ref{__FrameBundleOfABanachVectorBundle}).

Since $\mu_i^j$ and $\tau^\alpha_i$ are smooth maps, it follows that  $\mathbf{\Phi}^\alpha_n $ is smooth. Now, according to property (PSBVB 5) and the definition of $\mathcal{H}_{0}^{n}\left(  \mathbb{E}\text{,}\left(  E_{n}\right)  _{x_{n}}\right)$, the map $\mathbf{\Phi}^\alpha_n$ takes values in $\mathcal{H}_0^n(\mathbb{E})$.\\
  
Now the map
\[
\left(x_n;g_1,\dots g_n\right)\mapsto \left(x_n;\{(\bar{\tau}^\alpha_1)\}_{\mu_1^n(x_n)}^{-1}(g_1),\dots \{(\bar{\tau}^\alpha_n)\}_{\mu_n^n(x_n)}^{-1}(g_n)\right)
\] 
is the inverse of  $\mathbf{\Phi}^\alpha_n$ and this map is also smooth. It follows that   $\mathbf{\tau}_n^\alpha$ is a trivialization. \\

The transition maps $\mathbf{T}^{\alpha\beta}_n$ are given by
\[
\mathbf{T}^{\alpha\beta}_n (x_n)=\left(T^{\alpha\beta}_1(\mu_1^n(x_n)),\dots, T^{\alpha\beta}_n(\mu_n^n(x_n))\right).
\]
 
Again, property (PSBVB 5) implies that, for all $1\leq i\leq j$,  
\[ 
 T^{\alpha\beta}_i(\mu_i^j(x_j);\overline{\lambda_i^j}(u))=\overline{\lambda_i^j}\circ T^{\alpha\beta}_j(x_j;u)
\]
Thus  $\mathbf{T}^{\alpha\beta}_n (x_n)$ takes values in $\mathcal{H}_0^n(\mathbb{E})$ which ends the proof.\\
\end{proof}

\begin{definition}
\label{D_GeneralizedFrameBundle}
The quadruple $\left(  \mathbf{P}\left(E_{n}\right)  ,\mathbf{p}_{n},M_{n},\mathcal{H}_{0}^{n}\left(  \mathbb{E}\right)\right)  $ is called the 
generalized frame bundle of $E_{n}$ and is
denoted by $\mathbf{\ell}\left(  E_{n}\right)  $.
\end{definition}

%\begin{definition}
%$\left\{  \mathbf{\ell}\left(  E_{n}\right)  \right\}  _{n\in\mathbb{N}^{\ast
%}}=\left\{  \left(  \mathbf{P}\left(  E_{n}\right)  ,\mathbf{p}_{n},M_{n}%
%,\mathcal{H}_{0}^{n}\left(  \mathbb{E}\right)  \right)  \right\}
%_{n\in\mathbb{N}^{\ast}}$ is called a projective sequence of generalized frame
%bundles if $\left\{  \left(  \mathbf{P}\left(  E_{n}\right)  ,\mathbf{p}_{n}%
%,M_{n},\mathcal{H}_{0}^{n}\left(  \mathbb{E}\right)  \right)  \right\}
%_{n\in\mathbb{N}^{\ast}}$is a projective sequence of principal bundles.
%\end{definition}

For every $i\leq j$, we define the following projection 
\[
\begin{array}
[c]{cccc}%
\mathbf{r}_{i}^{j}: & \mathbf{P}\left(  E_{j}\right)  & \longrightarrow &
\mathbf{P}\left(  E_{i}\right) \\
& \left(  q_{1},\dots,q_{j}\right)  & \mapsto & \left(  q_{1},\dots
,q_{i}\right).
\end{array}
\]
Then, as for \cite{DGV}, Lemma 6.5.3, from the previous results and definitions  we obtain 

\begin{lemma}
\label{L_ProjectiveLimitOfSequenceGeneralizedFrameBundles}
The triple $\left(  \mathbf{r}_{i}^{j},\mu_{i}^{j},\left(  h_{0}\right)  _{i}^{j}\right)
$ is a principal bundle morphism of $\left(  \mathbf{P}\left(  E_{j}\right)
,\mathbf{p}_{j},M_{j},\mathcal{H}_{0}^{j}\left(  \mathbb{E}\right)  \right)  $ into
$\left(  \mathbf{P}\left(  E_{j}\right)  ,\mathbf{p}_{i},M_{i},\mathcal{H}_{0}%
^{i}\left(  \mathbb{E}\right)  \right)  $.
\end{lemma}

According to the proof of Proposition \ref{P_GeneralizedLinearFrameBundleOverMn} and Lemma \ref{L_ProjectiveLimitOfSequenceGeneralizedFrameBundles} and adapting the proof of \cite{DGV}, Proposition 6.5.4, we obtain:

\begin{theorem}
\label{T_StructureOnProjectiveLimitOfSequenceGeneralizedFrameBundles} The sequence 
$\left(  \mathbf{\ell}\left(  E_{n}\right)  \right)  _{n\in\mathbb{N}^{\ast}%
}=\left(  \mathbf{P}\left(  E_{n}\right)  ,\mathbf{p}_{n},M_{n},\mathcal{H}_{0}%
^{n}\left(  \mathbb{E}\right)  \right)  _{n\in\mathbb{N}^{\ast}}$ is a
projective sequence of principal bundles.  The projective
limit $\underleftarrow{\lim}\mathbf{P}\left(  E_{n}\right)  $ can be endowed
with a structure of smooth Fr\'{e}chet principal bundle over $\underleftarrow
{\lim}M_{n}$ whose structural group is $\mathcal{H}_{0}\left(  \mathbb{E}%
\right)  $.
\end{theorem}

\begin{remark}
\label{R_TransitionMapsAtlas}  
Adapting \cite{DGV}, Corollary 6.5.2 in the context of Theorem \ref{T_StructureOnProjectiveLimitOfSequenceGeneralizedFrameBundles} and according to the proof of Proposition \ref{P_GeneralizedLinearFrameBundleOverMn}, we have a Fr\'echet principal bundle atlas 
$\left\{ \left( U^\alpha=\underleftarrow{\lim}U_n^\alpha, \mathbf{\Phi}^\alpha=\underleftarrow{\lim}\mathbf{\Phi}_n^\alpha \right)\right\}_{\alpha\in A}$ 
such that, for each $n\in \mathbb{N}^*$,  each $   \mathbf{\ell}\left(  E_{n}\right)$  is trivializable over each $U_n^\alpha$  and  the transition maps are  
\[
\mathbf{T}^{\alpha\beta}_n (x_n)=\left(T^{\alpha\beta}_1(\mu_1^n(x_n)),\dots, T^{\alpha\beta}_n(\mu_n^n(x_n))\right)
\]
 where $T^{\alpha\beta}_j$ are the transition functions of the Banach bundle $\left(E_i,\pi_n, M_i\right)$ 
for the atlas $\{\left(U^\alpha_i,\tau^\alpha_i\right)\}_{\alpha\in A}$ where $1\leq i\leq n$.
\end{remark}

\subsection{Projective limits of $G$-structures}
\label{__ProjectiveLimitsOfGStructures}

Let
$\left(\mathcal{G}_{0}^{i}\left(  \mathbb{E}\right) ,\gamma_i^j\right) _{\left(
i,j\right)  \in\mathbb{N}^{\ast}\times\mathbb{N}^{\ast},\ i\leq j\ }$  
be a projective sequence of Banach Lie groups associated to a sequence $(G_n)$ of weak Banach Lie subgroup $G_n$ of $\operatorname{GL}(\mathbb{E}_{n})$ (cf. $\S$ \ref{__FrechetTopologicalSubgroupsOfH0E}). 
 
\begin{definition}
\label{D_ProjectiveLimitsG-Structures}
A sequence  
$\left\{ \left(  F_{n},\mathbf{p}_{n}|F_{n},M_{n},\mathcal{G}_{0}^{n}\left(  \mathbb{E}\right) \right) \right\}
_{n\in\mathbb{N}^{\ast}}$ 
is called a projective sequence of G-reductions of a  projective sequence of Banach principal bundles $\left\{  \left(  E_{n},\mathbf{p}_{n},M_{n}, \mathcal{H}%
_{0}^{n}\left(  \mathbb{E}\right) \right)  \right\}  _{n\in
\mathbb{N}^{\ast}}$  if,  for each $n\in\mathbb{N}^*$, $\left(  F_{n},\mathbf{p}_{n|F_{n}},M_{n},\mathcal{G}_{0}^{n}\left(  \mathbb{E}\right)  \right)$ is a  topological principal subbundle of $\left(  \mathbf{P}\left(  E_{n}\right)  ,\mathbf{p}_{n},M_{n},\mathcal{H}%
_{0}^{n}\left(  \mathbb{E}\right)  \right)$
 such that $\left(  F,\pi
_{|F},M,G\right)  $ has its own smooth principal structure, and the inclusion
is smooth.
\end{definition}

Then we have the following result:

\begin{theorem}
\label{T_FrechetStructureOnProjectiveLimitOfDescendingSequenceOfGStructures}
Consider a sequence  $\left(  \mathbf{\ell}\left(  E_{n}\right)  \right)_{n\in\mathbb{N}%
^{\ast}}=\left(  \mathbf{P}\left(  E_{n}\right)  ,\mathbf{p}_{n},M_{n},\mathcal{H}%
_{0}^{n}\left(  \mathbb{E}\right)  \right)  _{n\in\mathbb{N}^{\ast}}$ be a
 of generalized frame bundles over projective sequence of Banach
manifolds.\\
Let  $\left(  F_{n},\mathbf{p}_{n|F_{n}},M_{n},\mathcal{G}_{0}^{n}\left(  \mathbb{E}\right)  \right)  _{n\in\mathbb{N}^{\ast}}$ be a   projective system of G-reduction of $\left(  \mathbf{\ell}\left(  E_{n}\right)  \right)  _{n\in\mathbb{N}%
^{\ast}}$. 

Then $\underleftarrow{\lim}\mathbf{P}_{n}\left(  E_{n}\right)  $ can be
endowed with a structure of Fr\'{e}chet principal bundle over $\underleftarrow
{\lim}M_{n}$ whose structural group is $\mathcal{G}_{0}\left(  \mathbb{E}%
\right)  $.
\end{theorem}

\begin{proof}(summarized proof) According to Remark \ref{R_TransitionMapsAtlas}  where we consider the atlas
 $\left\{\left(U^\alpha=\underleftarrow{\lim}U_n^\alpha, \mathbf{\Phi}^\alpha=\underleftarrow{\lim}\mathbf{\Phi}_n^\alpha\right)\right\}_{\alpha\in A}$, 
 since each transition function $T^{\alpha\beta}_i$ belongs to $\mathcal{G}_0^n\left(\mathbb{E}\right)$ thus the transition maps $\mathbf{T}^{\alpha\beta}=\underleftarrow{\lim}\mathbf{T}_n^\alpha$  associated to the atlas  $\left\{\left(U^\alpha, \mathbf{\Phi}^\alpha\right)\right\}_{\alpha\in A}$ belongs to $\underleftarrow{\lim}\mathcal{G}_0^n\left(\mathbb{E}\right)$.\\
\end{proof}

\subsection{Projective limits of tensor structures}
\label{__ProjectiveLimitsOfTensorStructures}

\subsubsection{Projective limits of tensor structures of type $\left(1,1\right)$}
\label{___ProjectiveLimitsOfTensorStructuresOfType1-1}

Recall that a tensor of type $\left(1,1\right)  $ on a Banach space $\mathbb{E}$ is also an endomorphism of $\mathbb{E}$.

\begin{definition} 
\label{D_ProjectiveSequenceOfEndomorphisms}${}$
\begin{enumerate}
\item[1.] 
Let $\left\{  \left(\mathbb{E}_{i},\overline{\lambda_{i}^{j}}\right)  \right\}_{\left(i,j\right) \in\mathbb{N}^{\ast}\times\mathbb{N}^{\ast},\ i\leq j\ }$ be a
projective sequence of Banach spaces. For any $n\in\mathbb{N}^{\ast}$, let $A_{n}:\mathbb{E}_{n}\longrightarrow\mathbb{E}_{n}$ be an endomorphism of the Banach space $\mathbb{E}_{n}$. 
A sequence $\left(  A_{n}\right)  _{n\in\mathbb{N}^{\ast}}$ is called a projective sequence of endomorphisms, or is coherent for short, if, for any integer $j\geqslant i>0$, it fulfils the coherence condition:
\[
A_{i}\circ\overline{\lambda_{i}^{j}}=\overline{\lambda_{i}^{j}}\circ A_{j}%
\]
\item[2.] 
Let $\left\{  \left(  E_{n},\pi_{n},M_{n}\right)  \right\}  _{n\in
\mathbb{N}^{\ast}}$ be a projective sequence of Banach vector bundles.  
A sequence $\left(  \mathcal{A}_{n}\right)  _{n\in\mathbb{N}^{\ast
}}$  endomorphisms $\mathcal{A}_{n}$ of  $E_n$ is called a projective sequence of endomorphisms, or a coherent sequence for short, if, for each $x=\underleftarrow{\lim}x_n$, the sequence $\left(\left(\mathcal{A}_n\right)_{x_n}\right)$ is a coherent sequence of endomorphism of $\left(E_n\right)_{x_n}$.\\
\end{enumerate}
\end{definition}

We have the following properties:

\begin{proposition}
\label{P_ProjectiveLimitOfEndomorphisms}
Consider   a sequence $\left(  \mathcal{A}_{n}\right)  _{n\in\mathbb{N}^{\ast
}}$ of coherent   endomorphisms  $\mathcal{A}_{n}$  of $E_n$ on   a projective sequence of Banach vector bundles  $\left\{  \left(  E_{n},\pi_{n},M_{n}\right)  \right\}  _{n\in
\mathbb{N}^{\ast}}$. Then the projective limit $\mathcal{A}=\underleftarrow{\lim}({\mathcal{A}}_n)$ is well defined and is a smooth endomorphism of the Fr\'echet bundle  $E=\underleftarrow{\lim}{E_n}$.
\end{proposition}

\begin{proof}
 By definition, a coherent sequence of endomorphisms is nothing but a projective sequence of linear maps; So the projective limit ${A}=\underleftarrow{\lim}{{A}}_n$ is a well defined endomorphism of $\mathbb{E}=\underleftarrow{\lim}{\mathbb{E}}_n$. It follows that, for each $x=\underleftarrow{\lim}{x}_n\in M=\underleftarrow{\lim}{M}_n$, the projective limit 
$\mathcal{A}_x=\underleftarrow{\lim}(\mathcal{A}_n)_{x_n}$ is well defined. Now from property (PSBVB 5), there
exists a local trivialization $\tau_{n}:\pi_{n}^{-1}\left(  U_{n}\right)
\longrightarrow U_{n}\times\mathbb{E}_{n}$ such that the following diagram is
commutative:
\[
\begin{array}
[c]{ccc}
\left(  \pi_{i}\right)  ^{-1}\left(  U_{i}\right)  & \underleftarrow
{\lambda_{i}^{j}} & \left(  \pi_{j}\right)  ^{-1}\left(  U_{j}\right) \\
\tau_{i}\downarrow &  & \downarrow\tau_{j}\\
U_{i}\times\mathbb{E}_{i} & \underleftarrow{\mu_{i}^{j}\times\overline
{\lambda_{i}^{j}}} & U_{j}\times\mathbb{E}_{j}
\end{array}
\]
This implies that the restriction of $\mathcal{A}_n$ to $\left(  \pi_{i}\right)  ^{-1}\left(  U_{n}\right)$ is a projective sequence of smooth maps and so, from \cite{DGV}, Proposition 2.3.12, the restriction of $\mathcal{A}$ to $\pi  ^{-1}\left(  U\right)$ is a smooth map  where $U=\underleftarrow{\lim}U_n$, which ends the proof.\\
\end{proof} 

Let $\left\{  \left(  E_{n},\pi_{n},M_{n}\right)  \right\}  _{n\in
\mathbb{N}^{\ast}}$be a projective sequence of Banach vector bundles.  Consider a coherent  sequence $\left(  \mathcal{A}_{n}\right)  _{n\in\mathbb{N}^{\ast
}}$ of endomorphisms  $\mathcal{A}_n$ defined  on 
the  Banach bundle $\pi_{n}:E_{n}\longrightarrow M_{n}.$  \\
Fix  any  $x^0=\underleftarrow{\lim}{x_n^0}\in M=\underleftarrow{\lim}M_n$  and identify each typical fibre of $E_n$ with $(E_n)_{x_n^0}$. For all $n\in \mathbb{N}^*$, we denote by $G_n$ the isotropy group of $(\mathcal{A}_n)_{x_n^0}$ and by $\mathcal{G}_0^n$ the associated weak Banach Lie subgroup of $\mathcal{H}_0^n$ (cf. $\S$ \ref{__FrechetTopologicalSubgroupsOfH0E}).

Then from according to  Theorem
\ref{T_FrechetStructureOnProjectiveLimitOfDescendingSequenceOfGStructures}, we have:

\begin{theorem}
\label{T_ReductionProjectiveSequenceTensor1-1}
 Let $\left\{  \left(  E_{n},\pi_{n},M_{n}\right)  \right\}  _{n\in
\mathbb{N}^{\ast}}$ be a projective sequence of Banach vector bundles.  Consider a coherent  sequence $\left(  \mathcal{A}_{n}\right)  _{n\in\mathbb{N}^{\ast
}}$ of endomorphisms  $\mathcal{A}_n$ defined  on 
the  Banach bundle $\pi_{n}:E_{n}\longrightarrow M_{n}.$

Then $\underleftarrow{\lim}\mathbf{P}_{n}\left(  E_{n}\right)  $ can be
endowed with a structure of Fr\'{e}chet principal bundle over $\underleftarrow
{\lim}M_{n}$ whose structural group is $\mathcal{G}_{0}\left(  \mathbb{E}%
\right)  $ if and only if there exists a Fr\'echet atlas bundle
   $\left\{ U^\alpha= \underleftarrow{\lim} U_{n}^{\alpha}, \tau^{\alpha}=\tau_{n}^{\alpha}  \right\}
_{n\in\mathbb{N}^{\ast},\alpha\in A}$   such that  $\mathcal{A}_{n}^{\alpha}=\mathcal{A}_{n|U_{n}^{\alpha}}$ is locally modelled  $(\mathcal{A}_n)_{x_n^0}$
 and  the transition maps $T_{n}^{\alpha\beta}\left(x_{n}\right)  $ take values in  to the isotropy group $G_{n}$, for  all $n\in \mathbb{N}^*$ and all $\alpha,\beta \in A$ such that $U_n^\alpha\cap U_n^\beta\not=\emptyset$ .\\   
\end{theorem}

\begin{proof}  
According to Proposition \ref{P_ConditionTensorStructure}, for each $n\in \mathbb{N}^*$,  the Banach  principal bundle $\ell(E_n)$ can be
endowed with a structure of Banach principal bundle over $M_{n}$ whose structural group is $G_n$ if and only if there exists an atlas bundle $\left(U_n^\alpha,\tau_n^\alpha\right)_{\alpha\in A}$ such that $\mathcal{A}_{n}^{\alpha}$ is locally 
modelled  $(\mathcal{A}_n)_{x_n^0}$ and  the transition functions $T_{\alpha\beta
}\left(  x\right)  $ belong to the isotropy group $G_n$  for all  $x_n\in U_{\alpha}\cap U_{\beta}$. Now, assume that we have Fr\'echet atlas $\left\{ U^\alpha= \underleftarrow{\lim}\left(  U_{n}^{\alpha}, \tau^{\alpha}=\tau_{n}^{\alpha}\right)  \right\}
_{n\in\mathbb{N}^{\ast},\alpha\in A}$  for which such a property is true for all $n\in \mathbb{N}^*$.  For all $1\leq i\leq j$ we have

$\tau^\alpha_i\circ \lambda_i^j=\mu_i^j\times \overline{\lambda_i^j}\circ \tau_j^\alpha$ for all $\alpha\in A$

Thus we obtain

$$T^{\alpha\beta}_i\circ \overline{\lambda_i^j}=\overline{\lambda_i^j}\circ T^{\alpha\beta}_j$$ and so $(T^{\alpha\beta}_1,\dots,T^{\alpha\beta}_n)$ belongs to $\mathcal{G}_0^n\left(\mathbb{E}\right)$ which implies that $\underleftarrow{\lim}\mathbf{P}_{n}\left(  E_{n}\right)  $ can be
endowed with a structure of Fr\'{e}chet principal bundle over $\underleftarrow
{\lim}M_{n}$ whose structural group is $\mathcal{G}_{0}\left(  \mathbb{E}%
\right)  $  according to section \ref{__ProjectiveLimitsOfGStructures}.\\

Conversely, under such an assumption, according to  Theorem
\ref{T_FrechetStructureOnProjectiveLimitOfDescendingSequenceOfGStructures}, consider a projective sequence $\left(  F_{n},\mathbf{p}_{n|F_{n}},M_{n},\mathcal{G}_{0}^{n}\left(  \mathbb{E}\right)  \right)  _{n\in\mathbb{N}^{\ast}}$ of G-reductions for  $\left(\ell\left(\mathbb{E}_n\right)\right)_{n\in\mathbb{N}^*}$.  From the proof of this theorem and Remark \ref{R_TransitionMapsAtlas}, each transition map  $T^{\alpha\beta}_n$ takes values in $G_n$ and so the the frame bundle of $(E_n,\pi_n,M_n)$ can be
endowed with a structure of Banach principal bundle over $M_{n}$ whose structural group is $G_n$ (cf. $\S$ \ref{__GStructureAndTensorStructuresOnABanachVectorBundle}), which ends the proof.\\
\end{proof}

\subsubsection{Projective limits of tensor structures of type $\left(2,0\right)  $}
\label{___ProjectiveLimitsOfTensorStructuresOfType2-0}

This section is formally analogue to the previous one and we only give an adequate Definition and results without any proof.\\

\begin{definition}${}$
\label{D_ProjectiveSequenceOfBilinearForms}
\begin{enumerate}
\item[1.] Let $\left\{  \left(\mathbb{E}_{i},\overline{\lambda_{i}^{j}}\right)  \right\}  _{\left(
i,j\right)  \in\mathbb{N}^{\ast}\times\mathbb{N}^{\ast},\ i\leq j\ }$ be a
projective sequence of Banach spaces.

For any $n\in\mathbb{N}^{\ast}$, let $\omega_{n}:\mathbb{E}_{n}\times
\mathbb{E}_{n}\longrightarrow\mathbb{R}$ be a bilinear form of the Banach
space $\mathbb{E}_{n}.$\\
$\left(  \omega_{n}\right)_{n\in\mathbb{N}^{\ast}}$ is called a projective sequence of $2$-forms if it
fulfils the coherence condition 
\[
\forall (i,j):j \geq i \geq 1, \; \omega_{j}=\omega_{i}\circ\left(  \overline{\lambda_{i}^{j}}\times
\overline{\lambda_{i}^{j}}\right).
\]

\item[2.] Let $\left\{  \left(  E_{n},\pi_{n},M_{n}\right)  \right\}  _{n\in
\mathbb{N}^{\ast}}$ be a projective sequence of Banach vector bundles.  
A sequence $\left(  \Omega_{n}\right)  _{n\in\mathbb{N}^{\ast
}}$ of bilinear forms $\Omega_{n}$ on  $E_n$ is called a projective sequence of bilinear forms, or is a coherent sequence for short, if for each $x=\underleftarrow{\lim}x_n$, the sequence $\left(\left(\Omega_n\right)_{x_n}\right)$ is a coherent sequence of bilinear forms of $\left(E_n\right)_{x_n}$.
\end{enumerate}
\end{definition}

We have the following properties:\\

\begin{proposition}
\label{P_ProjectiveLimitOfBilinearForms}
Consider   a sequence $\left(  \mathcal{A}_{n}\right)  _{n\in\mathbb{N}^{\ast}}$ of coherent bilinear forms  
$\Omega_{n}$  of $E_n$ on   a projective sequence of Banach vector bundles  $\left\{  \left(  E_{n},\pi_{n},M_{n}\right)  \right\}  _{n\in\mathbb{N}^{\ast}}$.\\
Then the projective limit $\omega=\underleftarrow{\lim}(\Omega_n)$ is well defined and is a smooth endomorphism of the Fr\'echet bundle  $E=\underleftarrow{\lim}{E_n}$.
\end{proposition}

Let $\left\{  \left(  E_{n},\pi_{n},M_{n}\right)  \right\}  _{n\in
\mathbb{N}^{\ast}}$ be a projective sequence of Banach vector bundles.  Consider a coherent  sequence $\left( \Omega_n\right)  _{n\in\mathbb{N}^{\ast
}}$ of bilinear forms  $\Omega_n$ defined  on 
the  Banach bundle $\pi_{n}:E_{n}\longrightarrow M_{n}.$  \\
Fix  any  $x^0=\underleftarrow{\lim}{x_n^0}\in M=\underleftarrow{\lim}M_n$  and identify each typical fiber of $E_n$ with $(E_n)_{x_n^0}$. For all $n\in \mathbb{N}^*$, we denote by $G_n$ the isotropy group of $(\Omega_n)_{x_n^0}$ and by $\mathcal{G}_0^n$ the associated weak Banach Lie subgroup of $\mathcal{H}_0^n$ (cf. $\S$ \ref{__FrechetTopologicalSubgroupsOfH0E}).

\begin{theorem}
\label{T_ReductionAscendingSequenceTensor0-2}
 Let $\left\{  \left(  E_{n},\pi_{n},M_{n}\right)  \right\}  _{n\in
\mathbb{N}^{\ast}}$ be a projective sequence of Banach vector bundles.  Consider a coherent  sequence $\left(  \Omega_{n}\right)  _{n\in\mathbb{N}^{\ast}}$ of bilinear forms  $\Omega_n$ defined  on 
the  Banach bundle $\pi_{n}:E_{n}\longrightarrow M_{n}.$\\
Then $\underleftarrow{\lim}\mathbf{P}_{n}\left(  E_{n}\right)  $ can be
endowed with a structure of Fr\'{e}chet principal bundle over $\underleftarrow
{\lim}M_{n}$ whose structural group is $\mathcal{G}_{0}\left(  \mathbb{E}%
\right)  $ if and only if there exists a Fr\'echet bundle atlas $\left\{ U^\alpha= \underleftarrow{\lim} U_{n}^{\alpha}, \tau^{\alpha}=\tau_{n}^{\alpha}  \right\}
_{n\in\mathbb{N}^{\ast},\alpha\in A}$   such that  $\Omega_{n}^{\alpha}=\Omega_{n|U_{n}^{\alpha}}$ is locally 
modelled on $(\mathcal{A}_n)_{x_n^0}$ and  the transition maps $T_{n}^{\alpha\beta}\left(x_{n}\right)  $ take values in the isotropy group $G_{n}$, for  all $n\in \mathbb{N}^*$ and all $\alpha,\beta \in A$ such that $U_n^\alpha\cap U_n^\beta\not=\emptyset$ .   
\end{theorem}

\subsection{Examples of projective limits of tensor structures}
\label{__ExampleOfProjectiveLimitsOfTensorStructures}

\subsubsection{Projective limit of compatible weak almost K\"ahler and almost para-K\"ahler structures}
\label{___ProjectiveLimitOfCompatibleWeakAlmostKalherAndAlmostParaKalherStructures}

We consider the context and notations of $\S$ \ref{__CompatibilityBetwenDifferentStructures}\\
If $\Omega$  (resp. $g$,  resp. $\mathcal{I}$) is  a weak symplectic form  (resp. a weak Riemannian metric, resp. an almost complex structure) on a Fr\'echet bundle $(E,\pi_{E},M)$, the compatibility of any pair of the data $(\Omega, g,\mathcal{I})$ as given in Definition \ref{D_CompatibilityWeakSymplecticFormAlmostComplexStructureWeakRiemannianMetricOnABanachBundle} can be easily adapted to this Fr\'echet context. We can also easily  transpose the notion of  "local tensor $\mathcal{T}$ locally modelled on  a linear tensor $\mathbb{T}$" (cf. Definition \ref{D_TensorLocallyModelled}) in the Fr\'echet framework. We will see that such situations can be obtained by "projective limit"
 on a  sequence $\left\{  \left(  E_{n},\pi_{n},M_{n}\right)  \right\}  _{n\in
\mathbb{N}^{\ast}}$ of Banach vector bundles.  \\
 
Now, we introduce the following notations:
\begin{enumerate}
\item[-- ] 
Let $(\mathbb{T}_n)$ be a sequence of coherent tensors of type $(1,1)$ or $(2,0)$ on a projective sequence of Banach spaces $\left\{  \left(\mathbb{E}_{i},\overline{\lambda_{i}^{j}}\right)  \right\}  _{\left(
i,j\right)  \in\mathbb{N}^{\ast}\times\mathbb{N}^{\ast},\ i\leq j\ }$.\\
We denote by $G_n(\mathbb{T}_n)$ the isotropy group of $\mathbb{T}_n$. If $\mathbb{T}=\underleftarrow{\lim}\mathbb{T}_n$, then   $\mathcal{G}_0^n(\mathbb{T})$ is the Banach Lie group associated to the sequence $\left(G_n(\mathbb{T}_n)\right)_{n\in \mathbb{N}^*}$ (see section \ref{__FrechetTopologicalSubgroupsOfH0E}). We set $\mathcal{G}_0(\mathbb{T})=\underleftarrow{\lim}\mathcal{G}_0^n(\mathbb{T})$.
\item[--]  
Given  a sequence $\left\{  \left(  E_{n},\pi_{n},M_{n}\right)  \right\}  _{n\in
\mathbb{N}^{\ast}}$ of Banach vector bundles, we provide 
 each Banach bundle $E_n$ with a symplectic form $\Omega_n$ and/or  a weak Riemaniann metric $g_n$ and/or  an  complex structure 
 $\mathcal{I}_n$.
\end{enumerate}
 
According to this context, let us consider the following assumption:

\begin{description}
\item[\textbf{(H)}]
Assume that we have a sequence of weak Riemaniann metrics $(g_n)_{n\in \mathbb{N}^*}$ on the sequence $(E_n)_{n\in \mathbb{N}^*}$.  
For each $x_n\in M_n$, we denote by $(\widehat{E}_n)_{x_n}$ the completion of $(E_n)_{x_n}$ and we set $\widehat{E}_n=\bigcup\limits_{x_n\in M_n}(\widehat{E}_n)_{x_n}$ and $\widehat{E}=\bigcup\limits_{x_n\in M_n}  (\widehat{E}_n)_{x_n}$. Then  
 $\widehat{\pi}_n:\widehat{E}_n\longrightarrow M_n$ has a  Hilbert vector bundle structure and the inclusion $\iota_n:E_n\rightarrow \widehat{E}_n$ is a Banach bundle morphism.
\end{description}
 
From Theorem \ref{T_CompatibilityDarbouxFormWeakRiemannianMetricComplexStructureOnABanachBundle}, we then  obtain:
 
\begin{theorem}
\label{T_ ProjectiveLimitOfCompatibleWeakSymplecticFormsRiemannianMetricsComplexStructures}${}$
 \begin{enumerate}
\item[1.] 
Under the assumption (H),  $\left\{\left(\widehat{E}_n,\widehat{\pi}_n, M_n\right)\right\}_{n\in \mathbb{N}^*}$ is a projective sequence of Banach bundles and there exists an injective Fr\'echet bundle morphism $\iota: E=\underleftarrow{\lim}E_n\rightarrow \widehat{E}=\underleftarrow{\lim}\widehat{E}_n$ with dense range.\\ 
Let $\widehat{\Omega}_n$, (resp.$\widehat{g}_n$, resp. $\widehat{\mathcal{I}}_n$)  of ${\Omega}_n$ (resp. ${g}_n$, resp.  ${\mathcal{I}}_n$) to $\widehat{E}_n$ (cf. Theorem  \ref{T_CompatibilityDarbouxFormWeakRiemannianMetricComplexStructureOnABanachBundle}).
If the sequence  $({\Omega}_n)_{n\in \mathbb{N}^*}$, (resp. $({g}_n)_{n\in \mathbb{N}^*}$, resp. $(\mathcal{I}_n)_{n\in \mathbb{N}^*}$) is coherent, so is the sequence $(\widehat{\Omega}_n)_{n\in \mathbb{N}^*}$ (resp. $(\widehat{g}_n)_{n\in\mathbb{N}^*}$, resp.$(\widehat{\mathcal{I}}_n)_{n\in \mathbb{N}^*}$). 
  We set $\widehat{\Omega}=\underleftarrow{\lim}\widehat{\Omega}_n$ (resp. $\widehat{g}=\underleftarrow{\lim}\widehat{g}_n$, resp.
  $\widehat{\mathcal{I}}=\underleftarrow{\lim}\widehat{\mathcal{I}}_n$).  We then have
 ${\Omega}=\underleftarrow{\lim}{\Omega}_n=i^*\widehat{\Omega}$, $g=\underleftarrow{\lim}{g}_n=\iota^*\widehat{g}$ and 
 $\mathcal{I}=\underleftarrow{\lim}{\mathcal{I}}_n=\iota^*\widehat{\mathcal{I}}$. Moreover,  if any pair of the  data $(\Omega_n, g_n, \mathcal{I}_n)$ are compatible for all $n\in \mathbb{N}$, then  $(\widehat{\Omega},\widehat{g}, \widehat{\mathcal{I}})$  and  $({\Omega},{g}, {\mathcal{I}})$ on the Fr\'echet vector bundle  $\widehat{E}$ and $E$ respectively are compatible structures. 
 \item[2.]   
Let $\left(  \mathbf{\ell}\left(  E_{n}\right)  \right)  _{n\in\mathbb{N}%
^{\ast}}=\left(  \mathbf{P}\left(  E_{n}\right)  ,\mathbf{p}_{n},M_{n},\mathcal{H}%
_{0}^{n}\left(  \mathbb{E}\right)  \right)  _{n\in\mathbb{N}^{\ast}}$ be the associated
sequence of generalized frame bundles. Denote by $\mathcal{T}_n$ anyone tensor of the triple $(\Omega_n,g_n,\mathcal{I}_n)$. Assume that $\mathcal{T}_n$  defines a $\mathbb{T}$-tensor structure on $E_n$ for each $n\in \mathbb{N}^*$.  Then $\underleftarrow{\lim}\mathbf{P}_{n}\left(  E_{n}\right)  $  %and $\underleftarrow{\lim}\mathbf{P}_{n}\left(  \widehat{E}_{n}\right)  $ 
is a Fr\'{e}chet principal bundle over $\underleftarrow
{\lim}M_{n}$ whose structural group is $\mathcal{G}_{0}\left(  \mathbb{T}
\right)$.\\
Moreover,  if the sequence  $({\Omega}_n,g_n, {\mathcal{I}}_n)_{n\in \mathbb{N}^*}$ is coherent and each triple is a tensor structure  $({\Omega}_n,g_n, {\mathcal{I}}_n) $ on $E_n$, then 
  $({\Omega},g, {\mathcal{I}})$  is also a tensor structure on $E$ which compatible if $({\Omega}_n,g_n, {\mathcal{I}}_n)$ fulfils this property.
\item[3.] 
Under the assumption (H), all the properties of Point 2 are also valid for the sequence $(\widehat{\Omega}_n,\widehat{g}_n, \widehat{\mathcal{I}}_n)_{n\in \mathbb{N}^*} $ relatively to $\widehat{E}$.
\end{enumerate}
\end{theorem}
 
\begin{proof}
 1.  According to our assumption (H), by density of $\iota_n(E_n)$ (identified with $E_n$) in
$\widehat{E}_n$ and using compatible bundle atlases for $\widehat{E}$ and $E$
(cf. proof of Theorem \ref{T_GStructureKreinMetric}), for each $x_n\in M_n$ and each $n\in \mathbb{N}^*$, there exists a local trivialization  $\widehat{\tau}_{n}:\widehat{\pi}_{n}^{-1}\left(  U_{n}\right)\longrightarrow U_{n}\times\widehat{\mathbb{E}}_{n}$ and ${\tau}_{n}:{\pi}_{n}^{-1}\left(  U_{n}\right)
\longrightarrow U_{n}\times{\mathbb{E}}_{n}$ around $x_n$ such that $\widehat{\tau}_{n}\circ \iota_n={\tau}_{n}$ where $\widehat{\tau}_{n}\circ \iota_n(\widehat{\pi}_{n}^{-1}\left(  U_{n}\right))$ is dense in $U_{n}\times\widehat{\mathbb{E}}_{n}$.
Now, we can choose such a trivialization so that that the property (PSBVB 5) is satisfied for the sequence $\left\{\left( E_n,\pi_n, M_n \right)\right\}_{n\in \mathbb{N}^*}$. But since, for each $y_n\in U_n$, the fibre $(E_n)_{y_n}$ is dense in $(\widehat{E}_n)_{y_n}$, and each fibre inclusion $(\iota_n)_{y_n}$ is bounded, the bounding map $(\bar{\lambda}_i^j)_{y_j}:(E_j)_{y_j}\longrightarrow (E_i)_{y_i}$ is closed, $(\bar{\lambda}_i^j)_{y_j}$ can be extended to a bounding map   $(\widehat{\bar{\lambda}}_i^j)_{y_j}:(\widehat{E}_j)_{y_j}\longrightarrow (\widehat{E}_i)_{y_i}.$ Thus (PSBVB 5) is satisfied for the sequence $\{\widehat{E}_n,\widehat{\pi}_n, M_n\}_{n\in \mathbb{N}^*}$. Thus $\{\widehat{E}_n,\widehat{\pi}_n, M_n\}_{n\in \mathbb{N}^*}$ is a projective sequence of Banach bundles and of course 
$$\iota=\underleftarrow{\lim}\iota_n:E=\underleftarrow{\lim}E_n\longrightarrow \widehat{E}=\underleftarrow{\lim}\widehat{E}_n$$
is an injective Fr\'echet bundle morphism whose range is dense in $\widehat{E}$.\\

2 and 3. With the notations of Point 2,  assume that  the sequence $(\mathcal{T}_n)_{n\in \mathbb{N}^*}$ is coherent. Denote by $\widehat{T}_n$ the extension of $\mathcal{T}_n$ to $\widehat{E}_n$ (which is well defined by Theorem
 \ref{T_CompatibilityDarbouxFormWeakRiemannianMetricComplexStructureOnABanachBundle}). Using the same arguments as the ones used for the bonding maps, it follows that $(\widehat{\mathcal{T}}_n)_{n\in \mathbb{N}^*}$ is  also  a coherent sequence. Thus $\mathcal{T}=\underleftarrow{\lim}\mathcal{T}_n$ and $\widehat{\mathcal{T}}=\underleftarrow{\lim}\widehat{\mathcal{T}}_n$  are well defined and by construction we have $\mathcal{T}=\iota^*\widehat{\mathcal{T}}$.
 Moreover assume that  each $\mathcal{T}_n$ defines a $\mathbb{T}_n$-structure on $E_n$. Then from  Theorem \ref{T_CompatibilityDarbouxFormWeakRiemannianMetricComplexStructureOnABanachBundle},  $\widehat{\mathcal{T}}_n$  is a $\widehat{\mathbb{T}}$-structure on $\widehat{E}_n$. Thus from Theorem \ref{T_ReductionProjectiveSequenceTensor1-1} or Theorem  \ref{T_ReductionAscendingSequenceTensor0-2}, each principal bundle $\mathbf{P}_{n}\left(  E_{n}\right)$ (resp. $\mathbf{P}_{n}\left(  \widehat{E}_{n}\right)  $) is well defined and we can reduce its structural group $\mathcal{G}_{0}^n\left(  \mathbb{T}
\right)$  (resp.  $\mathcal{G}^n_{0}\left( \widehat{ \mathbb{T}}\right)$). The last part of the proof is a direct consequence of the previous properties.
\end{proof}

%\begin{remark}
%\label{R_ProjectiveLimitParaKahler}
%Given  a sequence $\left\{ \left( E_n,\pi_n,M_n \right) \right\}_{n\in \mathbb{N}^*}$ of Banach vector bundles, we provide each Banach bundle $E_n$ with a symplectic form $\Omega_n$ and/or a neutral metric $g_n$ and/or  an  almost  para-complex structure $\mathcal{I}_n$. Then by formal analogue arguments, we can obtain a similar Theorem as Theorem \ref{T_ProjectiveLimitOfCompatibleWeakSymplecticFormsRiemannianMetricsComplexStructures}.
%\end{remark}

\begin{corollary} 
\label{C_ProjectiveLimitOfAlmostKahlerStructures}
 Consider an almost K\"ahler (resp. para-K\"ahler) structure $({\Omega}_n,g_n, {\mathcal{I}}_n)_{n\in \mathbb{N}^*}$ on a projective $\{E_n,\pi_n,M_n\}_{n\in \mathbb{N}^*}$ of Banach vector bundles. If $({\Omega}_n,g_n, {\mathcal{I}}_n)_{n\in \mathbb{N}^*}$ is coherent, then $({\Omega}=\underleftarrow{\lim}{\Omega}_n,g=\underleftarrow{\lim}{g}_n,\mathcal{I}=\underleftarrow{\lim}{\mathcal{I}}_n)$ is an almost K\"ahler (resp. para-K\"aler) structure on the Fr\'echet bundle $E=\underleftarrow{\lim} E_n$.
\end{corollary}

\subsubsection{Application to sets of smooth maps}
\label{___ApplicationToSetsOfSmoothMaps}
Let $M$ be connected manifold of dimension $m$ and $N$ a connected compact manifold of dimension $n$. We denote by $\mathsf{C}^{k}(N,M)$ the set of maps $f:N\longrightarrow M$ of class $C^k$. It is well known that  $\mathsf{C}^{k}(N,M)$ is a Banach manifold  modelled on the Banach space $\mathsf{C}^{k}
(N,\mathbb{R}^m)$ (cf. \cite{Eli}). If $\lambda_i^j: \mathsf{C}^{j}(N,M)\longrightarrow \mathsf{C}^{i}(N,M)$ is the natural projection, the sequence $\left(\mathsf{C}^{k}(N,M)\right)_{k\in\mathbb{N}}$ is a projective sequence of Banach manifolds. By the way, the projective limit $\mathsf{C}^{\infty}(N,M)=\underleftarrow{\lim}\mathsf{C}^{k}(N,M)$ is provided with a Fr\'echet structure which is the usual Fr\'echet manifold structure on $\mathsf{C}^{\infty}(N,M)$.

Assume that $M$ is   a K\"ahler manifold. Let $\omega$, $\mathfrak{g}$ and $\mathfrak{I}$ be the associated symplectic form, Riemannian metric and complex structure on $M$ respectively.

If $\nu$ is a volume measure on $N$, on the tangent bundle $T\mathsf{C}^{k}(N,M)$ we introduce  (cf.\cite{Kum1} and \cite{Kum2}):
\begin{enumerate}
\item[--] $(\Omega_k)_f(X,Y)=\displaystyle\int_{N}\omega\left(X_{f(t)},Y_{f(t)}\right)d\nu(t)$;
\item[--] $(g_k)_f(X,Y)=\displaystyle\int_{N}\mathfrak{g}\left(X_{f(t)},Y_{f(t)}\right)d\nu(t)$;
\item[--] $(\mathcal{I}_k)_f(X)(t)=\mathfrak{I}\left(X_{f(t)}\right)$.
\end{enumerate}

It is clear that $\Omega_k$ (resp. $g_k$) is a weak bilinear form (resp. a weak Riemaniann metric) on $\mathsf{C}^{k}(N,M)$ and $\mathcal{I}_k$ is an almost complex structure on $T\mathsf{C}^{k}(N,M)$. Moreover, the compatibility of the triple $(\omega, \mathfrak{g}, \mathfrak{I})$ clearly implies the compatibility of $(\Omega_k, g_k,\mathcal{I}_k)$; In fact $(\Omega_k, g_k,\mathcal{I}_k)_{k\in \mathbb{N}^*}$ is a coherent sequence. Then, from Corollary \ref{C_ProjectiveLimitOfAlmostKahlerStructures} we get an almost K\"ahler structure on $\mathsf{C}^{\infty}(N,M)$. Note that this structure which can be directly defined on $\mathsf{C}^{\infty}(N,M)$ by the same formulae.

{\it If $M$ is a para-K\"ahler manifold, we can provide $\mathsf{C}^{\infty}(N,M)$ in the same way.}

\begin{remark} 
If $M$ is a K\"ahler manifold, then the almost complex structure is integrable and, in particular, $N_\mathfrak{I}\equiv 0$. It is clear that that we also have $N_{\mathcal{I}_k}\equiv 0$   for all $k\in \mathbb{N}^*$ and $N_\mathcal{I}\equiv 0$. Thus $\mathcal{I}$ is formally integrable but not integrable in general (cf. \cite{Kum1} in the case of sets of paths). On the opposite, if $M$ is a para-K\"ahler manifold then again the Nijenhuis tensor of the almost para-complex structure $\mathcal{J}_k$  is zero and, by Theorem \ref{T_CaracterizationOfIntegrableAlmostTangentParaComplexDecomposableComplexStructure}, this implies that the almost para-complex structure $\mathcal{J}$ on $\mathsf{C}^{\infty}(N,M)$ is also integrable.\\
\end{remark}

\section{Direct limits of tensor structures}
\label{_DirectLimitsOfTensorStructures}

In order to endow the direct limit of an ascending sequence of frame bundles
$\left(  \ell\left(  E_{n}\right)  \right)  _{n\in\mathbb{N}^{\ast}}=\left(
P_{n}\left(  E_{n}\right)  ,\mathbf{p}_{n},M_{n},G_{n}\right)  _{n\in\mathbb{N}%
^{\ast}}$ with a structure of convenient principal bundle, we first consider
the situation on the base $\underrightarrow{\lim}M_{n}$ where the models are
supplemented Banach manifolds. In order to get a convenient structure on the
direct limit $\underrightarrow{\lim}\ell\left(  E_{n}\right)  $ we have to
replace the pathological general linear group $\operatorname{GL}\left(
\mathbb{E}\right)  $ where $\mathbb{E=}\underrightarrow{\lim}\mathbb{E}_{n}$
(where $\mathbb{E}_{n}$ is the typical fibre of the bundle $E_{n}$) by another
convenient Lie group $\operatorname*{G}\left(  \mathbb{E}\right)  $.

The main references for this section are \cite{Bou}, \cite{CaPe},
\cite{Dah}, \cite{Glo}, \cite{Nee} and \cite{RoRo}.

\subsection{Direct limits}
\label{__DirectLimits}

\begin{definition}
\label{D_DirectSystem}
Let $\left(  I,\leq\right)  $ be a directed set and $\mathbb{A}$ a category.\\ 
$\mathcal{S}=\left\{  \left(  Y_{i},\varepsilon_{i}^{j}\right)  \right\}  _{i\in I,\ j\in I,\ i\leq j}$ is called a direct system if 

\begin{enumerate}
\item[--]
For all $i \in I$, $Y_i$ is an object of the category;
\item[--] 
For all $(i,j)\in I^2:i \leq j, \; \varepsilon_{i}^{j}:Y_{i}\longrightarrow Y_{j}$ is a morphism (bonding map)
\end{enumerate} 

such that:

\begin{description}
\item[\textbf{(DS 1)}] 
$\forall i \in I, \; \varepsilon_{i}^{i}=\operatorname{Id}_{Y_{i}}$;

\item[\textbf{(DS 2)}] 
$\forall (i,j,k) \in I^3:i\leq j\leq k, \; \varepsilon_{j}^{k}\circ\varepsilon_{i}^{j}=\varepsilon_{i}^{k}$.
\end{description}
\end{definition}

When $I=\mathbb{N}$ with the usual order relation, countable direct systems
are called \textit{direct sequences}.

\begin{definition}
\label{D_Cone}
$\mathcal{S}=\left\{\left(  Y,\varepsilon
_{i}\right)\right\}_{i\in I}$ is called a cone if $Y\in\operatorname*{ob}\mathbb{A}$ and
$\varepsilon_{i}:Y_{i}\longrightarrow Y$ is such that $\varepsilon_{i}%
^{j}\circ\varepsilon_{i}=\varepsilon_{j}$ whenever $i\leq j$.
\end{definition}

\begin{definition}
\label{D_DirectLimit}A cone 
$\left\{\left(  Y,\varepsilon_{i}\right)\right\}_{i\in I}$  is a
\textit{direct limit} of $\mathcal{S}$ if for every cone $\left\{ \left(  Z,\theta
_{i}\right) \right\}_{i\in I}$ over $\mathcal{S}$ there exists a unique morphism
$\psi:Y\longrightarrow Z$ such that $\psi\circ\varepsilon_{i}=\theta_{i}$. We
then write $Y=\underrightarrow{\lim}\mathcal{S}$ or $Y=\underrightarrow{\lim
}Y_{i}$.
\end{definition}

\subsection{Ascending sequences of supplemented Banach spaces}
\label{__AscendingSequencesOfSupplementedBanachSpaces}

Let $\mathbb{M}_{1}^{\prime}$, $\mathbb{M}_{2}^{\prime}$, ... be Banach spaces
such that:%
\[
\left\{
\begin{array}
[c]{c}%
\mathbb{M}_{1}=\mathbb{M}_{1}^{\prime}\\
\forall n\in\mathbb{N}^{\ast},\ \mathbb{M}_{n+1}\simeq\mathbb{M}_{n}%
\times\mathbb{M}_{n+1}^{\prime}%
\end{array}
\right.
\]

For $i,j\in\mathbb{N}^{\ast},\;i<j$, let us consider the injections
\[%
\begin{array}
[c]{cccc}%
\iota_{i}^{j}: & \mathbb{M}_{i}\backsimeq\mathbb{M}_{1}^{\prime}\times
\cdots\times\mathbb{M}_{i}^{\prime} & \rightarrow & \mathbb{M}_{j}%
\backsimeq\mathbb{M}_{1}^{\prime}\times\cdots\times\mathbb{M}_{j}^{\prime}\\
& (x_{1}^{\prime},\dots,x_{i}^{\prime}) & \mapsto & (x_{1}^{\prime}%
,\dots,x_{i}^{\prime},0,\dots,0)
\end{array}
\]

\begin{definition}
\label{D_AscendingSequenceOfSupplementedBanachSpaces}$\left\{  (\mathbb{M}%
_{n},\iota_{n}^{n+1})\right\}  _{n\in\mathbb{N}^{\ast}}$ is called an
ascending sequence of supplemented Banach spaces.
\end{definition}

\subsection{Direct limits of ascending sequences of Banach manifolds}
\label{__DirectLimitsOfAscendingSequencesOfBanachManifolds}

\begin{definition}
\label{D_AscendingSequenceOfBanachManifolds}\bigskip$\mathcal{M}%
=(M_{n},\varepsilon_{n}^{n+1})_{n\in\mathbb{N}^{\ast}}$ is called an ascending
sequence of Banach manifolds if for any $n\in\mathbb{N}^{\ast}$, $\left(
M_{n},\varepsilon_{n}^{n+1}\right)  $ is a weak submanifold of $M_{n+1}$ where
the model $\mathbb{M}_{n}$ is supplemented in $\mathbb{M}_{n+1}.$
\end{definition}

\begin{proposition}
\label{P_ConditionsDefiningDLCP}Let $\mathcal{M}=(M_{n},\varepsilon_{n}%
^{n+1})_{n\in\mathbb{N}^{\ast}}$ be an ascending sequence of Banach
manifolds.\newline Assume that for $x\in M=\underrightarrow{\lim}M_{n}$, there
exists a family of charts $(U_{n},\phi_{n})$ of $M_{n}$, for each
$n\in\mathbb{N}^{\ast}$, such that:

\begin{description}
\item[\textbf{(ASC 1)}] $(U_{n})_{n\in\mathbb{N}^{\ast}}$ is an ascending
sequence of chart domains;

\item[\textbf{(ASC 2)}] $\phi_{n+1}\circ\varepsilon_{n}^{n+1}=\iota_{n}%
^{n+1}\circ\phi_{n}$.\newline Then $U=\underrightarrow{\lim}U_{n}$ is an open
set of $M$ endowed with the $\mathrm{DL}$-topology and $\phi=\underrightarrow
{\lim}\phi_{n}$ is a well defined map from $U$ to $\mathbb{M}=\underrightarrow
{\lim}\mathbb{M}_{n}$. \newline Moreover, $\phi$ is a continuous homeomorphism
from $U$ onto the open set $\phi(U)$ of $\mathbb{M}$.
\end{description}
\end{proposition}

\begin{definition}
\label{DLChartProperty} We say that an ascending sequence $\mathcal{M}%
=(M_{n},\varepsilon_{n}^{n+1})_{n\in\mathbb{N}^{\ast}}$ of Banach manifolds
has the direct limit chart property \emph{\textbf{(DLCP)}} if $(M_{n}%
)_{n\in\mathbb{N}^{\ast}}$ satisfies \emph{\textbf{(ASC 1) }}$et$
\emph{\textbf{(ASC 2)}}.
\end{definition}

We then have the following result proved in \cite{CaPe}.

\begin{theorem}
\label{T_LBManifold} Let $\mathcal{M}=\left\{  (M_{n},\varepsilon_{n}%
^{n+1})\right\}  _{n\in\mathbb{N}^{\ast}}$ be an ascending sequence of Banach
manifolds, modelled on the Banach spaces $\mathbb{M}_{n}$. Assume that
$(M_{n})_{n\in\mathbb{N}^{\ast}}$ has the direct limit chart property at each
point $x\in M=\underrightarrow{\lim}M_{n}$.\newline Then there is a unique non
necessarly Haussdorf convenient manifold structure on $M=\underrightarrow
{\lim}M_{n}$ modeled on the convenient space $\mathbb{M}=\underrightarrow
{\lim}\mathbb{M}_{n}$ endowed with the $\mathrm{DL}$-topology. \newline
\end{theorem}

Moreover, if each $M_{n}$ is paracompact, then $M=\underrightarrow{\lim}M_{n}$
is provided with a Hausdorff convenient manifold structure.\newline

\subsection{The Fr\'{e}chet topological group $\operatorname{G}\left(
\mathbb{E}\right)  $}

\label{__TheFrechetTopologicalGroupG(E)}

Let $(\mathbb{E}_{n},\iota_{n}^{n+1})_{n\in\mathbb{N}^{\ast}}$ be an ascending
sequence of supplemented Banach spaces where $\mathbb{E}=\underrightarrow
{\lim}\mathbb{E}_{n}$ can be endowed with a structure of convenient vector
space. The group $\operatorname{GL}\left(  \mathbb{E}\right)  $ does not admit
any reasonable topological structure. So we replace it by the convenient Lie
group $\operatorname{G}(\mathbb{E})$ defined as follows.

Let $\mathbb{E}_{1}\subset\mathbb{E}_{2}\subset\cdots$ be the direct sequence
of supplemented Banach spaces; so there exist Banach subspaces $\mathbb{E}%
_{1}^{\prime},\mathbb{E}_{2}^{\prime},\dots$ such that:
\[
\left\{
\begin{array}
[c]{c}%
\mathbb{E}_{1}=\mathbb{E}_{1}^{\prime},\\
\forall i\in\mathbb{N}^{\ast},\mathbb{E}_{i+1}\backsimeq\mathbb{E}_{i}%
\times\mathbb{E}_{i+1}^{\prime}%
\end{array}
\right.  \;
\]
For $i,j\in\mathbb{N}^{\ast},\;i\leq j$, we have the injection
\[%
\begin{array}
[c]{cccc}%
\iota_{i}^{j}: & \mathbb{E}_{i}\backsimeq\mathbb{E}_{1}^{\prime}\times
\cdots\times\mathbb{E}_{i}^{\prime} & \rightarrow & \mathbb{E}_{j}%
\backsimeq\mathbb{E}_{1}^{\prime}\times\cdots\times\mathbb{E}_{j}^{\prime}\\
& (x_{1}^{\prime},\dots,x_{i}^{\prime}) & \mapsto & (x_{1}^{\prime}%
,\dots,x_{i}^{\prime},0,\dots,0)
\end{array}
\]
Any $A_{n+1}\in \operatorname{GL}\left(  \mathbb{E}_{n+1}\right)  $ is represented by
$\left(
\begin{array}
[c]{cc}%
A_{n} & B_{n+1}\\
A_{n}^{\prime} & B_{n+1}^{\prime}%
\end{array}
\right)  $ where%
\[
A_{n}\in\mathcal{L}\left(  \mathbb{E}_{n},\mathbb{E}_{n}\right)
,\ A_{n}^{\prime}\in\mathcal{L}\left(  \mathbb{E}_{n},\mathbb{E}_{n+1}%
^{\prime}\right)  ,\ B_{n+1}\in\mathcal{L}\left(  \mathbb{E}_{n+1}^{\prime
},\mathbb{E}_{n}\right)  \text{ and }B_{n+1}^{\prime}\in\mathcal{L}\left(
\mathbb{E}_{n+1}^{\prime},\mathbb{E}_{n+1}^{\prime}\right)  \text{.}%
\]
The group
\[
\operatorname{GL}_{0}\left(  \mathbb{E}_{n+1}|\mathbb{E}_{n}\right)  =\left\{  A\in \operatorname{GL}\left(
\mathbb{E}_{n+1}\right)  :A\left(  \mathbb{E}_{n}\right)  =\mathbb{E}%
_{n}\right\}
\]
can be identified with the Banach-Lie sub-group of operators of type $\left(
\begin{array}
[c]{cc}%
A_{n} & B_{n+1}\\
0 & B_{n+1}^{\prime}%
\end{array}
\right)  $ (cf. \cite{ChSt}).

The set
\[
\operatorname{G}(E_{n})=\left\{  A_{n}\in \operatorname{GL}(\mathbb{E}_{n}):\forall k\in\left\{  1,\dots
,n-1\right\}  ,A_{n}(\mathbb{E}_{k})=\mathbb{E}_{k}\right\}
\]
can be endowed with a structure of Banach-Lie subgroup.\newline An element
$A_{n}$ of $\operatorname{G}(E_{n})$ can be seen as
\[
\centering A_{n}=\left(  \;%
\begin{tabular}
[c]{ccccllll}\cline{3-4}\cline{6-6}\cline{8-8}%
$A_{1}$ & \multicolumn{1}{c|}{$B_{2}$} &
\multicolumn{1}{c|}{\multirow{2}{*}{$B_{3}$}} &
\multicolumn{1}{c|}{\multirow{3}{*}{$B_{4}$}} & \multicolumn{1}{l|}{} &
\multicolumn{1}{l|}{\multirow{5}{*}{$B_{i}$}} &  &
\multicolumn{1}{l|}{\multirow{7}{*}{$B_{n}$}}\\
0 & \multicolumn{1}{c|}{$B_{2}^{\prime}$} & \multicolumn{1}{c|}{} &
\multicolumn{1}{c|}{} & \multicolumn{1}{l|}{} & \multicolumn{1}{l|}{} &  &
\multicolumn{1}{l|}{}\\\cline{1-3}%
\multicolumn{2}{|c|}{0} & \multicolumn{1}{c|}{$B_{3}^{\prime}$} &
\multicolumn{1}{c|}{} & \multicolumn{1}{l|}{} & \multicolumn{1}{l|}{} &  &
\multicolumn{1}{l|}{}\\\cline{1-4}%
\multicolumn{3}{|c|}{0} & $B_{4}^{\prime}$ & \multicolumn{1}{l|}{} &
\multicolumn{1}{l|}{} &  & \multicolumn{1}{l|}{}\\\cline{1-3}%
\multicolumn{1}{l}{} & \multicolumn{1}{l}{} & \multicolumn{1}{l}{} &
\multicolumn{1}{l}{} & \multicolumn{1}{l|}{$\ddots$} & \multicolumn{1}{l|}{} &
& \multicolumn{1}{l|}{}\\\cline{1-6}%
\multicolumn{5}{|c|}{0} & $B_{i}^{\prime}$ &  & \multicolumn{1}{l|}{}%
\\\cline{1-5}%
\multicolumn{1}{l}{} & \multicolumn{1}{l}{} & \multicolumn{1}{l}{} &
\multicolumn{1}{l}{} &  &  & $\ddots$ & \multicolumn{1}{l|}{}\\\hline
\multicolumn{7}{|c}{0} & \multicolumn{1}{|l}{$B_{n}^{\prime}$}\\\cline{1-7}%
\end{tabular}
\ \ \ \;\right)
\]

For $1\leq i\leq j\leq k$, we consider the following diagram%

\[%
\begin{array}
[c]{ccc}%
\mathbb{E}_{k} & \underrightarrow{A_{k}} & \mathbb{E}_{k}\\
\iota_{j}^{k}\uparrow &  & \downarrow P_{j}^{k}\\
\mathbb{E}_{j} & \underrightarrow{A_{j}} & \mathbb{E}_{j}\\
\iota_{i}^{j}\uparrow &  & \downarrow P_{i}^{j}\\
\mathbb{E}_{i} & \underrightarrow{A_{i}} & \mathbb{E}_{i}%
\end{array}
\]

where $P_{i}^{j}:\mathbb{E}_{j}\longrightarrow\mathbb{E}_{i}$ is the
projection along the direction $\mathbb{E}_{i+1}^{\prime}\oplus\cdots
\oplus\mathbb{E}_{j}^{\prime}$.

The map%
\[%
\begin{array}
[c]{cccc}%
\theta_{i}^{j}: &\operatorname{G}(E_{j}) & \longrightarrow & \operatorname{G}(E_{i})\\
& A_{j} & \mapsto & P_{i}^{j}\circ A_{j}\circ\iota_{i}^{j}%
\end{array}
\]
is perfectly defined and we have%

\[%
\begin{array}
[c]{l}%
\left(  \theta_{i}^{j}\circ\theta_{j}^{k}\right)  \left(  A_{k}\right)
=\theta_{i}^{j}\left[  \theta_{j}^{k}\left(  A_{k}\right)  \right]
=\theta_{i}^{j}\left(  P_{j}^{k}\circ A_{j}\circ\iota_{j}^{k}\right)
=P_{i}^{j}\circ P_{j}^{k}\circ A_{j}\circ\iota_{j}^{k}\circ\iota_{i}^{j}%
\end{array}
\]

Because $P_{i}^{j}\circ P_{j}^{k}=P_{i}^{k}$ (projective system) and
$\iota_{j}^{k}\circ\iota_{i}^{j}=\iota_{i}^{k}$ (inductive system), we have%
\[
\left(  \theta_{i}^{j}\circ\theta_{j}^{k}\right)  \left(  A_{k}\right)
=P_{i}^{k}\circ A_{j}\circ\iota_{i}^{k}=\theta_{i}^{k}\left(  A_{k}\right)
\]

So $\left(  \operatorname{G}(E_{i}),\theta_{i}^{j}\right)  _{i\leq j}$ is a projective system of
Banach-Lie groups and the projective limit $\operatorname{G}\left(
\mathbb{E}\right)  =\underleftarrow{\lim}\operatorname{G}(E_{n})$ can be endowed with a
structure of Fr\'{e}chet topological group.

\subsection{Convenient Lie subgroups of $\operatorname{G}\left(\mathbb{E}\right) $}
\label{__ConvenientLieSubgroupsOfG(E)}

Let $\left(   \operatorname{H}_{n}\right)  _{n\in\mathbb{N}^{\ast}}$ be a sequence where each
$\operatorname{H}_{n}$ is a weak Banach-Lie subgroup of  $\operatorname{G}(E_{n})$ where%
\[
\forall\left(  i,j\right)  \in I^{2},i\leq j,\ \theta_{i}^{j}\left(
 \operatorname{H}_{j}\right)  = \operatorname{H}_{i}%
\]

Then $( \operatorname{H}_{i},\theta_{i}^{j})_{\left(  i,j\right)  \in I^{2},i\leq j}$ is a
projective sequence and the space $\operatorname{H}\left(  \mathbb{E}\right)
=\underleftarrow{\lim} \operatorname{H}_{n}$ can be endowed with a structure of Fr\'{e}chet group.

\subsection{Direct limits of Banach vector bundles}
\label{__DirectLimitsOfBanachVectorBundles}

\begin{definition}
\label{D_AscendingSequenceBanachVectorBundles} 
A sequence $\mathcal{E}=\left(E_{n},\pi_{n},M_{n}\right)  _{n\in\mathbb{N}^{\ast}}$ of Banach vector bundles is called a strong ascending sequence of Banach vector bundles if the
following assumptions are satisfied:

\begin{description}
\item[\textbf{(ASBVB 1)}] 
$\mathcal{M}=(M_{n})_{n\in\mathbb{N}^{\ast}}$ is an
ascending sequence of Banach $C^{\infty}$-manifolds, where $M_{n}$ is modeled
on the Banach space $\mathbb{M}_{n}$ such that $\mathbb{M}_{n}$ is a
supplemented Banach subspace of $\mathbb{M}_{n+1}$ and the inclusion
$\varepsilon_{n}^{n+1}:M_{n}\longrightarrow M_{n+1}$ is a $C^{\infty}$
injective map such that $(M_{n},\varepsilon_{n}^{n+1})$ is a weak submanifold
of $M_{n+1}$;

\item[\textbf{(ASBVB 2)}] 
The sequence $(E_{n})_{n\in\mathbb{N}^{\ast}}$ is an
ascending sequence such that the sequence of typical fibers $\left(
\mathbb{E}_{n}\right)  _{n\in\mathbb{N}^{\ast}}$ of $(E_{n})_{n\in
\mathbb{N}^{\ast}}$ is an ascending sequence of Banach spaces and
$\mathbb{E}_{n}$ is a supplemented Banach subspace of $\mathbb{E}_{n+1}$;

\item[\textbf{(ASBVB 3)}] 
For each $n\in\mathbb{N}^{\ast}$, $\pi_{n+1}%
\circ\lambda_{n}^{n+1}=\varepsilon_{n}^{n+1}\circ\pi_{n}$ where $\lambda
_{n}^{n+1}:E_{n}\longrightarrow E_{n+1}$ is the natural inclusion;

\item[\textbf{(ASBVB 4)}] 
Any $x\in M=\underrightarrow{\lim}M_{n}$ has the
direct limit chart property \emph{\textbf{(DLCP) }}for $(U=\underrightarrow
{\lim}U_{n},\phi=\underrightarrow{\lim}\phi_{n})$;

\item[\textbf{(ASBVB 5)}] 
For each $n\in\mathbb{N}^{\ast}$, there exists a
trivialization $\Psi_{n}:\left(  \pi_{n}\right)  ^{-1}\left(  U_{n}\right)
\longrightarrow U_{n}\times\mathbb{E}_{n}$ such that, for any $i\leq j$, the
following diagram is commutative:
\end{description}

\[
\begin{array}
[c]{ccc}%
\left(  \pi_{i}\right)  ^{-1}\left(  U_{i}\right)  & \underrightarrow
{\lambda_{i}^{j}} & \left(  \pi_{j}\right)  ^{-1}\left(  U_{j}\right) \\
\Psi_{i}\downarrow &  & \downarrow\Psi_{j}\\
U_{i}\times\mathbb{E}_{i} & \underrightarrow{\varepsilon_{i}^{j}\times
\iota_{i}^{j}} & U_{j}\times\mathbb{E}_{j}.
\end{array}
\]

\end{definition}

For example, the sequence $\left(  TM_{n},\pi_{n},M_{n}\right)  _{n\in
\mathbb{N}^{\ast}}$ is a strong ascending sequence of Banach vector bundles
whenever $(M_{n})_{n\in\mathbb{N}^{\ast}}$ is an ascending sequence which has
the direct limit chart property at each point of $x\in M=\underrightarrow
{\lim}M_{n}$ whose model $\mathbb{M}_{n}$ is supplemented in $\mathbb{M}%
_{n+1}$.

We then have the following result given in \cite{CaPe}.

\begin{proposition}
\label{P_StructureOnDirectLimitLinearBundles} 
Let $\left(  E_{n},\pi_{n},M_{n}\right)  _{n\in\mathbb{N}^{\ast}}$ be a strong ascending sequence of
Banach vector bundles. We have:

1. $\underrightarrow{\lim}E_{n}$ has a structure of not necessarly Hausdorff
convenient manifold modelled on the LB-space $\underrightarrow{\lim}%
\mathbb{M}_{n}\times\underrightarrow{\lim}\mathbb{E}_{n}$ which has a
Hausdorff convenient structure if and only if $M$ is Hausdorff.

2. $\left(  \underrightarrow{\lim}E_{n},\underrightarrow{\lim}\pi
_{n},\underrightarrow{\lim}M_{n}\right)  $ can be endowed with a structure of
convenient vector bundle whose typical fibre is $\underrightarrow{\lim
}\mathbb{\mathbb{E}}_{n}$ and whose structural group is a Fr\'{e}chet
topological group.
\end{proposition}

\subsection{Direct limits of frame bundles}
\label{__DirectLimitsOfFrameBundles}

\begin{definition}
\label{D_AscendingSequenceTangentFrameBundles}
A sequence $\left(  \ell\left(E_{n}\right)  \right)  _{n\in\mathbb{N}^{\ast}}=\left(  P_{n}\left(
E_{n}\right)  ,\mathbf{p}_{n},M_{n},G_{n}\right)  _{n\in\mathbb{N}^{\ast}}$ of
tangent frame bundles is called an ascending sequence of tangent frame bundles
if the following assumptions are satisfied:

\begin{description}
\item[\textbf{(ASTFB 1)}] 
$\mathcal{M}=(M_{n},\varepsilon_{n}^{n+1}%
)_{n\in\mathbb{N}^{\ast}}$ is an ascending sequence of Banach manifolds;

\item[\textbf{(ASTFB 2)}] 
For any $n\in\mathbb{N}^{\ast}$, $\left(
\lambda_{n}^{n+1},\varepsilon_{n}^{n+1},\gamma_{n}^{n+1}\right)  $ is a
morphism of principal bundles from $\left(  P_{n}\left(  E_{n}\right)
,\mathbf{p}_{n},M_{n},G_{n}\right)  $ to $\left(  P_{n+1}\left(  E_{n+1}\right)
,\mathbf{p}_{n+1},M_{n+1},G_{n}\right)  $ where $\lambda_{n}^{n+1}$ is defined via
local sections $s_{\alpha}^{n}$ of $P_{n}\left(  E_{n}\right)  $ and fulfils
the condition $\lambda_{n}^{n+1}\circ s_{\alpha}^{n}=s_{\alpha}^{n+1}%
\circ\varepsilon_{n}^{n+1}$;

\item[\textbf{(ASTFB 3)}] 
Any $x\in M=\underrightarrow{\lim}M_{n}$ has the direct limit chart property for 
$(U=\underrightarrow{\lim}U_{n},\phi=\underrightarrow{\lim}\phi_{n})$;

\item[\textbf{(ASTFB 4)}] 
For any $x\in M=\underrightarrow{\lim}M_{n}$, there
exists a trivialization $\Psi_{n}:\left(  \pi_{n}\right)  ^{-1}\left(
U_{n}\right)  \longrightarrow U_{n}\times\mathbb{E}_{n}$ such that the
following diagram is commutative:

\[
\begin{array}
[c]{ccc}%
\left(  \pi_{n}\right)  ^{-1}\left(  U_{n}\right)  & \underrightarrow
{\lambda_{n}^{n+1}} & \left(  \pi_{n+1}\right)  ^{-1}\left(  U_{n+1}\right) \\
\Psi_{n}\downarrow &  & \downarrow\Psi_{n+1}\\
U_{n}\times\mathbb{E}_{n} & \underrightarrow{\left(  \varepsilon_{n}%
^{n+1}\times\iota_{n}^{n+1}\right)  } & U_{n+1}\times\mathbb{E}_{n+1}.
\end{array}
\]

for each $n\in\mathbb{N}^{\ast}$.
\end{description}
\end{definition}

\begin{theorem}
\label{T_ConvenientStructureOnDirectLimitOfAscendingSequenceFrameBundles}
If $\left(  \ell\left(  E\right)  \right)  _{n\in\mathbb{N}^{\ast}}=\left(
P_{n}\left(  E\right)  ,\mathbf{p}_{n},M_{n},G_{n}\right)  _{n\in\mathbb{N}^{\ast}}$
is an ascending sequence of frame bundles, then the direct limit
\[
\left(  \underrightarrow{\lim}P_{n}\left(  E\right)  ,\underrightarrow{\lim
}\mathbf{p}_{n},\underrightarrow{\lim}M_{n},\operatorname*{G}(\mathbb{E})\right)
\]
can be endowed with a structure of convenient principal bundle.
\end{theorem}

\begin{proof}
According to (ASTFB 1) and (ASTFB 3) corresponding to (DLCP),
$M=\underrightarrow{\lim}M_{n}$ can be endowed with a structure of non
necessary Hausdorff convenient manifold (cf. \cite{CaPe}, \S \ Direct limit of
Banach manifolds).\newline From (ASTFB 2) and (PB 1) we deduce $\pi_{n+1}%
\circ\lambda_{n}^{n+1}=\varepsilon_{n}^{n+1}\circ\pi_{n}$. So the projection
$\pi=\underrightarrow{\lim}\pi_{n}:\underrightarrow{\lim}P_{n}\left(
E_{n}\right)  \longrightarrow\underrightarrow{\lim}M_{n}$ is well defined. Let
$u$ be an element of $P=\underrightarrow{\lim}P_{n}\left(  E_{n}\right)$; so
$u$ belongs to some $P_{n}\left(  E_{n}\right)  $. In particular,
$x=\pi\left(  u\right)  \in M_{n}$ where $\pi=\underrightarrow{\lim}\pi_{n}$.
According to (ASTFB 3) and (ASTFB 4), there exists a local chart $\left(
U,\phi\right)  $ of the convenient manifold $M$ whose domain contains $x$.
Moreover, there exists a local trivialization $\Psi_{n}:\pi_{n}^{-1}\left(
U_{n}\right)  \longrightarrow U_{n}\times G_{n}$ of $P_{n}\left(
E_{n}\right)  $.\newline For $i\geq n$, one can define the local trivilization
$\Psi_{i}:\pi_{i}^{-1}\left(  U_{i}\right)  \longrightarrow U_{i}\times G_{i}$
of $P_{i}\left(  E_{i}\right)  $. According to (ASTFB 2), the map $\Psi
:\pi^{-1}\left(  U\right)  \longrightarrow U\times$ $\operatorname*{G}%
(\mathbb{E})$ can be defined and it is one to one.\newline Let us study the
transition functions. Using the commutativity of the diagram%
\[%
\begin{array}
[c]{ccc}%
\left(  U_{i}^{\alpha}\cap U_{i}^{\beta}\right)  \times\operatorname*{G}_{i} &
\overset{\varepsilon_{i}^{j}\times J_{i}^{j}}{\longrightarrow} & \left(
U_{j}^{\alpha}\cap U_{j}^{\beta}\right)  \times\operatorname*{G}_{j}\\
\Psi_{i}^{\alpha}\circ\left(  \Psi_{i}^{\beta}\right)  ^{-1}\downarrow &  &
\downarrow\Psi_{j}^{\alpha}\circ\left(  \Psi_{j}^{\beta}\right)  ^{-1}\\
\left(  U_{i}^{\alpha}\cap U_{i}^{\beta}\right)  \times\operatorname*{G}_{i} &
\overset{\varepsilon_{i}^{j}\times J_{i}^{j}}{\longrightarrow} & \left(
U_{j}^{\alpha}\cap U_{j}^{\beta}\right)  \times\operatorname*{G}_{j}%
\end{array}
\]
we obtain transition functions $\Psi^{\alpha}\circ\left(  \Psi^{\beta}\right)
^{-1}$which are diffeomorphisms of $\left(  U^{\alpha}\cap U^{\beta}\right)
\times$ $\operatorname*{G}\left(  \mathbb{E}\right)  .$\newline We then have a
set of local trivializations $\left\{  \left(  u^{\alpha},\Psi^{\alpha
}\right)  \right\}  _{\alpha\in A}$ of $P$ where the projection $\pi
:P\longrightarrow M$ is smooth.\newline This structure is endowed with a right
action well defined according to {(ASTFB 2)}.
\end{proof}

\subsection{Direct limit of $G$-structures}

\label{__DirectLimitofGStructures}

Considering the reductions of the frame bundles, we then get the following
result whose proof is rather analogous to the previous one: we have to check
the different compatibility conditions and replace the structure group
$\operatorname*{G}(\mathbb{E})$ by $H\left(  \mathbb{E}\right)  $.

\begin{theorem}
\label{T_ConvenientStructureOnDirectLimitOfAscendingSequenceOfGStructures}Let
$\left(  \ell\left(  E_{n}\right)  \right)  _{n\in\mathbb{N}^{\ast}}=\left(
P_{n}\left(  E_{n}\right)  ,\mathbf{p}_{n},M_{n},G_{n}\right)  _{n\in\mathbb{N}%
^{\ast}}$ be an ascending sequence of frame bundles.\newline Let $\left(
F_{n},\mathbf{p}_{n|F_{n}},M_{n},H_{n}\right)  _{n\in\mathbb{N}^{\ast}}$ be a
sequence of associated $G$-structures where, for any $n\in\mathbb{N}^{\ast},$
$\left(  \left(  \lambda_{n}^{n+1}\right)  _{|F_{n}},\varepsilon_{n}%
^{n+1},\theta_{n}^{n+1}\right)  $ is a morphism of principal bundles.

Then $\left(  \underrightarrow{\lim}F_{n},\underrightarrow{\lim}\mathbf{p}_{n|F_{n}%
},\underrightarrow{\lim}M_{n},H(\mathbb{E})\right)  $ can be endowed with a
convenient $G$-structure.
\end{theorem}

\subsection{Direct limits of tensor structures}
\label{__DirectLimitsOfTensorStructures}

\index{DirectLimitsOfTensorStructures}

\subsubsection{Direct limits of tensor structures of type $\left(  1,1\right)$}
\label{___DirectLimitsOfTensorStructuresOfType1-1}

In a similar way in the projective system of tensor of type $(1,1)$, we also introduce the following definition and results.

\begin{definition} 
\label{D_DirectSequenceOfEndomorphisms}${}$
\begin{enumerate}
\item[1.] 
Let $\left\{  \left(  \mathbb{E}%
_{i},\overline{\lambda_{i}^{j}}\right)  \right\}  _{\left(  i,j\right)
\in\mathbb{N}^{\ast}\times\mathbb{N}^{\ast},\ i\leq j\ }$ be an ascending
sequence of supplemented Banach spaces.

For any $n\in\mathbb{N}^{\ast}$, let $A_{n}:\mathbb{E}_{n}\longrightarrow
\mathbb{E}_{n}$ be an endomorphism of the Banach space $\mathbb{E}_{n}%
.$\\
$\left(  A_{n}\right)  _{n\in\mathbb{N}^{\ast}}$ is called
an ascending sequence of endomorphisms if, for any integer $j\geqslant i>0$,
it fulfils the coherence condition,
\[
\overline{\lambda_{i}^{j}}\circ A_{i}=A_{j}\circ\overline{\lambda_{i}^{j}}.
\]
\item[2.] 
Let $\left\{  \left(  E_{n},\pi_{n},M_{n}\right)  \right\}  _{n\in
\mathbb{N}^{\ast}}$ be an ascending sequence of of supplemented  Banach vector bundles.  
A sequence $\left(  \mathcal{A}_{n}\right)  _{n\in\mathbb{N}^{\ast
}}$  endomorphisms $\mathcal{A}_{n}$ of  $E_n$ is called a coherent sequence of endomorphisms  if for each $x=\underleftarrow{\lim}x_n$, the sequence $\left(\left(\mathcal{A}_n\right)_{x_n}\right)$ is a coherent sequence of endomorphism of $\left(E_n\right)_{x_n}$.\\
\end{enumerate}
\end{definition}

\begin{proposition}
\label{P_DirectLimitOfEndomorphisms}
\begin{enumerate}
\item[1.] Consider   a sequence $\left(  {A}_{n}\right)  _{n\in\mathbb{N}^{\ast
}}$ of coherent   endomorphisms  ${A}_{n}$  of $E_n$ on   an ascending  sequence of supplemented Banach spaces   $\left\{  \left(  \mathbb{E}_{n}\right)  \right\}  _{n\in
\mathbb{N}^{\ast}}$. If $G(A _n)$ is the isotropy group of $A_n$, then $G(A_n)$ is a convenient subgroup of $\operatorname{G}(\mathbb{E})$.

\item[2.]Consider   a sequence $\left(  \mathcal{A}_{n}\right)  _{n\in\mathbb{N}^{\ast
}}$ of coherent   endomorphisms  $\mathcal{A}_{n}$  of $E_n$ on   an ascending sequence of supplemented Banach vector bundles  $\left\{  \left(  E_{n},\pi_{n},M_{n}\right)  \right\}  _{n\in
\mathbb{N}^{\ast}}$. Then the direct limit $\mathcal{A}=\underleftarrow{\lim}({\mathcal{A}}_n)$ is a well defined and is a smooth endomorphism of the Fr\'echet bundle  $E=\underleftarrow{\lim}{E_n}$.\\
\end{enumerate}
\end{proposition}

\begin{proof}
1. Under our assumption, we have $\mathbb{E}_{1}\subset\mathbb{E}_{2}\subset\cdots$  and Banach spaces $\mathbb{E}%
_{1}^{\prime},\mathbb{E}_{2}^{\prime},\dots$ such that:
\[
\left\{
\begin{array}
[c]{c}%
\mathbb{E}_{1}=\mathbb{E}_{1}^{\prime},\\
\forall i\in\mathbb{N}^{\ast},\mathbb{E}_{i+1}\backsimeq\mathbb{E}_{i}%
\times\mathbb{E}_{i+1}^{\prime}%
\end{array}
\right.  \;
\]
Of course ${G}(A_1)\subset \operatorname{GL}(\mathbb{E})= \operatorname{G}(\mathbb{E})$. Assume that for all $1\leq i\leq n$ we 
have ${G}(A_i) \subset  \operatorname{G}(\mathbb{E}_i)$. Since $ \mathbb{E}_{n+1}\backsimeq\mathbb{E}_{n}\times\mathbb{E}_{n+1}^{\prime}$ 
and the sequence $(A_n)_{n\in \mathbb{N}^*}$ is coherent, if $T_{n+1} $ belongs to $G(A_{n+1})$ it follows that 
$T_{n+1}\circ A_n=T_{n+1}\circ A_n$. Thus if $T_n$ is the restriction of $T_{n+1}$ to $E_n$, then $T_n$ belongs to  $\operatorname{GL}(\mathbb{E}_n)$ and so $T_n$ can be written as a matrix of type
$$\begin{pmatrix}
T_n&T'_n\\
0& S_{n+1}\\
\end{pmatrix}$$
Thus  $T_{n+1}$ belongs to $\operatorname{G}(\mathbb{E}_{n+1})$ which ends the proof of Point 1.\\

2. By definition, a coherent sequence of endomorphism is nothing but else an ascending system of linear map so the direct limit ${A}=\underrightarrow{\lim}{{A}}_n$ is a well  endomorphism of $\mathbb{E}=\underrightarrow{\lim}{\mathbb{E}}_n$. It follows that for each $x=\underrightarrow{\lim}{x}_n\in M=\underrightarrow{\lim}{M}_n$ the direct limit 
$\mathcal{A}_x=\underrightarrow{\lim}(\mathcal{A}_n)_{x_n}$ is well defined. Now from property (ASBVB 5)
for each $n\in\mathbb{N}^{\ast}$, there exists a
trivialization $\Psi_{n}:\left(  \pi_{n}\right)  ^{-1}\left(  U_{n}\right)
\longrightarrow U_{n}\times\mathbb{E}_{n}$ such that, for any $i\leq j$, the
following diagram is commutative:

\[
\begin{array}
[c]{ccc}%
\left(  \pi_{i}\right)  ^{-1}\left(  U_{i}\right)  & \underrightarrow
{\lambda_{i}^{j}} & \left(  \pi_{j}\right)  ^{-1}\left(  U_{j}\right) \\
\Psi_{i}\downarrow &  & \downarrow\Psi_{j}\\
U_{i}\times\mathbb{E}_{i} & \underrightarrow{\varepsilon_{i}^{j}\times
\iota_{i}^{j}} & U_{j}\times\mathbb{E}_{j}.
\end{array}
\]
This implies that the restriction of $\mathcal{A}_n$ to $\left(  \pi_{i}\right)  ^{-1}\left(  U_{n}\right)$ is an ascending system of smooth maps and so from Lemma 22 of \cite{CaPe}, the restriction of $\mathcal{A}$ to  $  \pi  ^{-1}\left(  U\right)$ is a smooth map  where $U=\underrightarrow{\lim}U_n$ which ends the proof.\\
\end{proof} 

Let $\left\{  \left(  E_{n},\pi_{n},M_{n}\right)  \right\}  _{n\in
\mathbb{N}^{\ast}}$ be an ascending sequence of supplemented Banach vector bundles.  Consider a coherent  sequence $\left(  \mathcal{A}_{n}\right)  _{n\in\mathbb{N}^{\ast
}}$ of endomorphisms  $\mathcal{A}_n$ defined  on 
the  Banach bundle $\pi_{n}:E_{n}\longrightarrow M_{n}.$  \\
Fix  any  $x^0=\underrightarrow{\lim}{x_n^0}\in M=\underrightarrow{\lim}M_n$  and identify each typical fiber of $E_n$ with $(E_n)_{x_n^0}$. For all $n\in \mathbb{N}^*$, we denote by $G((\mathcal{A}_n)_{x_n^0})$ the isotropy group of $(\mathcal{A}_n)_{x_n^0}$ . \\

Now by the same type of arguments as the ones used in the proof of Theorem \ref{T_ReductionProjectiveSequenceTensor1-1} but replacing "projective sequence" by "ascending sequence" and Fr\'echet atlas by convenient atlas, we obtain:

\begin{theorem}
\label{T_DirectReductionDirectLimAscendingSequenceOfTensors1-1} 
Let $\left\{  \left(  E_{n},\pi_{n},M_{n}\right)  \right\}  _{n\in
\mathbb{N}^{\ast}}$ be an ascending sequence of supplemented  Banach vector bundles.  Consider a coherent  sequence $\left(  \mathcal{A}_{n}\right)  _{n\in\mathbb{N}^{\ast
}}$ of endomorphisms  $\mathcal{A}_n$ defined  on 
the  Banach bundle $\pi_{n}:E_{n}\longrightarrow M_{n}.$

Then $\underrightarrow{\lim}\mathbf{P}_{n}\left(  E_{n}\right)  $ can be
endowed with a structure of convenient principal bundle over $\underrightarrow
{\lim}M_{n}$ whose structural group is $G\left(  \mathcal{A}_{x^0}%
\right) = \underrightarrow{\lim}G(\mathcal{A}_n)_{x^0_n}$ if and only if there exists a convenient atlas bundle   $\left\{ U^\alpha= \underleftarrow{\lim} U_{n}^{\alpha}, \tau^{\alpha}=\tau_{n}^{\alpha}  \right\}
_{n\in\mathbb{N}^{\ast},\alpha\in A}$   such that  $\mathcal{A}_{n}^{\alpha}=\mathcal{A}_{n|U_{n}^{\alpha}}$ is locally 
modelled  $(\mathcal{A}_n)_{x_n^0}$
 and  the transition maps $T_{n}^{\alpha\beta}\left(x_{n}\right)  $ take values in  to the isotropy group $G((\mathcal{A}_n)_{x_n^0})$, for  all $n\in \mathbb{N}^*$ and all $\alpha,\beta \in A$ such that $U_n^\alpha\cap U_n^\beta\not=\emptyset$ .\\   
\end{theorem}

\subsubsection{Direct limits of tensor structures of type $\left(  2,0\right)$}
\label{___DirectLimitsOfTensorStructuresOfType2-0}

We adapt the results obtained for coherent sequences tensor of type $(1,1)$ to such sequences of $(2,0)$-tensors.
 
\begin{definition}
\label{D_AscendingSequenceOfBilinearForms}${}$
\begin{enumerate}
\item[1.] 
Let $\left\{  \left(  \mathbb{E}%
_{i},\overline{\lambda_{i}^{j}}\right)  \right\}  _{\left(  i,j\right)
\in\mathbb{N}^{\ast}\times\mathbb{N}^{\ast},\ i\leq j\ }$ be an ascending
sequence of supplemented of Banach spaces.

For any $n\in\mathbb{N}^{\ast}$, let $\omega_{n}:\mathbb{E}_{n}\times
\mathbb{E}_{n}\longrightarrow\mathbb{R}$ be a bilinear form of the Banach
space $\mathbb{E}_{n}.$\\
$\left(  \omega_{n}\right)_{n\in\mathbb{N}^{\ast}}$ is called a projective sequence of $2$-forms if it fulfils the coherence condition, for any integer $j\geqslant i>0$%
\[
\omega_{i}=\left(  \overline{\lambda_{i}^{j}}\times\overline{\lambda_{i}^{j}%
}\right)  \circ\omega_{i}%
\]
\item[2.] 
Let $\left\{  \left(  E_{n},\pi_{n},M_{n}\right)  \right\}  _{n\in
\mathbb{N}^{\ast}}$ be a projective sequence of Banach vector bundles.  
A sequence $\left(  \Omega_{n}\right)  _{n\in\mathbb{N}^{\ast
}}$  bilinear form $\Omega_{n}$ on  $E_n$ is called a coherent sequence of bilinear  
   if for each $x=\underrightarrow{\lim}x_n$, the sequence $\left(\left(\Omega_n\right)_{x_n}\right)$ is a coherent sequence of bilinear form of $\left(E_n\right)_{x_n}$.\\
\end{enumerate}
\end{definition}

We have the following properties:

\begin{proposition}
\label{P_DirectLimitOfBilinearForms}
\item[1.] Consider   a sequence $\left(  \omega_{n}\right)  _{n\in\mathbb{N}^{\ast
}}$ of coherent   of bilinear forms  $\omega_{n}$  of $E_n$ on   an ascending  sequence of supplemented Banach spaces   $\left\{  \left(  \mathbb{E}_{n}\right)  \right\}  _{n\in
\mathbb{N}^{\ast}}$. If $G(\omega_n)$ is the isotropy group of $\omega_n$, then $G(\omega)=\underrightarrow{\lim} G(\omega_n)$ is a convenient subgroup of $\operatorname{G}(\mathbb{E})$.

 \item[2.] Consider   a sequence $\left(  \mathcal{A}_{n}\right)  _{n\in\mathbb{N}^{\ast
}}$ of coherent   bilinear forms  $\Omega_{n}$  of $E_n$ on   an ascending sequence of supplemented Banach vector bundles  $\left\{  \left(  E_{n},\pi_{n},M_{n}\right)  \right\}  _{n\in
\mathbb{N}^{\ast}}$. Then the projective limit $\Omega=\underrightarrow{\lim}(\Omega_n)$ is a well defined and is a smooth bilinear form of the convenient bundle  $E=\underrightarrow{\lim}{E_n}$.\\
\end{proposition}

Let $\left\{  \left(  E_{n},\pi_{n},M_{n}\right)  \right\}  _{n\in
\mathbb{N}^{\ast}}$ be an ascending sequence of supplemented Banach vector bundles.  Consider a coherent  sequence $\left( \Omega_n\right)  _{n\in\mathbb{N}^{\ast
}}$ of bilinear forms  $\Omega_n$ defined  on 
the  Banach bundle $\pi_{n}:E_{n}\longrightarrow M_{n}.$  \\
Fix  any  $x^0=\underrightarrow{\lim}{x_n^0}\in M=\underrightarrow{\lim}M_n$  and identify each typical fiber of $E_n$ with $(E_n)_{x_n^0}$. For all $n\in \mathbb{N}^*$, we denote by $G((\Omega_n)_{x_n^0})$ the isotropy group of $(\Omega_n)_{x_n^0}$. \\

\begin{theorem}
\label{T_ReductionAscendingTensor0-2}
 Let $\left\{  \left(  E_{n},\pi_{n},M_{n}\right)  \right\}  _{n\in
\mathbb{N}^{\ast}}$ be an ascending sequence of supplemented Banach vector bundles.  Consider a coherent  sequence $\left(  \Omega_{n}\right)  _{n\in\mathbb{N}^{\ast
}}$ of bilinear forms  $\Omega_n$ defined  on 
the  Banach bundle $\pi_{n}:E_{n}\longrightarrow M_{n}.$

Then $\underrightarrow{\lim}\mathbf{P}_{n}\left(  E_{n}\right)  $ can be
endowed with a structure of Fr\'{e}chet principal bundle over $\underleftarrow
{\lim}M_{n}$ whose structural group can be reduced to $G((\Omega)_{x^0})= \underleftarrow{\lim}G((\Omega_n)_{x_n^0})
  $ if and only if there exists a convenient atlas bundle  $\left\{ U^\alpha= \underleftarrow{\lim} U_{n}^{\alpha}, \tau^{\alpha}=\tau_{n}^{\alpha}  \right\}
_{n\in\mathbb{N}^{\ast},\alpha\in A}$   such that  $\Omega_{n}^{\alpha}=\mathcal{A}_{n|U_{n}^{\alpha}}$ is locally 
modelled  $(\Omega_n)_{x_n^0}$
 and  the transition maps $T_{n}^{\alpha\beta}\left(x_{n}\right)  $ take values in  to the isotropy group $G((\Omega_n)_{x_n^0})$, for  all $n\in \mathbb{N}^*$ and all $\alpha,\beta \in A$ such that $U_n^\alpha\cap U_n^\beta\not=\emptyset$ .\\   
\end{theorem}

\subsection{Examples of direct limits of tensor structures}
\label{__ExampleDirectLimitsOfTensorStructures}

\subsubsection{Direct limit of compatible weak almost K\"ahler and almost para-K\"ahler structures}
\label{___DirectLimitOfCompatibleWeakAlmostKalherAlmostParaKalherStructures}
We consider the context and notations of $\S$ \ref{__CompatibilityBetwenDifferentStructures}.\\
If $\Omega$  (resp. $g$,  resp. $\mathcal{I}$) is  a weak symplectic
form  (resp. a weak Riemannian metric, resp. an almost complex structure) on a convenient  bundle $(E,\pi_{E},M)$, the compatibility of any pair of the data $(\Omega, g,\mathcal{I})$ as given in Definition \ref{D_CompatibilityWeakSymplecticFormAlmostComplexStructureWeakRiemannianMetricOnABanachBundle} can be easily adapted to this convenient  context. We can also easily  transpose the notion of  "local tensor $\mathcal{T}$ locally modelled on  a linear tensor $\mathbb{T}$" (cf. Definition \ref{D_TensorLocallyModelled}) in the convenient framework. We will see that such  situations can be obtained by some "direct limit"
 of  a  sequence $\{\left(E_n,\pi_n,M_n\right)\}_{n\in \mathbb{N}^*}$ of Banach vector bundles. Once again, this context is very similar to the context of $\S$ \ref{__ExampleOfProjectiveLimitsOfTensorStructures} and we will only write down the details without any proof.   \\
 
First we introduce the following corresponding notations:
\begin{enumerate}
\item[--] Let $(\mathbb{T}_n)$ be a sequence of coherent tensor of type $(1,1)$ or $(2,0)$ on an ascending sequence of supplemented Banach spaces $\left\{  \left(\mathbb{E}_{i},\overline{\lambda_{i}^{j}}\right)  \right\}  _{\left(
i,j\right)  \in\mathbb{N}^{\ast}\times\mathbb{N}^{\ast},\ i\leq j\ }$. We denote by $G(\mathbb{T}_n)$ the isotropy group of $\mathbb{T}_n$. If $\mathbb{T}=\underrightarrow{\lim}\mathbb{T}_n$,    ${G}(\mathbb{T})=\underrightarrow{\lim}G(\mathbb{T}_n)$ is convenient    sub-group of $\operatorname{G}(\mathbb{E})$
\item[--]  Given  an ascending sequence $\{E_n,\pi_n,M_n\}_{n\in \mathbb{N}^*}$ of supplemented Banach vector bundles, we provide 
 each Banach bundle $E_n$ with a symplectic form $\Omega_n$,and/or  a weak Riemaniann metric $g_n$ and/or  an almost  complex structure  $\mathcal{I}_n$.
\end{enumerate} 

According to this context, let us consider the following assumption:

\begin{description}
\item[\textbf{(H)}]
Assume that we have a sequence of weak Riemaniann metrics $(g_n)_{n\in \mathbb{N}^*}$ on the sequence $(E_n)_{n\in \mathbb{N}^*}$.  For each $x_n\in M_n$, we denote by $(\widehat{E}_n)_{x_n}$ and we set $\widehat{E}_n=\bigcup\limits_{x_n\in M_n}(\widehat{E}_n)_{x_n}$ and we set $\widehat{E}=\bigcup\limits_{x_n\in M_n}  (\widehat{E}_n)_{x_n}$. Then  $\widehat{\pi}_n:\widehat{E}_n\longrightarrow M_n$ has a  Hilbert vector bundle structure and the inclusion $\iota_n:E_n\rightarrow \widehat{E}_n$ is a Banach bundle morphism.
\end{description}

From Theorem \ref{T_CompatibilityDarbouxFormWeakRiemannianMetricComplexStructureOnABanachBundle}, we then  obtain
 
\begin{theorem}
\label{T_ InductiveLimitOfCompatibleWeakSymplecticFormsRiemannianMetricsComplexStructures}${}$
\begin{enumerate}
\item[1.] 
Under the assumption (H),  $\left\{ \left( \widehat{E}_n,\widehat{\pi}_n, M_n \right) \right\}_{n\in \mathbb{N}^*}$ is  an ascending  sequence of supplemented Banach bundles and there exists an injective convenient bundle morphism $\iota: E=\underrightarrow{\lim}E_n\rightarrow \widehat{E}=\underrightarrow{\lim}\widehat{E}_n$ with dense range.\\ 
Let $\widehat{\Omega}_n$, (resp. $\widehat{g}_n$, resp. $\widehat{\mathcal{I}}_n$) the extension of ${\Omega}_n$ (resp. ${g}_n$, resp.  ${\mathcal{I}}_n$) to $\widehat{E}_n$  (cf. Theorem \ref{T_CompatibilityDarbouxFormWeakRiemannianMetricComplexStructureOnABanachBundle}). If $({\Omega}_n)_{n\in \mathbb{N}^*}$, (resp. $({g}_n)_{n\in \mathbb{N}^*}$, resp. $(\mathcal{I}_n)_{n\in \mathbb{N}^*}$) is a coherent sequence, so is the sequence $(\widehat{\Omega}_n)_{n\in \mathbb{N}^*}$ (resp. $(\widehat{g}_n)_{n\in\mathbb{N}^*}$, resp.$(\widehat{\mathcal{I}}_n)_{n\in \mathbb{N}^*}$). 
  We set $\widehat{\Omega}=\underrightarrow{\lim}\widehat{\Omega}_n$ (resp. $\widehat{g}=\underrightarrow{\lim}\widehat{g}_n$, resp.
  $\widehat{\mathcal{I}}=\underrightarrow{\lim}\widehat{\mathcal{I}}_n$).  We then have
 ${\Omega}=\underrightarrow{\lim}{\Omega}_n=i^*\widehat{\Omega}$, $g=\underrightarrow{\lim}{g}_n=\iota^*\widehat{g}$ and 
 $\mathcal{I}=\underrightarrow{\lim}{\mathcal{I}}_n=\iota^*\widehat{\mathcal{I}}$. Moreover,  if any pair of the  data $(\Omega_n, g_n, \mathcal{I}_n)$ is compatible for all $n\in \mathbb{N}$, then  $(\widehat{\Omega},\widehat{g},\widehat{\mathcal{I}})$ (resp. $({\Omega},g, {\mathcal{I}})$) on the Fr\'echet vector bundle  $\widehat{E}$ (resp. $E$) respectively is compatible. 
\item[2.] 
Let $\left(  \mathbf{\ell}\left(  E_{n}\right)  \right)  _{n\in\mathbb{N}%
^{\ast}}=\left(  \mathbf{P}\left(  E_{n}\right)  ,\pi_{n},M_{n},\mathcal{H}%
_{0}^{n}\left(  \mathbb{E}\right)  \right)  _{n\in\mathbb{N}^{\ast}}$ be the associated
sequence of generalized frame bundles. Denote by $\mathcal{T}_n$ anyone tensor among $(\Omega_n,g_n,\mathcal{I}_n)$. Assume that $\mathcal{T}_n$  defines a $\mathbb{T}$-tensor structure on $E_n$ for each $n\in \mathbb{N}^*$.  Then $\underrightarrow{\lim}\mathbf{P}_{n}\left(  E_{n}\right)  $  %and $\underrightarrow{\lim}\mathbf{P}_{n}\left(  \widehat{E}_{n}\right)  $ 
is a convenient principal bundle over $\underrightarrow
{\lim}M_{n}$ whose structural group can be reduced to ${G}\left(  \mathbb{T}
\right)$. Moreover,  if the sequence  $({\Omega}_n,g_n, {\mathcal{I}}_n)_{n\in \mathbb{N}^*}$ is coherent and each triple is a tensor structure on $E_n$, then %(\widehat{\Omega},\widehat{g}, \widehat{\mathcal{I}})$ and 
  $({\Omega},g, {\mathcal{I}})$  is also a tensor structure on $E$ which is compatible if $({\Omega}_n,g_n, {\mathcal{I}}_n)$ are so.
\item[3.] 
Under the assumption (H),  all the properties of Point 2 are also valid for the sequence $(\widehat{\Omega}_n,\widehat{g}_n, \widehat{\mathcal{I}}_n)_{n\in \mathbb{N}^*} $ relatively to $\widehat{E}$.
\end{enumerate}
\end{theorem}
 
%\begin{remark}
%\label{R_DirectLimitParaKahler} 
%Given  an ascending  sequence $\{E_n,\pi_n,M_n\}_{n\in \mathbb{N}^*}$ of supplemented Banach vector bundles, we provide 
% each Banach bundle $E_n$ with a symplectic form $\Omega_n$ and/or  a neutral metric $g_n$ and/or  an  almost  para-complex structure $\mathcal{I}_n$. Then by formal analogue arguments, we can obtain a similar Theorem as Theorem \ref{T_ DirectLimitOfCompatibleWeak}.
%\end{remark}

%Now  Theorem \ref{T_ ProjectiveLimitOfCompatibleweaksymplecticFormsRiemannian metricscomplex structures} and Remark \ref{R_ProjectiveLimitParaKahler} we obtain directly

\begin{corollary} 
\label{C_DirectLimitAlmostKahlerStructures}
 Consider an almost K\"ahler (resp. para-K\"ahler) structure $({\Omega}_n,g_n, {\mathcal{I}}_n)_{n\in \mathbb{N}^*}$ on an ascending sequence $\left\{ \left( E_n,\pi_n,M_n \right) \right\}_{n\in \mathbb{N}^*}$ of Banach vector bundles. If $({\Omega}_n,g_n, {\mathcal{I}}_n)_{n\in \mathbb{N}^*}$ is coherent, then $({\Omega}=\underrightarrow{\lim}{\Omega}_n,g=\underrightarrow{\lim}{g}_n,\mathcal{I}=\underrightarrow{\lim}{\mathcal{I}}_n)$ is an almost K\"ahler (resp. para-K\"alher) structure on the Fr\'echet bundle $E=\underrightarrow{\lim} E_n$.
\end{corollary}

\subsubsection{Application to Sobolev loop spaces}
\label{___ApplicationToSobolevSpaces}

It is well known that if $M$ is a connected manifold of dimension $m$, the set $\mathsf{L}_{k}^p(\mathbb{S}^1,M)$ of loops $\gamma:\mathbb{S}^1\longrightarrow M$ of Sobolev  class $\mathsf{L}_k^p$ is a Banach manifold  modelled on the Sobolev space $\mathsf{L}_k^p(\mathbb{S}^1,\mathbb{R}^m)$ (see \cite{Pel2}).
 
\medskip

Let us consider an ascending sequence $\{\left(  M_{n},\omega_{n}\right)
\}_{n\in\mathbb{N}^{\ast}}$ of finite dimensional manifolds. Then  the direct limit $M=\underrightarrow{\lim}M_{n}$ is modelled on the convenient space $\mathbb{R}^{\infty}$. According to \cite{Pel2}, Proposition 34,  the direct limit  $\mathsf{L}_{k}^{p}(\mathbb{S}^{1},M)=\underrightarrow{\lim}(\mathsf{L}_{k}^{p}(\mathbb{S}^{1},M_{n}))$  has a Hausdorff  convenient manifold structure modelled on the convenient space $\mathsf{L}_{k}^{p}(\mathbb{S}^{1},\mathbb{R}^{\infty
})=\underrightarrow{\lim}(\mathsf{L}_{k}^{p}(\mathbb{S}^{1},\mathbb{R}^{m}))$.\\

Assume now that each $M_n$ is a K\"ahler manifold and let  $\omega_n$, $\mathfrak{g}_n$ and $\mathfrak{I}_n$ the associated symplectic form, Riemannian metric and complex structure on $M_n$ respectively. Moreover, we suppose that each sequence  $(\omega_n)_{n\in \mathbb{N}^*}$, $(\mathfrak{g}_n)_{n\in \mathbb{N}^*}$ and $(\mathfrak{I}_n)_{n\in \mathbb{N}^*}$ is a coherent sequence on the ascending sequence $(M_n)_{n\in \mathbb{N}^*}$.

For $p$ and $k$ fixed, on the tangent bundle  $T\mathsf{L}_{k}^p(\mathbb{S}^1,M_n)$,  we introduce  (cf. \cite{Kum1} and \cite{Kum2}):
\begin{enumerate}
\item[--] 
$(\Omega_n)_\gamma(X,Y)=\displaystyle\int_{\mathbb{S}^1}\omega_n\left(X_{\gamma(t)},Y_{\gamma(t)}\right)dt$;
\item[--] 
$(g_n)_\gamma(X,Y)=\displaystyle\int_{\mathbb{S}^1}\mathfrak{g}_n\left(X_{\gamma(t)},Y_{\gamma(t)}\right)d\nu(t)$;
\item[--] 
$(\mathcal{I}_n)_\gamma(X)(t)=\mathfrak{I}_n\left(X_{\gamma(t)}\right)$.
\end{enumerate}

It is clear that $\Omega_n$ (resp. $g_n$) is a weak bilinear form (resp. a Riemaniann metric) on $\mathsf{L}_{k}^p(\mathbb{S}^1,M)$ and $\mathcal{I}_n$ is an almost complex structure on $T\mathsf{L}_{k}^p(\mathbb{S}^1,M)$. The compatibility of the triple $(\omega, \mathfrak{g}, \mathfrak{I})$ implies clearly the compatibility of any $(\Omega_n, g_n,\mathcal{I}_n)$. Moreover, $(\Omega_n, g_n,\mathcal{I}_n)_{n\in \mathbb{N}^*}$ is a coherent sequence. Then, from Corollary \ref{C_DirectLimitAlmostKahlerStructures}, we get an almost K\"ahler structure on $\mathsf{L}_{k}^{p}(\mathbb{S}^{1},M).$

\textit{If each $M_n$ is a para-K\"ahler manifold, we can provide  $\mathsf{L}_{k}^{p}(\mathbb{S}^{1},M)$ with an almost para-K\"ahler structure in the same way.}

\begin{remark}
As in the projective limit context, if each $M_n$ is a K\"ahler manifold, then the almost complex structure is integrable and, in particular, the Nijenhuis tensor $N_{\mathfrak{I}_n}$ is zero. It is clear that we also have $N_{\mathcal{I}_n}\equiv 0$ for all $n\in \mathbb{N}^*$ and also $N_\mathcal{I}\equiv 0$, by classical direct limit argument; thus $\mathcal{I}$ is formally integrable but not integrable in general. On the opposite, if  each $M_n$ is a para-K\"ahler manifold, then, again, the Nijenhuis tensor of the almost para-complex structure $\mathcal{J}_n$  is zero and by Theorem \ref{T_CaracterizationOfIntegrableAlmostTangentParaComplexDecomposableComplexStructure}, this implies that the almost para-complex structure $\mathcal{J}$ on $\mathsf{L}_{k}^{p}(\mathbb{S}^{1},M)$ is also integrable.
\end{remark}

\section{Projective limits of adapted connections to $G$-structures}
\label{_ProjectiveLimitsOfAdaptedConnectionsToGStructures}

Let $\left(  F,\pi_{|F},M,G\right)  $ be a Banach $G$-structure, reduction of
the generalized tangent frame bundle $\ell\left(  TM\right)  =\left(  P\left(
TM\right)  ,\pi,M,\operatorname{GL}\left(  \mathbb{M}\right)  \right)  $.

\begin{definition}
\label{D_AdaptedConnectionToGStructure}A principal connection $\Phi$ on the
principal bundle $\ell\left(  TM\right)  $ is called an \textit{adapted
connection to the }$G$\textit{-structure} $\left(  F,\pi_{|F},M,G\right)  $ if
its Lie algebra valued connection form $\omega$ takes values in the Lie
algebra $\mathcal{G}$ of the Lie group $G$.
\end{definition}

\begin{theorem}
\label{T_FrechetConnectionCorrespondingToDescendingSequenceOfAdaptedConnections}%
Let $\left(  F_{n},\mathbf{p}_{n|F_{n}},M_{n},\mathcal{K}_{0}^{n}\left(
\mathbb{M}\right)  \right)  _{n\in\mathbb{N}^{\ast}}$ be a sequence of
$G$-structures, where $F_{n}$ is a $H_{n}$-reduction of the generalized
tangent frame bundle $\mathbf{\ell}\left(  TM_{n}\right)  =\left(
\mathbf{P}\left(  TM_{n}\right)  ,\mathbf{p}_{n},M_{n},\mathcal{H}_{0}^{n}\left(
\mathbb{M}\right)  \right)  $ where $\left(  T\mu_{n}^{n+1},\mu_{n}%
^{n+1},\gamma_{n}^{n+1}\right)  $ is a principal bundle morphism. \newline Let
$\left(  \omega_{n}\right)  _{n\in\mathbb{N}^{\ast}}$ be a sequence of
connections where, for any $n\in\mathbb{N}^{\ast}$, $\omega_{n}$ is adapted to
the $G$-structure of $F_{n}$. \newline Assume that for $i\leq j$, we have:%
\begin{equation}
\left(  T\mu_{i}^{j}\right)  ^{\ast}\left(  \omega_{i}\right)  =\overline
{\gamma_{i}^{j}}.\omega_{j}
\label{Rel_CompatibilityAdaptedDescendingSequenceOfConnections}%
\end{equation}
where $\overline{\gamma_{i}^{j}}:\mathfrak{k}_{j}\longrightarrow
\mathfrak{k}_{i}$ is the morphism of Lie algebras induced by the morphism
$\gamma_{i}^{j}:\mathcal{K}_{0}^{j}\left(  \mathbb{M}\right)  \longrightarrow
\mathcal{K}_{0}^{i}\left(  \mathbb{M}\right)  .$\newline Then $\omega
=\underleftarrow{\lim}\omega_{n}\in\Omega^{1}\left(  \underleftarrow{\lim
}F_{n},\underleftarrow{\lim}\mathfrak{k}_{n}\right)  $ is a connection form
adapted to the convenient bundle $\left(  \underleftarrow{\lim}F_{n}%
,\underleftarrow{\lim}\pi_{n|F_{n}},\underleftarrow{\lim}M_{n},\mathcal{K}%
_{0}^{n}\left(  \mathbb{M}\right)  \right)  $.
\end{theorem}

\section{Direct limits of adapted connections to $G$-structures}
\label{_DirectLimitsOfAdaptedConnectionsToGStructures}

\begin{theorem}
\label{T_ConvenientConnectionCorrespondingToAscendingSequenceOf AdaptedConnections}%
Let $\left(  F_{n},\mathbf{p}_{n|F_{n}},M_{n},H_{n}\right)  _{n\in\mathbb{N}^{\ast}}$
be the sequence of associated $G$-structures, where $F_{n}$ is a $H_{n}%
$-reduction of the tangent frame bundle $\ell\left(  TM_{n}\right)  =\left(
P_{n}\left(  TM_{n}\right)  ,\mathbf{p}_{n},M_{n},\operatorname{GL}\left(
\mathbb{M}_{n}\right)  \right)  $, where $\left(  \lambda_{n}^{n+1},\mu
_{n}^{n+1},\gamma_{n}^{n+1}\right)  $ is a principal bundle morphism. \newline
Let $\left(  \omega_{n}\right)  _{n\in\mathbb{N}^{\ast}}$ be a sequence of
connections where, for any $n\in\mathbb{N}^{\ast}$, $\omega_{n}$ is adapted to
the $G$-structure of $F_{n}$. Assume that for $i\leq j$, we have:%
\begin{equation}
\left(  \lambda_{i}^{j}\right)  ^{\ast}\left(  \omega_{j}\right)
=\overline{I_{i}^{j}}.\omega_{i}
\label{Rel_CompatibilityAdaptedAscendingSequenceOfConnections}%
\end{equation}
where $\overline{I_{i}^{j}}:\mathfrak{h}_{i}\longrightarrow\mathfrak{h}_{j}$
is the morphism of Lie algebras induced by the morphism $I_{i}^{j}%
:H_{i}\longrightarrow H_{j}.$\newline Then $\omega=\underrightarrow{\lim
}\omega_{n}\in\Omega^{1}\left(  \underrightarrow{\lim}F_{n},\underrightarrow
{\lim}\mathfrak{h}_{n}\right)  $ is a connection form adapted to the
convenient bundle $\left(  \underrightarrow{\lim}F_{n},\underrightarrow{\lim
}\mathbf{p}_{n|F_{n}},\underrightarrow{\lim}M_{n},H\left(  \mathbb{M}\right)
\right)  $.
\end{theorem}

\begin{proof}
We use the results of Theorem
\ref{T_ConvenientStructureOnDirectLimitOfAscendingSequenceOfGStructures}, the
Proposition \ref{P_SC_LieAlgebraValuedConnextion} and the relations of
compatibility (\ref{Rel_CompatibilityAdaptedAscendingSequenceOfConnections}%
).\newline The details are left to the reader.
\end{proof}

\appendix
\section{The convenient framework}
\label{_ConvenenientFramework}

\subsection{Convenient calculus}
\label{__ConvenientCalculus}

The convenient calculus discovered by A. Fr\"{o}licher and A. Kriegl about
1982 gives an adapted framework for differentiation in non-normable spaces
(cf. \cite{FrKr}) which coincides with the classical G\^{a}teaux approach on
Fr\'{e}chet spaces.

The references for this section are the tome \cite{KrMi} which includes some
further results and the papers \cite{CaPe}, \cite{FiTe} and \cite{Gra}.

To develop differential calculus on a vector spaces $E$ endowed with a locally
convex topology $\tau$, the basic idea is to test smoothness along smooth
curves since this notion is a concept without problems.

\begin{definition}
\label{D_SmoothCurve} Let $E$ be a locally convex vector space. A curve
$c:\mathbb{R}\longrightarrow E$ is called smooth if, for all $t$, the
derivative $c^{\prime}\left(  t\right)  $ exists where
\[
c^{\prime}\left(  t\right)  =\underset{h\longrightarrow0}{\lim}\dfrac{1}%
{h}\left(  c\left(  t+h\right)  -c\left(  t\right)  \right)  .
\]
It is \textit{smooth} if all iterative derivatives exist.
\end{definition}

\begin{notation}
\label{N_SetSmoothCurves}The set of smooth curves is denoted by $C^{\infty
}\left(  \mathbb{R},E\right)  $.
\end{notation}

The notion of smoothness is a bornological concept, i.e., $C^{\infty}\left(
\mathbb{R},E\right)  $ does not depend on the locally convex topology on $E$
but only on its associated \textit{bornology}. This can be deduced from the
characterization of smoothness in terms of $Lip^{k}$ for any $k$ (cf.
Proposition \ref{P_LipschitzSmoothCurves}). So smoothness can be translated in
terms of Lipschitz conditions which involves bounded sets.

In any locally convex space there is a natural notion of bounded sets called
the \textit{von Neumann bornology}.

Using the following classical theorem of Mackey in functional analysis, one
can deduce that if one varies the topology of $E$ while keeping the same
topological dual $E^{\prime}$, the bounded sets remain the sames.

\begin{theorem}
\label{T_CharacterizationBoundedSubsetByImageByContinuousLinearMappings}A
subset $B$ of a locally convex space $E$ is called bounded if $\ell\left(
B\right)  $ is bounded for all $\ell\in E^{\prime}$.
\end{theorem}

A multilinear map $E_{1}\times\cdots\times E_{n}\longrightarrow F$ is bounded
if bounded sets $B_{1}\times\cdots\times B_{n}$ are mapped onto bounded
subsets of $F$.

Continuous linear functional are clearly bounded.

\begin{notation}
\label{N_SetLCSBoundedLinearOperators}The locally convex vector space of
bounded linear operators with uniform convergence on bounded sets is denoted
by $L^{\times}\left(  E,F\right)  $.
\end{notation}

\begin{definition}
\label{D_Lipschitz}Let $E$ be a locally convex vector space. If $J$ is an open
subset of $\mathbb{R}$, a curve $c:\mathbb{R}\longrightarrow E$ is called
\textit{Lipschitz} on $J$ if the set
\[
\{\dfrac{c\left(  t_{2}\right)  -c\left(  t_{1}\right)  }{t_{2}-t_{1}}%
;t_{1},t_{2}\in J,t_{1}\not =t_{2}\}
\]
is bounded in $E$. The curve $c$ is \textit{locally Lipschitz} if every point
in $\mathbb{R}$ has a neighbourhood on which $c$ is Lipschitz. For
$k\in\mathbb{N}$, the curve $c$ is of class $Lip^{k}$ if $c$ is derivable up
to order $k$, and if the $k^{th}$-derivative $c:\mathbb{R}\longrightarrow E$
is locally Lipschitz.
\end{definition}

We then have the following link between both these notions (\cite{KrMi}
section 1.2):

\begin{proposition}
\label{P_LipschitzSmoothCurves} Let $E$ be a locally convex vector space and
let $c:\mathbb{R}\longrightarrow E$ be a curve. Then $c$ is $C^{\infty}$ if
and only if $c$ is $Lip^{k}$ for all $k\in\mathbb{N}$.
\end{proposition}

Note that the topology can vary considerably without changing the bornology;
the \textit{bornologification} $E_{\text{born}}$ of $E$ is the finest locally
convex structure having the same bounded sets.

Thanks to these objects, we can define another topology on $E$.

\begin{definition}
\label{D_c-InfinityOpen} A subset $U\subset E$ is called $c^{\infty}$-open if
$c^{-1}\left(  U\right)  $ is open in $\mathbb{R}$ for all $c\in C^{\infty
}\left(  \mathbb{R},E\right)  $.\newline The generated topology is called
$c^{\infty}$-topology and $E$ equipped with this topology is denoted by
$c^{\infty}E$.
\end{definition}

The $c^{\infty}$-topology is in general finer than the original topology and
$E$ is not a topological vector space when equipped with the $c^{\infty}$-topology.

\begin{remark}
\label{R_CoincidenceTopologiesFrechetSoBanachSpaces} For Fr\'{e}chet spaces
and so Banach spaces, this topology coincides with the given locally convex topology.
\end{remark}

For every absolutely convex closed bounded set $B$, the linear span $E_{B}$,
spanned by $B$ in $E$, is equipped with the Minkowski functional
\[
p_{B}\left(  v\right)  =\inf\left\{  \lambda>0:v\in\lambda.B\right\}
\]
which is a norm on $E_{B}$.

We then have the following characterization of $c^{\infty}$-open sets
(\cite{KrMi} Theorem 2.13):

\begin{proposition}
\label{P_c-InfinityOpen} $U\subset E$ is $c^{\infty}$-open if and only if
$U\cap E_{B}$ is open in $E_{B}$ for all absolutely convex bounded subsets
$B\subset E$.
\end{proposition}

The notion of smooth mapping between locally convex spaces $E$ and $F$ given
in the following definition corresponds to the usual notion when $E$ and $F$
are finite dimensional (cf. \cite{Bom}).

\begin{definition}
\label{D_SmoothMapping}Let $E$ and $F$ be locally convex spaces. A mapping
$f:E\longrightarrow F$ is called smooth if it maps smooth curves to smooth
curves:%
\[
\forall c\in C^{\infty}\left(  \mathbb{R},E\right)  ,\ f\circ c\in C^{\infty
}\left(  \mathbb{R},F\right)
\]

\end{definition}

\begin{notation}
The space of smooth mappings from $E$ to $F$ is denoted by $c^{\infty}\left(
E,F\right)  $.
\end{notation}

In order to define a category of adapted differential structure for direct
limits of geometric objects, we need a weak completeness assumption linked
with the notion of boundness.

\begin{definition}
\label{D_MackeyCauchySequence} A sequence $\left(  x_{n}\right)  $ in $E$ is
called Mackey-Cauchy if there exists a bounded absolutely convex subset $B$ of
$E$ such that $\left(  x_{n}\right)  $ is a Cauchy sequence in the normed
space $E_{B}$.
\end{definition}

\begin{definition}
\label{D_ConvenientVectorSpace} A locally convex vector space is said to be
$c^{\infty}$-complete or \textit{convenient } if any Mackey-Cauchy sequence
converges ($c^{\infty}$-completeness).
\end{definition}

\begin{definition}
Let $E$ be a locally convex space. A curve $c:\mathbb{R}\longrightarrow E$ is
said to be weakly smooth if $\lambda\circ c$ is smooth for all $\lambda\in
E^{\prime}.$
\end{definition}

We then have the following characterizations of a convenient vector space.

\begin{proposition}
A locally convex vector space is \textit{convenient if one of the following
equivalent conditions is satisfied:}

\begin{enumerate}
\item A curve $c:\mathbb{R}\longrightarrow E$ is smooth if and only if it is
weakly smooth;

\item For every absolutely convex closed bounded set $B$ the linear span
$E_{B}$ of $B$ in $E$, equipped with the norm $p_{B}$ is complete.

\item Any Lipschitz curve in $E$ is locally Riemann integrable.
\end{enumerate}

\begin{example}
\label{Ex_RInfinity}${}$The vector space $\mathbb{R}^{\infty},$ also denoted
by $\mathbb{R}^{\left(  \mathbb{N}\right)  }$ or $\Phi$, of all finite
sequences is a countable convenient vector space (\cite{KrMi}, 47.1) which is
not metrizable. $\mathbb{R}^{\infty}=$ $\coprod\limits_{n\in\mathbb{N}^{\ast}%
}\mathbb{R}^{n}$ is endowed with the finest locally convex topology for which
the inclusions $\varepsilon_{n}:\mathbb{R}^{n}\longrightarrow\coprod
\limits_{n\in\mathbb{N}^{\ast}}\mathbb{R}^{n}$ are continuous. A subset $U$ of
$\mathbb{R}^{\infty}$ is open for this topology if and only if the
intersections of $U$ with the finite dimensional subspaces of $\mathbb{R}%
^{\infty}$ are open in these subspaces. A basis of neighbourhoods of $0$
corresponds to the sets $B_{0}^{\infty}\left(  \varepsilon_{n}\right)
=\left\{  \left(  x_{n}\right)  _{n\in\mathbb{N}^{\ast}}\in\mathbb{R}^{\infty
}:\left\vert x_{n}\right\vert <\varepsilon_{n}\right\}  $ where $\left(
\varepsilon_{n}\right)  _{n\in\mathbb{N}^{\ast}}$ is a sequence of non
negative reals (\cite{Jar}, 4.1.4).
\end{example}
\end{proposition}

\begin{proposition}
\label{P_ConvenientStructureSetOfBoundedLinearOperators}If $E$ and $F$ are
convenient then the space $L^{\times}\left(  E,F\right)  $ is also convenient.
\end{proposition}

The products and limits of convenient vector spaces are convenient. Moreover,
all Fr\'{e}chet spaces are convenient.\newline In general, an inductive limit
of $c^{\infty}$-complete spaces need not be $c^{\infty}$-complete (cf.
\cite{KrMi}, 2.15, example).

The following theorem shows that the differential calculus on convenient
vector spaces is an appropriate extension of analysis on such spaces.

\begin{theorem}
\label{T_PropertiesDifferentialConvenientSpaces}Let $U$ be a $c^{\infty}$-open
set of a convenient vector space $E$ and let $F$ and $G$ be convenient vector spaces.

\begin{enumerate}
\item The space $c^{\infty}\left(  U,F\right)  $ may be endowed with a
structure of convenient vector space. The subspace $L^{\times}\left(
E,F\right)  $ is closed in $c^{\infty}\left(  E,F\right)  $;

\item The category is cartesian closed, i.e. we have the natural
diffeomorphism:%
\[
C^{\infty}\left(  E\times F,G\right)  \simeq C^{\infty}\left(  E,C^{\infty
}\left(  F,G\right)  \right)
\]

\item The differential operator%
\[
d:C^{\infty}\left(  E,F\right)  \longrightarrow C^{\infty}\left(  E,L\left(
E,F\right)  \right)
\]%
\[
df\left(  x\right)  v=\underset{t\longrightarrow0}{\lim}\dfrac{f\left(
x+tv\right)  -f\left(  x\right)  }{t}%
\]
exists and is linear and smooth.

\item The chain rule holds:%
\[
d\left(  f\circ g\right)  \left(  x\right)  v=df\left(  g\left(  x\right)
\right)  dg\left(  x\right)  v.
\]

\end{enumerate}
\end{theorem}

\subsection{Convenient geometric structures}
\label{__ConvenientGeometricStructures}

\subsubsection{Convenient manifolds}
\label{___ConvenientManifolds}

According to  \cite{KrMi} section  27.1, a $C^\infty$-{\it atlas} for a set $M$  modelled on a convenient space $X$ is a family $\{(U_\alpha, u_\alpha) \}_{\alpha\in A}$ of subsets $U_\alpha$ of $M$ and maps $u_\alpha$ from $U_\alpha$ to $X$  such that:

\begin{enumerate}
\item[--]
$u_\alpha$ is a bijection of $U_\alpha$ onto a $c^\infty$-open  subset of  $X$ for all $\alpha\in A$;

\item[--]
$M=\displaystyle\bigcup_{\alpha\in A}U_\alpha$;

\item[--]For any $\alpha$ and $\beta$ such that  $U_{\alpha\beta}=U_\alpha\cap U_\beta\not=\emptyset$,\\
$u_{\alpha\beta}=u_\alpha\circ u_\beta^{-1}:u_\beta(U_{\alpha\beta})\longrightarrow u_\alpha(U_{\alpha\beta})$ is a conveniently  smooth map.
\end{enumerate}

Classically, we have a notion of equivalent  $C^\infty$-atlases on  $M$. An equivalent class of $C^\infty$-atlases on $M$ is a maximal  $C^\infty$-atlas.
Such an atlas defines a topology on $M$ which is not in general Hausdorff.

\begin{definition}
\label{D_nnHConvenientManifold} 
A maximal  $C^\infty$-atlas on $M$  is called a non necessary Hausdorff convenient manifold structure on $M$ (n.n.H. convenient manifold $M$ for short); it is called a Hausdorff convenient manifold structure on $M$  when the topology defined by this  atlas is a Hausdorff topological space.
\end{definition}
Following the classical framework, when $X$ is a Banach space (resp. a Fr\'echet space)  we  say that $M$ is a \textit{Banach manifold} (resp. \textit{Fr\'echet manifold}) if $M$ is provided with a $C^\infty$- atlas (modelled on $X$) which generates a Hausdorff topological space.\\

The notion of vector bundle modelled on a convenient space over a n.n.H. convenient manifold  is defined in a classic way (cf.  (\cite{KrMi}, 29). Note that since a convenient space is Hausdorff, a vector bundle modelled on a convenient space has a natural structure of n.n.H. convenient manifold which is Hausdorff if and only if the base is a Hausdorff convenient manifold.

\subsubsection{Convenient principal bundles}
\label{__ConvenientPrincipalBundles}

\begin{definition}
\label{D_ConvenientPrincipalBundle}A convenient principal bundle is a
quadruple $\left(  P,\pi,M,G\right)  $ where $P$ and $M$ are convenient
manifolds, $\pi:P\longrightarrow M$ is a smooth map and $G$ is a convenient
Lie group acting on $P$ on the right via $%
\begin{array}
[c]{cccc}%
\widehat{R}: & P\times G & \longrightarrow & P\\
& \left(  p,g\right)  & \longmapsto & p\cdot g
\end{array}
$ such that: \newline For every $x\in$ $M$, there is an open $U$ with $x\in U$
and a diffeomorphism $\tau:$ $P_{|U}=\pi^{-1}\left(  U\right)  \longrightarrow
U\times G$ (local trivialization) satisfying the following properties:

\begin{description}
\item[\textbf{(PB 1)}] $\operatorname*{pr}_{1}\circ\tau=\pi$ where
$\operatorname*{pr}_{1}:U\times G\longrightarrow U$ is the first projection;

\item[\textbf{(PB 2)}] $\forall y\in U,\forall\left(  g,g^{\prime}\right)  \in
G^{2},\left(  \tau^{-1}\left(  y,g\right)  \right)  \cdot g^{\prime}=\tau
^{-1}\left(  y,g.g^{\prime}\right)  $.
\end{description}
\end{definition}

Given a trivializing cover $\left\{  \left(  U_{\alpha},\tau_{\alpha}\right)
\right\}  _{\alpha\in A}$ of the principal bundle $\left(  P,\pi,M,G\right)
$, the transition maps of this bundle are the smooth maps%
\[%
\begin{array}
[c]{cccc}%
\mathrm{T}_{\alpha\beta}: & U_{\alpha}\cap U_{\beta} & \longrightarrow & G\\
& x & \longmapsto & \left(  \overline{\tau}_{\alpha,x}\circ\left(
\overline{\tau}_{\beta,x}\right)  ^{-1}\right)  \left(  e\right)
\end{array}
\]

A principal bundle is completely determined by its cocycles (cf. \cite{KrMi}, 37.8).

For a cocycle $\left\{  g_{\alpha\beta}:U_{\alpha}\cap U_{\beta}%
\longrightarrow G\right\}  _{\left(  \alpha,\beta\right)  \in A^{2}}$ over an
open cover $\left\{  U_{\alpha}\right\}  _{\alpha\in A}$ of $M$, we consider
the set%
\[
S=\bigcup\limits_{\alpha\in A}\left(  \left\{  \alpha\right\}  \times
U_{\alpha}\times G\right)
\]

and define the equivalence relation:%
\[
\left(  \alpha,x,g\right)  \sim\left(  \beta,y,h\right)  \Longleftrightarrow
\left\{
\begin{array}
[c]{c}%
x=y\\
h=g_{\beta\alpha}\left(  x\right)  .y
\end{array}
.\right.
\]

We then obtain a principal bundle $\left(  P,\pi,M,G\right)  $ where
$P=S/\sim$ with the projection $\pi\left(  \widetilde{\left(  \alpha
,x,g\right)  }\right)  =x$. The local trivializations $\tau_{\alpha}:\pi
^{-1}\left(  U_{\alpha}\right)  \longrightarrow U_{\alpha}\times G$ are given
by:%
\[
\tau_{\alpha}\left(  \widetilde{\left(  \beta,x,g\right)  }\right)  =\left(
x,g_{\alpha\beta}\left(  x\right)  .g\right)  \text{.}%
\]

\begin{definition}
\label{D_MorphismConvenientPrincipalBundles}A morphism between the principal
bundles $\left(  P_{1},\pi_{1},M_{1},G_{1}\right)  $ and $\left(  P_{2}%
,\pi_{2},M_{2},G_{2}\right)  $ is a triple $\left(  \lambda,\mu,\gamma\right)
$ where $\lambda:P_{1}\longrightarrow P_{2}$ and $\mu:M_{1}\longrightarrow
M_{2}$ are smooth maps and $\gamma:G_{1}\longrightarrow G_{2}$ a morphism of
convenient Lie groups satisfying the following conditions:

\begin{description}
\item[\textbf{(PBM 1)}] $\pi_{2}\circ\lambda=\mu\circ\pi_{1};$

\item[\textbf{(PBM 2)}] $\forall p\in P_{1},\forall g\in G_{1},\ \lambda
\left(  p.g\right)  =\lambda\left(  p\right)  .\gamma\left(  g\right)  $.
\end{description}
\end{definition}

\subsubsection{Connections on a convenient principal bundle}
\label{___ConnectionsOnAConvenientPrincipalBundle}

Let $\left(  P,\pi,M,G\right)  $ be a convenient principal fibre bundle. We
denote the Lie algebra of the Lie group $G$ by $\mathcal{G}$ and its neutral
element by $e$.

\begin{definition}
\label{D_Connection}A \textit{connection} on the principal bundle $\left(
P,\pi,M,G\right)  $ is a fiber projection $\Phi:TP\longrightarrow VP$ viewed
as a vector valued $1$-form $\Phi\in\Omega^{1}\left(  P,VP\right)
\subset\Omega^{1}\left(  P,TP\right)  $ such that $\Phi\circ\Phi=\Phi$ and
$\operatorname*{im}\Phi=VP$.
\end{definition}

\begin{definition}
\label{D_PrincipalConnection}A principal connection is a connection on $P$
which is $G$-equivariant for the principal right action $\widehat{R}:P\times
G\longrightarrow G$%
\[
\forall g\in G,\ T\widehat{R}_{g}\left(  \Phi\right)  =\Phi.T\widehat{R}_{g}%
\]

\end{definition}

Because the vertical subbundle of $P$ is trivialized as a vector bundle
$VP\simeq P\times\mathcal{G}$ over $P$, we get a $\mathcal{G}$-valued $1$-form
$\omega\in\Omega^{1}\left(  P,\mathfrak{g}\right)  $ defined, for any $p\in
P$, by
\[
\omega\left(  X_{p}\right)  =T_{e}\left(  \widehat{R^{p}}\right)  ^{-1}\left(
X_{p}\right)
\]
where $\widehat{R^{p}}\left(  g\right)  =p\cdot g$. \newline$\omega$ is called
the \textit{Lie algebra valued connection form}.\newline We then have:%
\[
\Phi\left(  X_{p}\right)  =\zeta_{\omega\left(  X_{p}\right)  }\left(
p\right)
\]

where the \textit{fundamental vector field} $\zeta_{X}$ is defined for any
$X\in\mathcal{G}$ by%
\[
\forall p\in P,\ \zeta_{X}\left(  p\right)  =T_{\left(  p,e\right)  }%
\widehat{R}\left(  \mathbf{0}_{p},X\right)  \text{.}%
\]

\begin{proposition}
Let $\Phi\in\Omega^{1}\left(  P,VP\right)  $ be a principal connection on the
principal bundle $\left(  P,\pi,M,G\right)  .$Then $\omega$ fulfils the relations:

\begin{description}
\item[\textbf{(VLACF 1)}] {$\forall X\in\mathfrak{g},\ \omega\left(  \zeta
_{X}\left(  p\right)  \right)  =X$~; }

\item[\textbf{(VLACF 2)}] {$\forall X\in\mathfrak{g},\forall X_{p}\in
T_{p}P,\ \left(  \widehat{R_{g}}^{\ast}\omega\right)  \left(  X_{p}\right)
=\operatorname*{Ad}\left(  g^{-1}\right)  .\omega\left(  X_{p}\right)  $~; }

\item[\textbf{(VLACF 3)}] {$\forall X\in\mathfrak{g},\ \mathcal{L}_{\zeta_{X}%
}\omega=-\operatorname*{ad}\left(  X\right)  .\omega$. }
\end{description}
\end{proposition}

The following proposition is used to get a structure of convenient connection
on a direct limit of some principal connections on Banach bundles.

\begin{proposition}
\label{P_SC_LieAlgebraValuedConnextion}A $1$-form $\omega\in\Omega^{1}\left(
P,\mathfrak{g}\right)  $ fulfilling condition \emph{(VLACF 1) }defines a
connection $\Phi$ on $P$ as%
\[
\Phi\left(  X_{p}\right)  =T_{e}\widehat{R^{p}}.\omega\left(  X_{p}\right)  .
\]
This connection is a principal connection if and only if the condition
\emph{(VLACF 2)} is satisfied.
\end{proposition}

\section{Linear operators}
\label{_LinearOperators}

\subsection{Operators from a Banach space to its dual and bilinear forms}
\label{__OperatorsFromBanachSpaceToItsDual}

Let $\mathbb{E}$ be a Banach space equipped with a norm $\left\Vert
\ \right\Vert $.

In this appendix, the topological dual will be denoted $\mathbb{E}^{\ast}$
(instead of $\mathbb{E}^{\prime}$). It will be equipped with the associated
norm $\left\Vert \ \right\Vert ^{\ast}$defined by%
\begin{equation}
\left\Vert \alpha\right\Vert ^{\ast}=\sup_{\left\Vert u\right\Vert
=1}\left\vert \alpha.u\right\vert \label{eq_NormOnDuaOfBanachSpace}%
\end{equation}

Given any continuous injective operator $A:\mathbb{E}\longrightarrow
\mathbb{E}^{\ast}$, we consider the norm $\left\Vert \ \right\Vert _{A}$ on
$\mathbb{E}$ defined by%
\[
\left\Vert u\right\Vert _{A}=\left\Vert A.u\right\Vert ^{\ast}=\sup
_{\left\Vert v\right\Vert =1}\left\vert \left\langle A.u,v\right\rangle
\right\vert =\sup_{v\in\mathbb{E}\setminus\left\{  0\right\}  }\dfrac
{\left\vert \left\langle A.u,v\right\rangle \right\vert }{\left\Vert
v\right\Vert }%
\]

where $\left\langle \ ,\ \right\rangle $ is the duality bracket.

The operator norm of $A$ is given by%
\begin{equation}
\left\Vert A\right\Vert ^{\mathrm{op}}=\sup_{u\in\mathbb{E}\setminus\left\{
0\right\}  }\dfrac{\left\Vert A.u\right\Vert ^{\ast}}{\left\Vert u\right\Vert
}=\sup_{u\in\mathbb{E}\setminus\left\{  0\right\}  }\sup_{v\in\mathbb{E}%
\setminus\left\{  0\right\}  }\dfrac{\left\vert \left\langle
A.u,v\right\rangle \right\vert }{\left\Vert u\right\Vert \left\Vert
v\right\Vert } \label{eq_NormOperatorFromBanachSpaceToItsDual}%
\end{equation}

Because we have%
\[
\forall u\in\mathbb{E},\ \left\Vert u\right\Vert _{A}\leqslant\left\Vert
A\right\Vert ^{\mathrm{op}}\left\Vert u\right\Vert
\]
then the Identity map $\operatorname{Id}:\left(  \mathbb{E},\left\Vert
\ \right\Vert \right)  $ $\longrightarrow\left(  \mathbb{E},\left\Vert
\ \right\Vert _{A}\right)  $ is continuous.

Denote by $\widehat{\mathbb{E}_{A}}$ the completion of $\mathbb{E}$ for the
norm $\left\Vert \ \right\Vert _{A}$.

To such an operator $A$ is associated a canonical continuous bilinear form
$\mathcal{B}_A$ on $\mathbb{E}$ defined by%
\[
\mathcal{B}_A\left(  u,v\right)  =\left\langle A.u,v\right\rangle
\]

Conversely, any continuous bilinear form $B$ on $\mathbb{E}$ gives rise to the
following continuous linear mappings $B_{L}^{\flat}:\mathbb{E}\longrightarrow
\mathbb{E}^{\ast}$ and $B_{R}^{\flat}.u:\mathbb{E}\longrightarrow\mathbb{R}$ where%

\[%
\begin{array}
[c]{cccc}%
B_{L}^{\flat}.u: & \mathbb{E} & \longrightarrow & \mathbb{R}\\
& v & \mapsto & B\left(  u,v\right)
\end{array}
,%
\begin{array}
[c]{cccc}%
B_{R}^{\flat}.v: & \mathbb{E} & \longrightarrow & \mathbb{R}\\
& u & \mapsto & B\left(  u,v\right)
\end{array}
\]

\begin{definition}
\label{D_WeaklyStronglyNondegenerateBilinearForm} 
The continuous bilinear form $B:\mathbb{E}\times\mathbb{E}\longrightarrow\mathbb{R}$ is said to be left (resp. right) weakly non degenerate\index{bilinear form!weakly non degenerate}
if $B_{L}^{\flat}:\mathbb{E}\longrightarrow\mathbb{E}^{\ast}$ (resp.
$B_{R}^{\flat}:\mathbb{E}\longrightarrow\mathbb{E}^{\ast}$) is
injective.\newline 
If $B_{L}^{\flat}$ (resp. $B_{R}^{\flat}$) is also surjective, then $B$ is said to be left (resp. right) strongly non degenerate.
\end{definition}

Note that if $B$ is symmetric or skew-symmetric then $B$ is left weakly (resp.
strongly) non degenerate if and only if $B$ is right weakly (resp. strongly)
non degenerate. In this case, we will simply say that $B$ is weakly (resp.
strongly) non degenerate.

If we assume that $B$ is right and left weak non degenerate, as previously, we
can define two new norms $\left\Vert u\right\Vert ^{L}=\left\Vert B_{L}%
^{\flat}\left(  u\right)  \right\Vert ^{\ast}$ and $\left\Vert u\right\Vert
^{R}=\left\Vert B_{R}^{\flat}\left(  u\right)  \right\Vert ^{\ast}$ and denote
$\mathbb{E}_{L}$ (resp. $\mathbb{E}_{R}$) the normed space $\left(
E,\left\Vert \ \right\Vert ^{L}\right)  $ (resp. $\left(  E,\left\Vert
\ \right\Vert ^{R}\right)  $).

\begin{remark}
\label{R_IndependanceOfTheNorms}If $\widehat{\left\Vert \ \right\Vert }$ is an
equivalent norm of $\left\Vert \ \right\Vert $on $\mathbb{E}$, then the
corresponging norms $\widehat{\left\Vert \ \right\Vert }^{\ast}$ and
$\left\Vert \ \right\Vert ^{\ast}$ are also equivalent norms on $\mathbb{E}%
^{\ast}$ and so $\widehat{\left\Vert \ \right\Vert }^{L}$ and $\left\Vert
\ \right\Vert ^{L}$ (resp. $\widehat{\left\Vert \ \right\Vert }^{R}$ and
$\left\Vert \ \right\Vert ^{R}$). Therefore the completion of $\mathbb{E}_{L}$
(resp. $\mathbb{E}_{R}$) only depends on the Banach structure on $\mathbb{E}$.
\end{remark}

If $B$ is symmetric or skew-symmetric, then we have $\left\Vert \ \right\Vert
^{L}=\left\Vert \ \right\Vert ^{R}$ and $\mathbb{E}_{L}=\mathbb{E}_{R}$. In
this case, we set $\left\Vert \ \right\Vert _{B}=\left\Vert \ \right\Vert
^{L}=\left\Vert \ \right\Vert ^{R}$ and the associate Banach space will be
denoted $\mathbb{E}_{B}$.

\begin{proposition}
\label{P_ExtensionLinearMappingLinearSpaceIntoDual}Let $A:\mathbb{E}%
\longrightarrow\mathbb{E}^{\ast}$ be an injective continuous operator.

\begin{enumerate}
\item There exists a unique continuous linear mapping $\widehat{A}%
:\widehat{\mathbb{E}_{A}}\longrightarrow\mathbb{E}^{\ast}$ whose restriction
to $\mathbb{E}$ is $A$.\newline Moreover, the range of $\widehat{A}$ is closed
and $\widehat{A}$ is an isometry
\index{isometry}
from $\widehat{\mathbb{E}_{A}}$ onto this range.

\item The associated bilinear form $\mathcal{B}_A$ is a continuous left weakly
non degenerate bilinear form which can be extended to a form of the same type
$\widehat{\mathcal{B}_A}$ on $\mathbb{E}\times\widehat{\mathbb{E}_{A}}$.\newline
Moreover, if $\mathcal{B}_A$ is symmetric (resp. skew-symmetric) then
$\widehat{\mathcal{B}_A}$ is symmetric (resp. skew-symmetric).

\item If $\mathbb{E}$ is reflexive then $\widehat{A}$ is an isomorphism from
$\widehat{\mathbb{E}_{A}}$ to $\mathbb{E}^{\ast}$ and $\widehat{A}^{\ast}$ is
an isomorphism from $\mathbb{E}^{\ast\ast}\simeq\mathbb{E}$ to $\widehat
{\mathbb{E}_{A}}$.
\end{enumerate}
\end{proposition}

\subsection{Krein inner products}
\label{__KreinInnerProducts}

In this section, the reader is refereed to Bognar's book (\cite{Bog}).

\begin{definition}
\label{D_IndefiniteInnerProduct} Let $\mathbb{E}$ be a Banach space.

\begin{enumerate}
\item An indefinite inner product%
\index{indefinite inner product}
on $\mathbb{E}$ is a continuous weak symmetric bilinear form $g$ on
$\mathbb{E}$.

\item An element $u$ of $\mathbb{E}$ is called

\begin{enumerate}
\item positive if $g\left(  u,u\right)  \geqslant0$;

\item negative if $g\left(  u,u\right)  \leqslant0$;

\item neutral if $g\left(  u,u\right)  =0$.
\end{enumerate}
\end{enumerate}
\end{definition}

We denote by

\begin{enumerate}
\item[--] $\mathbb{E}^{+}$ the set of positive elements of $\mathbb{E}$;

\item[--] $\mathbb{E}^{-}$ the set of negative elements of $\mathbb{E}$;

\item[--] $\mathbb{E}^{0}$ the set of neutral elements of $\mathbb{E}$.
\end{enumerate}

We set $\mathbb{E}^{++}=\mathbb{E}^{+}\setminus\left\{  0\right\}  $,
$\mathbb{E}^{--}=\mathbb{E}^{-}\setminus\left\{  0\right\}  $ and
$\mathbb{E}^{00}=\mathbb{E}^{0}\setminus\left\{  0\right\}  $.

For $\left(  u,v\right)  \in\mathbb{E}^{2}$, we will denote $g\left(
u,v\right)  =\left(  u|v\right)  $.

If $\mathbb{E}^{++}$ and $\mathbb{E}^{-}$ are non empty sets, $g$ is called an
\textit{indefinite metric}%
\index{indefinite metric}%
. If $\mathbb{E}^{-}$ $=\emptyset$ (resp.$\mathbb{E}^{++}=\emptyset$) the $g$
is called an \textit{definite positive} (resp. \textit{define negative}). In
these cases, $\mathbb{E}^{00}=\emptyset$.

\begin{definition}
\label{D_OrthogonalSubspaceDefiniteMetricBanachSpace} 
Let $\mathbb{F}$ be a linear subspace of a Banach space $\mathbb{E}$ and let $g$ be an indefinite
inner product on $\mathbb{E}$. The orthogonal of $\mathbb{F}$ is the set
\[
\mathbb{F}^{\bot}=\left\{  u\in\mathbb{E}:\forall v\in\mathbb{F},g\left(
u,v\right)  =0\right\}
\]

\end{definition}

$\mathbb{F}^{\bot}$ is a closed linear subspace of $\mathbb{E}$. If $g$ is
indefinite, in general, $F\cap F^{\bot}\neq\emptyset$. For instance, if
$\mathbb{F}$ is the vector space generate by $\mathbb{E}^{+}$ or
$\mathbb{E}^{-}$, then $(\;|\;)_{F}$ is positive definite or definite negative
(cf. \cite{Bog}, Lemma 2.5). A subspace $\mathbb{F}$ is \textit{isotropic}%
\index{isotropic}
if $\mathbb{F}\subset\mathbb{F}^{\perp}$.

\begin{definition}
\label{D_KreinInnerProduct} Let $\mathbb{E}$ be a Banach space.

\begin{enumerate}
\item An indefinite inner product $g$ on $\mathbb{E}$ is called a Krein inner
product
\index{Krein inner product}
if there exists a decomposition $\mathbb{E}=\mathbb{E}^{+}\oplus\mathbb{E}%
^{-}$ of $\mathbb{E}$ such that the restriction of $g$ to $\mathbb{E}^{+}$
(resp. $\mathbb{E}^{-}$) is definite positive(resp. definite negative). For
each $u\in\mathbb{E}$, we set $u^{+}$ and $u^{-}$ the component of $u$ on
$\mathbb{E}^{+}$ and $\mathbb{E}^{-}$ respectively. 

\item a Krein inner product is called a neutral inner product \index{neutral inner product} if there exists an associated  decomposition  $\mathbb{E}=\mathbb{E}^{+}\oplus\mathbb{E}%
^{-}$ of $\mathbb{E}$  such that $\mathbb{E}^+$ and $\mathbb{E}^-$ are isomorphic.

\item An indefinite inner product $g$ on $\mathbb{E}$ is called a indefinite
Krein inner product if there exists a decomposition $\mathbb{E}=\mathbb{E}%
_{1}\oplus\mathbb{E}_{2}$ of $\mathbb{E}$ such that the restriction of $g$ to
$\mathbb{E}_{1}$ (resp. $\mathbb{E}_{2}$ ) is an indefinite inner product
where $\mathbb{E}_{1}$ is the $g$-orthogonal of $\mathbb{E}_{2}$.
\end{enumerate}
\end{definition}

\begin{remark}
\label{R_NonTrivialDecompositionKreinProduct} Any indefinite inner product $g$
on a Banach space need not have a non trivial decomposition as in Definition
\ref{D_KreinInnerProduct}. The reader will find such an example in \cite{Bog},
Example 11.3.
\end{remark}

Consider an inner Krein product $g$ on $\mathbb{E}$ and a decomposition
$\mathbb{E}=\mathbb{E}^{+}\oplus\mathbb{E}^{-}$. For each $u\in\mathbb{E}$, we
set $u^{+}$ and $u^{-}$ the component of $u$ on $\mathbb{E}^{+}$ and
$\mathbb{E}^{-}$ respectively. To $g$ is canonically associate a pre-Hilbert
product $\gamma$ on $\mathbb{E}$ given by
\[
\gamma_{\mathbb{E}^{+}}=g_{\mathbb{E}^{+}},\quad\gamma_{\mathbb{E}^{-}%
}=-g_{\mathbb{E}^{-}}
\]

where the ${\gamma}$-orthogonal of $\mathbb{E}^{+}$ is $\mathbb{E}^{-}$.

\bigskip We denote by ${\mathbb{E}}_{g}$ the Hilbert space generated by the
pre-Hilbert product ${\gamma}$, i.e. the ${\gamma}$-dual of $\mathbb{E}$. We
then have an orthogonal decomposition ${\mathbb{E}}_{g}={\mathbb{E}}_{g}%
^{+}\oplus{\mathbb{E}}_{g}^{-}$, $\mathbb{E}^{+}=\mathbb{E}\cap{\mathbb{E}%
}_{g}^{+}$, $\mathbb{E}^{-}=\mathbb{E}\cap{\mathbb{E}}_{g}^{-}$. Then
$\mathbb{E}=\mathbb{E}^{+}\oplus\mathbb{E}$. \newline If we consider the map
\[%
\begin{array}
[c]{cccc}%
J: & \mathbb{E} & \longrightarrow & \mathbb{E}\\
& u=u^{+}+u^{-} & \mapsto & u^{+}-u^{-}%
\end{array}
\]
we have $g(u,v)={\gamma}(u,Jv)$ and ${\gamma}(u,v)=g(u,Jv)$. Moreover
$\gamma(u,Jv)=\gamma(Ju,v).$

$J$ is called the \textit{fundamental symmetry}%
\index{fundamental symmetry}
of $g$ and $\gamma$ the \textit{fundamental Hilbert product}%
\index{fundamental inner product}
of $g$.

From the previous subsection, $g$ and $\gamma$ can be extended to symmetric
bilinear forms $\widehat{g}$ and $\widehat{\gamma}$ on $\mathbb{E}%
\times\mathbb{E}_{g}$. In fact, we have

\begin{proposition}
\label{P_ExtensionKreinInnerProduct} 
The Krein inner product $g$ can be extended to a continuous Krein inner product $\widehat{g}$ on $\mathbb{E}_{g}$ and $\widehat{g}$ is strongly non degenerate. \newline Conversely, let
$\widehat{g}$ be a Krein inner product on a Hilbert space $\mathbb{H}$. For
any Banach space which is embedded in $\mathbb{H}$ by a linear continuous map
$\iota:\mathbb{E}\rightarrow\mathbb{H}$, $g=\iota^{\ast}\widehat{g}$ is a
Krein inner product on $\mathbb{E}$.
\end{proposition}

We then have the following results.

\begin{corollary}
\label{C_IsomorphismKreinProducts} Let $\mathbb{E}$ be a Banach space.

\begin{enumerate}
\item Let $g_{1}$ and $g_{2}$ be two pre-Hilbert products on $\mathbb{E}$ such
that the norms $||\;||_{g_{1}}$ and $||\;||_{g_{1}}$ are equivalent on
$\mathbb{E}$. Then there exists a linear isomorphism $\varphi$ of $\mathbb{E}$
such that $\varphi^{\ast}g_{2}=g_{1}$.

\item Let $g$ and $g^{\prime}$ be two Krein inner products on $\mathbb{E}$
such that the norms $||\;||_{\gamma_{1}}$ and $||\;||_{\gamma_{2}}$ are
equivalent on $\mathbb{E}$ where $\gamma_{i}$ is the pre-hilbert product
associated to $g_{i}$, for $i=1,2$. If $\mathbb{E}=\mathbb{E}_{1}^{+}%
\oplus\mathbb{E}_{1}^{-}$ and $\mathbb{E}=\mathbb{E}_{2}^{+}\oplus
\mathbb{E}_{2}^{-}$ are the decompositions associated to $g_{1}$ and $g_{2}$
respectively such that $E_{1}^{+}$ and $E_{2}^{+}$ (resp. $E_{1}^{-}$ and
$E_{2}^{-}$) are isomorphic, then there exists an isomorphism $\varphi$ of
$\mathbb{E}$ such that $\varphi^{\ast}g_{2}=g_{1}$.
\end{enumerate}
\end{corollary}

\begin{definition}
\label{D_IsometryKreinInnerProduct} An isometry%
\index{isometry}
of a Krein inner product $g$ on a Banach space is an operator $I:\mathbb{E}%
\rightarrow\mathbb{E}$ such that
\[
\forall\left(  u,v\right)  \in\mathbb{E}^{2},\ g(Iu,Iv)=g(u,v)
\]

\end{definition}

From the Proposition \ref{P_ExtensionKreinInnerProduct}, it follows that such
an isometry is the restriction to $\mathbb{E}$ of an isometry $\widehat{I}$ of
$\widehat{g}$ on $\mathbb{E}_{g}$. According to a decomposition $\mathbb{E}%
_{g}=\mathbb{E}_{g}^{+}\oplus\mathbb{E}_{g}^{-}$, $\widehat{I}$ is a sum of an
isometry $\widehat{I}^{+}$ and $\widehat{I}^{-}$ of $\mathbb{E}_{g}^{+}$ and
$\mathbb{E}_{g}^{-}$ respectively. It follows that $\widehat{I}_{g}%
=\widehat{I}^{+}-\widehat{I}^{-}$ is an isometry of $\widehat{\gamma}$.
Conversely each isometry $\widehat{I}_{g}$of $\widehat{\gamma}$ can be
decomposed in a sum of an isometry $\widehat{I}_{g}^{+}$ and $\widehat{I}%
_{g}^{-}$ of $\mathbb{E}_{g}^{+}$ and $\mathbb{E}_{g}^{-}$ respectively and
$\widehat{I}=\widehat{I}_{g}^{+}-\widehat{I}_{g}^{-}$ is an isometry of $g$.

Finally we obtain:

\begin{proposition}
\label{P_IsometryKreinProduct} Let $g$ be a Krein inner product $g$ on
$\mathbb{E}$ and $\iota$ the densely embedding of $\mathbb{E}$ in the Hilbert
space $\mathbb{H}$. We consider a decomposition $\mathbb{H}=\mathbb{H}%
^{+}\oplus\mathbb{H}^{-}$.

\begin{enumerate}
\item Each element of the isometry group%
\index{isometry group}
of $\mathbb{H}$ can be written as a matrix $\left(
\begin{array}
[c]{cc}%
\widehat{I}^{+} & 0\\
0 & \widehat{I}^{-}%
\end{array}
\right)  $ where $\widehat{I}^{+}$ and $\widehat{I}^{-}$ is an isometry of
$\mathbb{H}^{+}$ and $\mathbb{H}^{-}$ respectively.

\item The isometry group of the Krein inner product $\widehat{g}$ is the set
of matrix $\left(
\begin{array}
[c]{cc}%
\widehat{I}^{+} & 0\\
0 & -\widehat{I}^{-}%
\end{array}
\right)  $ where $\widehat{I}^{+}$ and $\widehat{I}^{-}$ is an isometry of
$\mathbb{H}^{+}$ and $\mathbb{H}^{-}$ respectively.

\item If $\mathbb{E}^{+}=\mathbb{H}^{+}\cap\mathbb{E}$ and $\mathbb{E}%
^{+}=\mathbb{H}^{-}\cap\mathbb{E}$, then each isometry is of type $\iota
^{\ast}\widehat{I}$ where $\widehat{I}$ is an isometry of $\widehat{g}$.
\end{enumerate}
\end{proposition}

\subsection{Symplectic structures}

\label{__SymplecticStructures}

In this section, $\mathbb{E}$ is a Banach space equipped with a continuous
skew-symmetric bilinear form $\Omega:\mathbb{E}\times\mathbb{E\longrightarrow
R}$. Consider the associated continuous linear map $\Omega^{\flat}%
:\mathbb{E}\longrightarrow\mathbb{E}^{\ast}$ as defined in
\ref{__OperatorsFromBanachSpaceToItsDual}.

\begin{definition}
\label{D_WeakAndLinearSymplecticStructureOnBanachSpace}$\Omega^{\flat}$ is
called weakly non-degenerate if it is injective, and the pair $\left(
\mathbb{E},\Omega\right)  $ is called a weak linear symplectic structure.
\end{definition}

\begin{definition}
\label{D_StrongAndLinearSymplecticStructureOnBanachSpace}$\Omega^{\flat}$ is
called (strongly) non-degenerate if it is an isomorphism, and the pair
$\left(  \mathbb{E},\Omega\right)  $ is called a (strong) linear symplectic
structure%
\index{symplectic structure}%
.
\end{definition}

\begin{example}
\label{Ex_WeakSymplecticFormOnProductSpaceByDual} If $\mathbb{E}$ is any
Banach space, the product $\mathbb{E}\times\mathbb{E}^{\ast}$ is naturally
equipped with a weak symplectic structure defined by%
\[
\Omega\left(  \left(  u,\alpha\right)  ,\left(  v,\beta\right)  \right)
=\left\langle \alpha,v\right\rangle -\left\langle \beta,u\right\rangle
\]

\end{example}

\begin{remark}
The existence of a linear symplectic structure $\left(  \mathbb{E}%
,\Omega\right)  $ implies not only that $\mathbb{E}$ is isomorphic to
$\mathbb{E}^{\ast}$ but also that $\mathbb{E}$ is reflexive.
\end{remark}

Weinstein has treated lagrangian subspaces%
\index{lagrangian subspace}
of a symplectic Banach space in \cite{Wei}.

\begin{definition}
\label{D_IsotropicSpaceForSymplecticStructure}Let $\Omega$ be a symplectic
structure on the Banach space $\mathbb{E}$. \newline A subspace $\mathbb{F}$
is isotropic if
\[
\forall\left(  u,v\right)  \in\mathbb{F},\ \Omega\left(  u,v\right)  =0
\]

\end{definition}

An isotropic subspace is always closed.

If
\[
\mathbb{F}^{\bot}=\left\{  u\in\mathbb{E}:\forall v\in\mathbb{F},\Omega\left(
u,v\right)  =0\right\}
\]
is the symplectic orthogonal%
\index{symplectic orthogonal}
of the subspace $\mathbb{F}$, then $\mathbb{F}$ is isotropic if and only if
$\mathbb{F}\subset\mathbb{F}^{\bot}$.

\begin{definition}
\label{D_MaximalIsotropicSpaceForSymplecticStructure}Let $\Omega$ be a
symplectic structure on the Banach space $\mathbb{E}$. \newline A subspace
$\mathbb{F}$ is maximal isotropic if $\mathbb{F}=\mathbb{F}^{\bot}$.
\end{definition}

Unfortunately, in the Banach framework, a maximal isotropic subspace
$\mathbb{L}$ need not be supplemented. Following Weinstein's terminology
(\cite{Wei}), we introduce the notion of Lagrangian space.

\begin{definition}
\label{D_LagrangianSpaceForSymplecticStructure}Let $\Omega$ be a symplectic
structure on the Banach space $\mathbb{E}$. \newline An isotropic space
$\mathbb{L}$ is called a \textit{Lagrangian space} if there exists an
isotropic space $\mathbb{L}^{\prime}$ such that $\mathbb{E}=\mathbb{L}%
\oplus\mathbb{L}^{\prime}$.
\end{definition}

Since $\Omega$ is (strongly) non degenerate, this implies that $\mathbb{L}$
and $\mathbb{L}^{\prime}$ are maximal isotropic and then Lagrangian spaces
(see \cite{Wei}).

Unfortunately, in general, given a symplectic structure, Lagrangian subspaces
need not exist (cf. \cite{KaSw}). Even for a strong symplectic
structure on a Banach space which is not Hilbertizable the non existence of
Lagrangian subspace is an open problem to our knowledge.

\bigskip

Following \cite{Wei}, we introduce the fundamental notion of Darboux form.

\begin{definition}
\label{D_DarbouxForm}A symplectic form $\Omega$ on a Banach space $\mathbb{E}$
is a \textit{Darboux (linear) form}%
\index{Darboux form}
if there exists a Banach subspace $\mathbb{L}$ and an isomorphism
$A:\mathbb{E}\rightarrow\mathbb{L}\oplus\mathbb{L}^{\ast}$ such that
$\Omega=A^{\ast}\Omega_{\mathbb{L}}$ where $\Omega_{\mathbb{L}}$ is defined by%
\[
\Omega_{\mathbb{L}}\left(  \left(  u,\eta\right)  ,\left(  v,\xi\right)
\right)  =\left\langle \eta,v\right\rangle -\left\langle \xi,u\right\rangle
\]

\end{definition}

Note that, in this case, $\mathbb{E}$ must be reflexive.

\subsection{Cotangent structures}
\label{__CotangentStructures}

As in finite dimension, a cotangent structure is a particular case of
symplectic structure.

\begin{definition}
\label{D_CotangentStructure}A (linear) cotangent structure%
\index{cotangent structure}
on a Banach space $\mathbb{E}$ is a weak symplectic form $\Omega$ on
$\mathbb{E}$ such that there exists a maximal isotropic space $\mathbb{L}$ of
$\mathbb{E}$ giving rise to a decomposition $\mathbb{E}=\mathbb{L}+\mathbb{K}$
where $\mathbb{K}$ is isomorphic to $\mathbb{L}$. \newline Such a space will
be called a weak Lagrangian subspace.
\end{definition}

Note that $\mathbb{L}$ is not assumed to have a Lagrangian complement, even an
isotropic complement, which is equivalent to Lagrangian (cf. \cite{ThSc}).
Nevertheless, if $\mathbb{E}$ is a Hilbert space, then a cotangent structure
must be a (linear) Darboux form (see \cite{Wei}).

We do not know if there exists a cotangent structure which is not a Darboux
form\textit{.}

\begin{proposition}
\label{P_DarbouxFormCotangentStructure} If $\Omega$ is a linear Darboux form,
then $\Omega$ is a cotangent structure if we have a Lagrangian decomposition
$\mathbb{E}=\mathbb{L}\oplus\mathbb{L}^{\prime}$ where $\mathbb{L}$ and
$\mathbb{L}^{\prime}$ are isomorphic.
\end{proposition}

Then for any other complement subspace $\mathbb{K}$ of $\mathbb{L}$, then
$\mathbb{E}=\mathbb{L}\oplus\mathbb{K}$ is also an admissible decomposition
for $\Omega$. If $\mathbb{E}$ is a Hilbert space it is always true for any
Lagrangian decomposition (cf. \cite{Wei}).

\subsection{Tangent structures}

\label{__TangentStructures}

\begin{definition}
\label{D_TangentStructure}A (linear) tangent structure%
\index{tangent structure}
on a Banach space $\mathbb{E}$ is an endomorphism $J$ of $\mathbb{E}$ such
that
\[
\operatorname*{im}J=\ker J
\]
and there exists a decomposition $\mathbb{E}=\ker J+\mathbb{K}$.
\end{definition}

If $J$ is a tangent structure on $\mathbb{E}$, consider the decomposition
$\mathbb{E}=\ker J\oplus\mathbb{K}$. The restriction $J_{\mathbb{K}}$ of $J$
to $\mathbb{K}$ is an isomorphism onto $\ker J$. According to this
decomposition, $J$ can be identify with the matrix%
\[
\left(
\begin{array}
[c]{cc}%
0 & J_{\mathbb{K}}\\
0 & 0
\end{array}
\right)
\]

On the product of Banach spaces $\mathbb{E}^{2}=\mathbb{E}\times\mathbb{E}$,
we have a canonical tangent structure $J_{\mathrm{can}}$ defined by
$J_{\mathrm{can}}(u_{1},u_{2})=u_{1}$ and the restriction of $J_{\mathrm{can}%
}$ to $\mathbb{E}\times\{0\}$ vanishes. Consider $\mathbb{E}_{1}%
=\mathbb{E}\times\{0\}$ and $\mathbb{E}_{2}=\{0\}\times\mathbb{E}$ and let
$\iota_{i}:\mathbb{E}\rightarrow\mathbb{E}_{i}$ be the canonical isomorphism
from $\mathbb{E}$ to $\mathbb{E}_{i}$, for $i=1,2$. Then $J_{\mathrm{can}}$
can be identified with the matrix%

\[
\left(
\begin{array}
[c]{cc}%
0 & \iota_{1}\circ\operatorname{Id}_{\mathbb{E}}\circ\iota_{2}^{-1}\\
0 & 0
\end{array}
\right)
\]

\begin{proposition}
\label{P_IsomorphismDecompositionTangentStructure} 
If $\mathbb{E}=\ker J\oplus
K$ is a decomposition associated to a tangent structure $J$ on $\mathbb{E}$,
then $A:\mathbb{E} \longrightarrow\ker J\oplus\ker J$ is an isomorphism such
that $A^{*} J_{can}=J $ where $J_{can}$ is the canonical tangent structure on
$\ker J\oplus\ker J$.
\end{proposition}

\subsection{Complex and para-complex structures}
\label{__ComplexParaComplexStructures}

\subsubsection{Complex Structures}
\label{___ComplexStructures}

\begin{definition}
\label{D_ComplexStructure}
A complex  structure\index{complex structure}
on a Banach space $\mathbb{E}$ is an endomorphism $\mathcal{I}$ of
$\mathbb{E}$ such that $\mathcal{I}^{2}=-\operatorname{Id}_{\mathbb{E}}$.
\end{definition}

If $\mathcal{I}$ is a complex structure on $\mathbb{E}$ then $\mathcal{I}$ is
an isomorphism of $\mathbb{E}$ and this space can be provided with a structure of
complex Banach space $\mathbb{E}_{\mathbb{C}}$ defined by
\[
(\lambda+\mathrm{i}\mu)u=\lambda u+\mu\mathcal{I}(u).
\]
In this way, $\mathbb{E}_{C}$, as a real space, is isomorphic to $\mathbb{E}%
\oplus\mathbb{E}$. Moreover, $\mathcal{I}$ can be extended to a complex
isomorphism of $\mathbb{E}_{C}$, again denoted $\mathcal{I}$, which has two
complex eigenvalues $\mathrm{i}$ and $-\mathrm{i}$. Thus we have the
decomposition $\mathbb{E}_{C}=\mathbb{E}^{+}\oplus\mathbb{E}^{-}$ where
$\mathbb{E}^{+}$ is the eigenspace associated to $\mathrm{i}$ and
$\mathbb{E}^{-}$ is the eigenspace associated to $-\mathrm{i}$.

Moreover, the restriction of $-\mathcal{I}$ to $\mathbb{E}^{+}$ is an
isomorphism onto $\mathbb{E}^{-}$whose inverse is the restriction of
$-\mathcal{I}$ to $\mathbb{E}^{-}$. According to this decomposition,
$\mathcal{I}$ can be identify with the matrix
\[
\left(
\begin{array}
[c]{cc}%
\mathrm{i}.\operatorname{Id}_{\mathbb{E}^{+}} & 0\\
0 & -\mathrm{i}.\operatorname{Id}_{\mathbb{E}^{-}}%
\end{array}
\right)
\]

As well known, on the product of Banach spaces  $\hat{\mathbb{E}}={\mathbb{E}}\times{\mathbb{E}}$, we have a canonical complex structure $\mathcal{I}_\textrm{{can}}(u_1,u_2)=(-u_2,u_1)$. If ${\mathbb{E}}_1={\mathbb{E}}\times \{0\}$ and ${\mathbb{E}}_2=\{0\}\times {\mathbb{E}}$  and let $\iota_i:{\mathbb{E}}\longrightarrow \hat{{\mathbb{E}}}$ be the canonical isomorphism from ${\mathbb{E}}$ to ${\mathbb{E}}_i$, for $i \in \left\{1,2\right\}$.  Therefore we have $\mathcal{I}_\textrm{{can}}\circ \iota_1=-\iota_2$ and  $\mathcal{I}_\textrm{{can}}\circ \iota_2=-\iota_1$. According to the decomposition $\widehat{{\mathbb{E}}}={\mathbb{E}}_1\oplus {\mathbb{E}}_2$,  $\mathcal{I}_\textrm{{can}}$  can be identified with the matrix
\[
\begin{pmatrix} 0&- \operatorname{Id}_{{\mathbb{E}}_1}\cr
\operatorname{Id}_{{\mathbb{E}}_2}& 0
 \end{pmatrix}.
\]

\begin{definition}
\label{D_DecompositionComplexStructure} 
A complex structure $\mathcal{I}$ on a Banach space is decomposable\index{complex structure!decomposable} if there exists supplemented isomorphic subspaces
$\mathbb{E}_{1}$ and $\mathbb{E}_{2}$ such that $\mathcal{I}$ can be identify
with a matrix of type
\[%
\begin{pmatrix}
0 & -I\cr
 I^{-1} &0
\end{pmatrix}
.
\]
where $I$ is an isomorphism from $\mathbb{E}_{2}$ to $\mathbb{E}_{1}$.
\end{definition} 

The canonical complex structure on ${\mathbb{E}}\oplus {\mathbb{E}}$  is always  decomposable. \\

\begin{remark}
\label{R_NonDecomposableComplexStructure} 
There exist non decomposable complex structures on Banach spaces.  Indeed,  in \cite{Fer}, Theorem 3, 
%Uniqueness of complex structure and real hereditarily indecomposable Banach spaces Advances in Mathematics 213 (2007) 462?488
 the reader can find  an example of a real Banach space with only two complex structures not isomorphic and not decomposable. (see also \cite{Wen}, $\S$ 3.2). 
% (1995) REAL AND COMPLEX OPERATOR IDEALS, Quaestiones Mathematicae, 18:1-3, 271-285,) V.Ferenczi} %Uniqueness of complex structure and real hereditarily indecomposable Banach spaces
% Advances in Mathematics 213 (2007) 462?488
\end{remark}

For an even dimensional space ${\mathbb{E}}$, there always exists a complex structure and a para-complex structure which is decomposable and which  is isomorphic to the canonical one. However, in the Banach context  the existence or the uniqueness of the isomorphism class  of complex structures  is much more difficult than in the finite dimensional framework.  Indeed,  it is well known that there exist Banach spaces without complex structure (see for instance \cite{Die}). %Dieudonné, J. Complex structures on real Banach spaces. Proc. Amer. Math. Soc. (3) 1(1952), 162?164.}.
A real Banach space whose complexification is a primary space has at most one complex structure ( \cite{FeGa},  Theorem 22). %Ferenczi, V., Galego, E. Countable groups of isometries on Banach spaces. Trans. Amer. Math. Soc. (362) 8 (2010), 4385?4431
In particular, the classical spaces of sequences $c_0$,  $\ell_p$ or functions  $L^p([0,1])$ with  $1\leq p\leq \infty$ and $C([0,1] )$ have a unique complex structure. 
But there exist complex structures which are not necessary isomorphic to the canonical one.  Even more, there can exist an infinity number of isomorphic classes of complex structures on Banach spaces (for references about such different situations, see for instance \cite{Car}). %arXiv:1401.1781v1  [math.FA]  8 Jan 2014 A BANACH SPACE WITH A COUNTABLE INFINITE NUMBER OF COMPLEX STRUCTURES W. CUELLAR CARRERA}.
 Now,  the  existence of a decomposable complex structure on a Banach space ${\mathbb{E}}$ is equivalent to   the existence of a  decomposition  ${\mathbb{E}}={\mathbb{E}}_1\oplus {\mathbb{E}}_2$ where ${\mathbb{E}}_1$ and ${\mathbb{E}}_2$ are isomorphic.
 
Such a  Banach space will be called \textit{decomposable}.   %Note that in previous classical spaces  the complex structure is decomposable.  We have the following about decomposable complex structure

\subsubsection{Para-complex  structures}
\label{___ParaComplexStructures}
\begin{definition}
\label{D_ProductStructure}
 A product structure\index{product structure} on a Banach space $\mathbb{E}$ is an endomorphism $\mathcal{J}$ of
$\mathbb{E}$ such that $\mathcal{J}^{2}=\operatorname{Id}_{\mathbb{E}}$.
\end{definition}

If $\mathcal{J}$  is a product structure we have a decomposition $\mathbb{E}=\mathbb{E}^+\oplus\mathbb{E}^-$  where $\mathbb{E}^+$ (resp. $\mathbb{E}^-$) is the eigenspace  associated to the eigenvalue $+1$ (resp. $-1$). In general  $\mathbb{E}^+$ and  $\mathbb{E}^-$  are not  isomorphic.

\begin{definition}
\label{D_ParaComplexStructure}
A product structure $\mathcal{J}$ will be called a  para-complex structure if in the decomposition  $\mathbb{E}=\mathbb{E}^+\oplus\mathbb{E}^-$,  the subspaces $\mathbb{E}^+$ and  $\mathbb{E}^-$ are isomorphic. 
\end{definition}
In the one hand, a para-complex  structure $\mathcal{J}$ on $\mathbb{E}$,   can be written as a matrix of type
\[%
\begin{pmatrix}
\operatorname{Id}_{\mathbb{E}^+}&0\cr 0 & -\operatorname{Id}_{\mathbb{E}^-}
\end{pmatrix}
.\]

On the other hand, consider a decomposition $\mathbb{E}=\mathbb{E}_1\oplus \mathbb{E}_2$ such that  $\mathbb{E}_1$ and $\mathbb{E}_2$ are isomorphic subspaces of $\mathbb{E}$, then the operator  $\mathcal{J}$  is defined by  a matrix of type
\[
\begin{pmatrix}
0 & I\cr I^{-1} & 0
\end{pmatrix}
.\]
where  $I$, isomorphism from $\mathbb{E}_1$ onto $\mathbb{E}_2$, is a para-complex structure. In this case we have
\[ 
\mathbb{E}^+=\{u-I^{-1}u,\in u\in \mathbb{E}_1\}=\{Iv+v,\; v\in \mathbb{E}_2\}\\
\mathbb{E}^-=\{u-I^{-1}u,\in u\in \mathbb{E}_1\}=\{Iv-v,\; v\in \mathbb{E}_2\}.
\]
Now if $\mathcal{S}$  is the symmetry\index{symmetry} on $\mathbb{E}$  defined by the matrix 
\begin{eqnarray}\label{S}
\begin{pmatrix} -\operatorname{Id}_{{\mathbb{E}}_1}&0\cr
0&\operatorname{Id}_{{\mathbb{E}}_2}
 \end{pmatrix},
\end{eqnarray}
then $\mathcal{I}=\mathcal{S}\mathcal{J}$ is a decomposable complex structure on $\mathbb{E}$  associated to the decomposition $\mathbb{E}=\mathbb{E}_1\oplus \mathbb{E}_2$. \\
Conversely, if $\mathcal{I}$ is a decomposable complex structure and $\mathbb{E}=\mathbb{E}_1\oplus \mathbb{E}_2$ is an associated decomposition, $\mathcal{J}=\mathcal{S}\mathcal{I}$ is a para-complex structure on $\mathbb{E}$..
 
Thus  from this last argument,  on the product of Banach spaces  $\hat{\mathbb{E}}={\mathbb{E}}\times{\mathbb{E}}$, we have a "canonical" para-complex structure $\mathcal{J}_\textrm{can}(u_1,u_2)=(u_2,u_1)$. If ${\mathbb{E}}_1={\mathbb{E}}\times \{0\}$ and ${\mathbb{E}}_2=\{0\}\times {\mathbb{E}}$  and let $\iota_i:{\mathbb{E}}\longrightarrow \hat{{\mathbb{E}}}$ be the canonical isomorphism from ${\mathbb{E}}$ onto ${\mathbb{E}}_i$, for $i=1,2$.  Therefore, we have $\mathcal{J}_\textrm{can}\circ \iota_1=\iota_2$ and  $\mathcal{J}_{can}\circ \iota_2=\iota_1$. According to the decomposition $\widehat{{\mathbb{E}}}={\mathbb{E}}_1\oplus {\mathbb{E}}_2$,  $\mathcal{J}_{can}$  can be identified with the matrix
\[
\begin{pmatrix} 0&Id_{{\mathbb{E}}_1}\cr
Id_{{\mathbb{E}}_2}& 0
 \end{pmatrix}.
\]
%Let $g$ be a Krein inner product on $\E$ and consider an associated decomposition $\E=\E^+\oplus \E^-$. Then the associated fundamental symmetry $J$ is a linear product structure. More generally if $g$ is a Krein indefinite product and $\E=E_1\oplus \E_2$ is the associated decomposition then then if $u=u_1+u_2$    is the associated decomposition then $J(u)=u_1-u_2$ is also a linear product structure.
%On the other hand,
 
%\begin{definition}
%\label{D_DecompositionComplexStructure} A complex structure $\mathcal{I}$ on a
%Banach space is decomposable if there exists supplemented isomorphic subspaces
%$\mathbb{E}_{1}$ and $\mathbb{E}_{2}$ such that $\mathcal{I}$ can be identify
%with the matrix of type
%\[%
%\begin{pmatrix}
%0 & -J\cr J^{-1} & 0
%\end{pmatrix}
%.
%\]
%where $J$ is an isomorphism from $\mathbb{E}_{2}$ to $\mathbb{E}_{1}$
%\end{definition}
Now, if $\mathcal{I}$ (resp. $\mathcal{J}$) is a decomposable complex (resp. a para-complex) structure on ${\mathbb{E}}$ and if ${\mathbb{E}}={\mathbb{E}}_1\oplus {\mathbb{E}}_2$ is an associated decomposition, then we can associate to $\mathcal{I}$ (resp. $\mathcal{J}$) a tangent structure ${J}_\mathcal{I}$ (resp. ${J}_\mathcal{J}$) given by

$J_\mathcal{I}(u)=-\mathcal{I}(u_2)$ if $u=u_1+u_2$ and  $\ker J_\mathcal{I}={\mathbb{E}}_1.$

$(\textrm{resp. } J_\mathcal{J}(u)=\mathcal{I}(u_2)$ if $u=u_1+u_2$ and  $\ker J_\mathcal{I}={\mathbb{E}}_1).$ \\ 
Of course  ${J}_\mathcal{I} = {J}_\mathcal{J}$ if $\mathcal{I}=\mathcal{S}\mathcal{J}$.\\
Conversely, if $J$ is a tangent structure on $\mathbb{E}$ with $\mathbb{E}=\ker J\oplus \mathbb{K}$ and if $J$ is defined by the matrix
\[
\begin{pmatrix} 0&J_\mathbb{K}\cr
0& 0
 \end{pmatrix},
\]
then, in an evident way, we have a unique decomposable complex (resp. para-complex) structure $\mathcal{I}_J$ (resp. $\mathcal{J}_J$) associated to $J$; we also have $\mathcal{I}_J=\mathcal{S}\mathcal{J}_J$ with a clear  adequate operator $\mathcal{S}$.

In fact,  we have the following result about decomposable complex, para-complex and tangent structures:
\begin{proposition}
\label{P_IsomorphismDecomposableComplexStructure} 
Let $\mathcal{I}$ (resp. $\mathcal{J}$)  be a decomposable complex (resp. para-complex) structure on a Banach space $\mathbb{E}$. If $\mathbb{E}=\mathbb{E}_{1}\oplus\mathbb{E}_{2}$ is a decomposition for
$\mathcal{I}$  and $\mathcal{J}$ and $I$ is the associated isomorphism from $\mathbb{E}_{2}$ to
$\mathbb{E}_{1}$, then $A: \mathbb{E}\longrightarrow\mathbb{E}_{1}%
\oplus\mathbb{E}_{1}$ given by $A(u,v)=(u,Iv)$ is an isomorphism such that
$A^{*} \mathcal{I}_\textrm{can}=\mathcal{I}$ (resp.  $A^{*} \mathcal{J}_\textrm{can}=\mathcal{J}$)  where $\mathcal{I}_\textrm{can}$ (resp. $\mathcal{J}_\textrm{can}$)  is the
canonical complex (resp. para-complex) structure on $\mathbb{E}_{1}\oplus\mathbb{E}_{1}$.
\end{proposition}

\subsection{Compatibilities between different structures}

\label{__CompatibilityBetwenDifferentStructures}

\subsubsection{Tangent structure compatible with a cotangent structure}
\label{___TangentStructureCompatibleCotangentStructure}

\begin{definition}
\label{D_CompatibilityTangentStructureWithCotangentStructure}
A cotangent structure\index{cotangent structure} $\Omega$ on the Banach space $\mathbb{E}$ is compatible with a
tangent structure $J$ if $\ker J$ is Lagrangian and
\[
\forall\left(  u,v\right)  \in\mathbb{E}^{2},\ \Omega(Ju,v)+\Omega(u,Jv)=0.
\]
\end{definition}

Given a decomposition $\mathbb{E}=\mathbb{L}\oplus\mathbb{K}$ associated to a
tangent structure $J$ where $\mathbb{L}=\ker J$, the restriction
$J_{\mathbb{K}}$ of $J$ to $\mathbb{K}$ is an isomorphism $J_{K}%
:\mathbb{K}\rightarrow\mathbb{L}$ and, in this decomposition, $J$ can be
written as a matrix of the type $\left(
\begin{array}
[c]{cc}%
0 & J_{\mathbb{K}}\\
0 & 0
\end{array}
\right)  $. According to such a decomposition, to a cotangent structure
$\Omega$ is associated an operator $\Omega^{\flat}$ which can be written as a matrix
$\left(
\begin{array}
[c]{cc}%
0 & \Omega_{\mathbb{LK}}\\
\Omega_{\mathbb{KL}} & \Omega_{\mathbb{KK}}%
\end{array}
\right)  $. With these notations, $\Omega$ is compatible with $J$ if and only
if we have 
\[
\Omega^{\flat}\circ J+J^{\ast}\circ\Omega^{\flat}=0
\] %which, in the
%reflexive is equivalent to
%\[
%\Omega^{\flat}\circ\Omega_{\mathbb{KL}}\circ J_{K}=(\Omega_{\mathbb{LK}}\circ
%J_{K})^{\ast}.
%\]
The following proposition gives a link between tangent structures, cotangent
structures and indefinite inner products.

\begin{proposition}
\label{P_CompatibilityBetweenTangentCotangentStructuresAndIndefineInnerProduct}
Let $\mathbb{E}$ be a Banach space. Then we have
\begin{enumerate}
\item[1.] 
Assume that there exists a tangent structure $J$ on $\mathbb{E}$ and
an indefinite inner product $g$ on $\ker J$. Then there exists a cotangent
structure $\Omega$ on $\mathbb{E}$ compatible with $J$ such that $\ker J$ is a
weak Lagrangian space for $\Omega$.

\item[2.] 
Assume that there exists a cotangent structure $\Omega$ on
$\mathbb{E}$ and let $\mathbb{E}=\mathbb{L}\oplus\mathbb{K}$ be an associated
decomposition. Then there exists a tangent structure $J$ on $\mathbb{E}$ such
that $\ker J=\mathbb{L}$.
%and  there exists a indefinite product on $\mathbb{L}$ such that  $g(u,v)=\O(\mathcal{I}u,v)-\O(u, \mathcal{I}v)$ where $\mathcal{I}$ is the complex structure associated to $J$ as previously.
\item[3.] 
Assume that there exists a tangent structure $J$ on $\mathbb{E}$
which is compatible with a cotangent structure $\Omega$. Then, if
$\mathbb{E}=\ker J\oplus\mathbb{K}$, $\ker J$ is a weak Lagrangian space for
$\Omega $. Moreover, $g(u,v)=\Omega(Ju,v)$ is an indefinite inner product on $\ker
J$
\end{enumerate}
\end{proposition}

\subsubsection{Tangent structure, cotangent structure and indefinite inner product}
\label{___TangentCotangentStructureIndefiniteInnerProduct}

\begin{proposition}
\label{P_CompatibilityTangentCotangentInnerProduct} Let $\mathbb{E}$ be a
Banach space.

\begin{enumerate}
\item If there exists a tangent structure $J$ on $\mathbb{E}$ and an
indefinite inner product on $\ker J$, then there exists a cotangent structure
$\Omega$ on $\mathbb{E}$ compatible with $J$ such that $\ker J$ is a weak
Lagrangian space for $\Omega$.

\item If there exists a cotangent structure $\Omega$ on $\mathbb{E}$ and if
$\mathbb{E}=\mathbb{L}\oplus\mathbb{K}$ is an associated decomposition, then
there exists a tangent structure $J$ on $\mathbb{E}$ such that $\ker
J=\mathbb{L}$.

\item Assume that there exists a tangent structure $J$ on $\mathbb{E}$ which
is compatible with a cotangent structure $\Omega$ where $\mathbb{E}=\ker
J\oplus\mathbb{K}$ is an associated decomposition ($\ker J$ is a weak
Lagrangian space for $\Omega$). Then $g(u,v)=\Omega(Ju,v)$ is an indefinite
inner product on $\ker J$.
\end{enumerate}
\end{proposition}

\subsubsection{Symplectic structures, inner products and  decomposable complex structures }
\label{___SymplecticInnerProductDecomposableComplexStructures}

\begin{definition}
\label{D_CompatibilitySymplecticStructurePreHilbertProductComplexStructures}
Let $\mathbb{E}$ be a Banach space.

\begin{enumerate}
\item We say that a weak symplectic structure $\Omega$ and a complex structure 
$\mathcal{I}$ on $\mathbb{E}$ are compatible if $(u,v)\mapsto\Omega
(u,\mathcal{I}v)$ is a pre-Hilbert product on $\mathbb{E}$ and%
\[
\forall\left(  u,v\right)  \in\mathbb{E}^{2},\ \Omega(\mathcal{I}%
u,\mathcal{I}v)=\Omega(u,v)
\]

\item We say that a pre-Hilbert product $g$ and a complex structure 
$\mathcal{I}$ on $\mathbb{E}$ are compatible, if
\[
\forall\left(  u,v\right)  \in\mathbb{E}^{2},\ g(\mathcal{I}u,\mathcal{I}%
v)=g(u,v)
\]

\item We say that a pre-hilbert product $g$ and a weak symplectic structure
$\Omega$ on $\mathbb{E}$ are compatible if $\mathcal{I}=(g^{\flat})^{-1}%
\circ\Omega^{\flat}$ is a well defined complex structure or para-complex structure on $\mathbb{E}$.\\
\end{enumerate}
\end{definition}

 If $\mathbb{E}$ is isomorphic to a Hilbert space and $\Omega$  is a symplectic form then there exist a compatible complex structure $\mathcal{I}$ on $\mathbb{E}$ and an Hilbert  product compatible with $\mathcal{I}$ and $\Omega$ (see for instance \cite{ChMa}).  In our context we have:

\begin{proposition}
\label{P_CompatibilitySymplecticStructureInnerProductDecomposableComplexStructures}
Let $\mathbb{E}$ be a Banach space.
\begin{enumerate}
\item[1.]  
Assume  that  there exists a   decomposable complex structure $\mathcal{I}$ on $\mathbb{E}$ and a pre-Hilbert  product $g$ on 
$\mathbb{E}$ which are compatible.  Then there exists a weak symplectic structure $\Omega$ on $\mathbb{E}$ given by 
{$\Omega( u,v)={g}(\mathcal{I}u, v)$ for all $u,v\in \mathbb{E}$}
which  is compatible with $\mathcal{I}$. Moreover if  $\mathbb{E}=\mathbb{E}_1\oplus \mathbb{E}_2$ is a decomposition such that $\mathbb{E}_1$ %(resp. $ \mathbb{E}_2$ is the eigenspace of $\mathcal{J}$ associated to the eigenvalue $+1$ (resp. $-1$ such that $\mathbb{E}_1$
 and $\mathbb{E}_2$ are $g$-orthogonal, then $\Omega$ is a Darboux form and $\mathbb{E}_1$ and $\mathbb{E}_2$ are supplemented Lagrangian subspaces.
\item[2.]  
Assume that there exists a pre-hilbert inner product $g$ on $\mathbb{E}$. Let $\widehat{\mathbb{E}}$ the Hilbert space which is the completion of the pre-Hilbert space  $(\mathbb{E},g)$. For any weak symplectic structure $\Omega$ (not necessarily compatible with $g$) such that  $\Omega^{\flat}(\mathbb{E})\subset g^{\flat}(\mathbb{E})$ assume that   $\Omega^{\flat}(\mathbb{E})$ is dense in $\widehat{\mathbb{E}}$. Then there exists  a Hilbert product $\bar{g}$  and  a decomposable complex structure $\mathcal{I}$ on $\mathbb{E}$ which are  compatible with $\Omega$. Moreover there exists a decomposition $\mathbb{E}=\mathbb{E}_1\oplus \mathbb{E}_2$  associated to $\mathcal{I}$ such that $\mathbb{E}_1$ and $\mathbb{E}_2$ are $\bar{g}$-orthogonal, then $\Omega$ is a Darboux form and $\mathbb{E}_1$ and $\mathbb{E}_2$ are supplemented Lagrangian subspaces.
\item[3.] 
Assume that there exists a decomposable complex structure $\mathcal{I}$ on $\mathbb{E}$ compatible with a weak symplectic structure $\Omega$. Then  $g(u,v)=\Omega(u, \mathcal{I}v)$ is a neutral inner product on $\mathbb{E}$ which is compatible with $\mathcal{I}$. Moreover,   there exists a decomposition $\mathbb{E}=\mathbb{E}_1\oplus \mathbb{E}_2$ which is associated to $\mathcal{I}$ and to $\Omega$ then $\mathbb{E}_1$ and $\mathbb{E}_2$ are $g$-orthogonal and Lagrangian.
\end{enumerate}
\end{proposition}

\begin{proof}
1. Let $\mathcal{I}$ be a  complex structure compatible with a pre-Hilbert product $g$. Then if we set $\Omega( u,v)={g}(\mathcal{I}u, v)$ for all $u,v\in \mathbb{E}$ then $\Omega$ is a continuous skew symmetric bilinear form. Since $g$ is weakly non degenerate, so is $\Omega$. Then the space $\mathbb{E}_g$ has a Hilbert structure and we can extend $g$ to a Hilbert product $\hat{g}$ on $\mathbb{E}_g$. From the compatibility of $\mathcal{I}$ with $g$, it follows that $\mathcal{I}$ is an isometry for the norm $||\;||_g$; So we can extend ${\mathcal{I}}$ to a complex structure  $\hat{\mathcal{I}}$ on $\mathbb{E}_g$ and  so $\Omega$ can also be extended by the same formula to a symplectic form $\hat{\Omega}$ on $\mathbb{E}_g$. By classical results (cf. \cite{ChMa})    %P.  Chernoff   and  J.  Marsden,  Properties  of  Infinite   Dimensional  Hamiltonian  Systems, Springer-Verlag  Lecture  Notes  in  Mathematics  425  (1974).  
 in a Hilbert space, the announced result is true in $\mathbb{E}_g$. Now  since $\mathcal{I}$ leaves $\mathbb{E}$ invariant, the announced results  are also true in restriction to $\mathbb{E}$.\\

2. Fix some pre-Hilbert product $g$ in $\mathbb{E}$ and some symplectic structure $\Omega$ on $\mathbb{E}$ such that  $\Omega^\flat(\mathbb{E})\subset g^\flat(\mathbb{E})$ and  $\Omega^\flat(\mathbb{E})$ is dense in the completion $\widehat{\mathbb{E}}$ of the pre-Hilbert space $(\mathbb{E},g)$. We denote by $\widehat{g}$ the extension of $g$ to $\widehat{\mathbb{E}}$. Since $\Omega$ is bounded on $\mathbb{E}$ (relative to its original norm $||\;||$) and the inclusion of  $(\mathbb{E},g)$ in $(\mathbb{E},||\;||)$ is continuous, we can extend $\Omega$ to a $2$-form $\widehat{\Omega}$ on $\widehat{\mathbb{E}}$. We set $A=(g^\flat)^{-1}\circ \Omega^\flat$. Thus $A$ is an injective endomorphism of $\mathbb{E}$ and we have 
\begin{eqnarray}\label{eq_OgA}
g(Au,v)=<\Omega^\flat(u),v>=\Omega(u,v).
\end{eqnarray}
Thus, if we have $A^2=-\operatorname{Id}_\mathbb{E}$ then the proof is completed.\\
However, in general, it is not true. Therefore, we shall build from $A$  a new endomorphism which will be the announced complex structure. At first, for the sake of simplicity, we denote by $\ll\;,\;\gg$ the Hilbert product on $\mathbb{E}_g$ defined by $g$. From our assumptions on $\Omega$, it follows that the range of $A$ is dense. 
 Thus we  can define the adjoint $A^*:\widehat{\mathbb{E}}^*\equiv \widehat{\mathbb{E}}\longrightarrow \widehat{\mathbb{E}}$ and, by definition, for all $u, v\in \mathbb{E}$  we have:
\[
\ll A^*u,v\gg=\ll u,Av\gg=\Omega(v,u)=-\Omega(u,v)=\ll -Au,v\gg.
\]
and since $A$ is injective with dense range, $A^*$ is an isomorphism. Now the extension  $\widehat{A}$ of $A$ to $\widehat{\mathbb{E}}$ is well defined  and we have 
$\widehat{A}=(\widehat{g}^\flat)^{-1}\circ \widehat{\Omega}^\flat$. Since the restriction of $\widehat{A}$ to $\mathbb{E}=-A$ we must have $A^*=\widehat{A}^*$ and so $\widehat{A}$ is an isomorphism and in particular $\widehat{\Omega}$ is a strong symplectic form.

%Again from our assumption, this implies that the restriction of $A^*$ to $\mathbb{E}$ is $-A$ and so $A$ can be extended to $\widehat{\mathbb{E}}$ by an isomorphism  $\widehat{A}$ of $\widehat{\mathbb{E}}$. 
%Thus, from the relation (\ref{OgA}), $\Omega$ can be also extended to a  non degenerate skew symmetric bilinear form  on $\mathbb{E}_g$ also denoted $\Omega$.  %We claim claim that $\bar{g}(u,v)=\ll AA^*u,v\gg$ is a Hilbert product on $\E_g$. Indeed 
Now we have 
\begin{eqnarray}\label{eq_AA*}
\ll \widehat{A}\widehat{A}^*v,u\gg=\ll \widehat{A}^*v,\widehat{A}^*u\gg=\ll v,\widehat{A}\widehat{A}^*u\gg.
\end{eqnarray}
This implies that $\widehat{A}{A}^*$ is symmetric. From (\ref{eq_AA*}) we also have 
\[
\forall u \in \widehat{\mathbb{E}}, u\not=0, \; \ll \widehat{A}\widehat{A}^*u,u\gg=\ll \widehat{A}^*u,\widehat{A}^*u\gg \;> 0.
\]  
Because  $\widehat{A}^*=-\widehat{A}$, $\widehat{A}\widehat{A}^*$ is a symmetric  auto-adjoint definitive positive operator. Therefore there exits $\widehat{R}$ such that 
$\widehat{R}^2=\widehat{A}\widehat{A}^*$ which is auto-adjoint. Moreover, since $\widehat{A}$ and $\widehat{A}^*$ leave $\mathbb{E}$ invariant the same is true for $\widehat{R}$. 
We set $\widehat{\mathcal{I}}=\widehat{R}^{-1}\widehat{A}$.\\

In finite dimension (cf. \cite{Can}, $\S$ 12) as in this Hilbert context (cf. \cite{Wei} proof of Proposition  5.1), we can show that $\widehat{\mathcal{I}}$ commutes with $\widehat{A}$ and $\widehat{R}$, $\widehat{\mathcal{I}}$ is compatible  with $\widehat{g}$, and we have the following properties:
\[
\widehat{\mathcal{I}}\widehat{\mathcal{I}}^*=Id,\; \widehat{\mathcal{I}}^*=-\widehat{\mathcal{I}},\; \widehat{\mathcal{I}}^2= -Id \;\textrm{ and } \widehat{\Omega}(\widehat{\mathcal{I}}u,\widehat{\mathcal{I}}v)=\widehat{\Omega}(u,v)
\]
We can remark that  $\widehat{\Omega}(u,\widehat{\mathcal{I}}v)= \widehat{g}(\widehat{R}u,v)$ and so $\bar{g}$ defined by $\overline{\widehat{g}}(u,v)=\widehat{g}(\widehat{R}u,v)$ is a Hilbert product compatible with $\widehat{\Omega}$.  Since $\widehat{R}$ and $\widehat{A}^*$ leaves $\mathbb{E}$ invariant , the same is true for $\widehat{\mathcal{I}}$ and so we have the same property in restriction to $\mathbb{E}$.\\
 
 On the other hand since $\widehat{\Omega}$ is a symplectic form on the Hilbert space $\widehat{\mathbb{E}}$ we have a decomposition $\widehat{\mathbb{E}}=\widehat{\mathbb{L}}\oplus \widehat{ \mathbb{L}}'$ in two orthogonal Lagangian spaces relative to the Hilbert product defined by $\ll\;,\;\gg$  (cf. proof of Proposition 5.1 in \cite{Wei}). As $\widehat{\mathcal{I}}=(\overline{\widehat{g}}^\flat)^{-1} \circ\widehat{ \Omega}^\flat$, the restriction of $\mathcal{I}$ to $\mathbb{L}=\mathbb{E}\cap \widehat{L}$ is an isomorphism onto $\mathbb{L}'=\mathbb{E}\cap \widehat{L}'$ and so $\mathcal{I}=\widehat{\mathcal{I}}_{| \mathbb{E}}$ is decomposable.  %Now since $A$ and $A^*$ leaves $\mathbb{E}$ invariant, the same is true for $R$ and so for $\mathcal{I}$. This implies   all the previous properties  are valid in restriction to  $\mathbb{E}$.\\
 
3. Consider a complex structure $\mathcal{I}$ on $\mathbb{E}$ compatible with a weak symplectic structure $\Omega$ on $\mathbb{E}$. We set $g(u,v)=\Omega(u, \mathcal{I}v)$.  From the definition of compatibility between $\mathcal{I}$ and $\Omega$, it follows that $g$ is symmetric and  positive definite  and so is a pre-Hilbert product.
 
  Now $g(\mathcal{I}u,\mathcal{I}v)=\Omega(\mathcal{I}u, -v)=\Omega(v,\mathcal{I}u)=g(v,u)=g(u,v)$ and so $g$ is compatible with $\mathcal{I}$. The result is then a consequence of Point 1.\\
 \end{proof}

 % \begin{rem}\label{contextJOg}\normalfont 
 % It follows from the previous proof that in  context of Point 1, Point 2 and  Point 3,  the Banach space $\E$ is continuously at densely embedded  in a Hilbert space $\mathbb{H}$
   %and $\O$, $g$ and $\mathcal{I}$  are  the respective  restriction to $\E$ of strong Darboux form, Hilbert product and complex structure on $\mathbb{H}$. 
    % \end{rem}
 
Consider the direct sum $\mathbb{H}_{0}=\mathbb{L}\oplus\mathbb{L}$, where
$\mathbb{L}$ is a Banach space, provided with the Hilbert product
\[
<u_{1}+u_{2},v_{1}+v_{2}>_{\mathbb{H}_{0}}=<u_{1},v_{1}>_{\mathbb{L}}%
+<u_{2},v_{2}>_{\mathbb{L}}.
\]
We denote by $\Omega_{D}$ the canonical Darboux form on $\mathbb{L}%
\oplus\mathbb{L}$, $g_{\mathrm{can}}$ the previous Hilbert product and
$\mathcal{I}_{\mathrm{can}}$ the unique complex structure on $\mathbb{H}_{0}$
compatible with $g_{\mathrm{can}}$ and $\Omega_{D}$.

\begin{corollary}
\label{C_DarbouxFormPreHilbertProductComplexStructure} 
Consider a weak Darboux form $\Omega$, a pre-Hilbert product $g$ and a decomposable complex structure $\mathcal{I}$ on a Banach space $\mathbb{E}$. Assume that any pair of such a triple is
compatible. Then the third one is compatible with anyone of the given
pair.\newline There exists a Hilbert space $\mathbb{H}$ such that $\mathbb{E}$
is continuously and densely embedded in $\mathbb{H}$ and an isomorphism
$A_{\mathbb{H}}:\mathbb{H}\rightarrow\mathbb{H}_{0}=\mathbb{L}\oplus
\mathbb{L}$ such that the restriction of $A$ of $A_{\mathbb{H}}$ to
$\mathbb{E}$ has the following properties:
\[
A^{\ast}\Omega_{D}=\Omega,A^{\ast}g_{\mathrm{can}}=gA^{\ast}\mathcal{I}%
_{\mathrm{can}}=\mathcal{I}.
\]
\end{corollary}

\begin{remark}
\label{R_DecompositionInLagrangianOrthogonalspaces}  
Under the assumption of Corollary \ref{C_DarbouxFormPreHilbertProductComplexStructure}, we have a decomposition $\mathbb{E}=\mathbb{E}_1\oplus \mathbb{E}_2$ such that $\mathbb{E}_1$ and $ \mathbb{E}_2$ are isomorphic, Lagrangian and orthogonal subspaces of $\mathbb{E}$
\end{remark}

\begin{definition}
\label{D_WeakKahlerBanachSpace} A Banach space $\mathbb{E}$ is called a  weak  (resp. strong) K\"ahler Banach space if there exists on $\mathbb{E}$ a weak (resp. strong) Darboux  form $\Omega$, a weak (resp. strong)  inner product $g$ and a decomposable complex structure such that any pair of such a triple data is
compatible
\end{definition}

\subsubsection{Symplectic structures, neutral inner products and  para-complex structures }
\label{___SymplecticInnerProductParaComplexStructures}

\begin{definition}
\label{D_CompatibilitySymplecticStructureNeutralInnerProductParaComplexStructures}
Let $\mathbb{E}$ be a Banach space.${}$

\begin{enumerate}
\item 
We say that a weak symplectic structure $\Omega$ and a para-complex structure 
$\mathcal{J}$ on $\mathbb{E}$ are compatible if $(u,v)\mapsto\Omega
(u,\mathcal{J}v)$ is a neutral inner product on $\mathbb{E}$ and%
\[
\forall\left(  u,v\right)  \in\mathbb{E}^{2},\ \Omega(\mathcal{J}%
u,\mathcal{J}v)=-\Omega(u,v)
\]
\item 
We say that a neutral inner product $g$ and a para-complex  structure 
$\mathcal{J}$ on $\mathbb{E}$ are compatible, if
\[
\forall\left(  u,v\right)  \in\mathbb{E}^{2},\ g(\mathcal{J}u,\mathcal{J}%
v)=-g(u,v)
\]
\item 
We say that a neutral inner product $g$ and a weak symplectic structure
$\Omega$ on $\mathbb{E}$ are compatible if $\mathcal{J}=(g^{\flat})^{-1}%
\circ\Omega^{\flat}$ is a well defined para-complex structure  on $\mathbb{E}$.
\end{enumerate}
\end{definition}

\begin{remark}
\label{R_AlmostDecomposableComplexAlmostParaComplex} 
Let  $\mathcal{J}=\mathcal{S}\mathcal{I}$ be the para-complex structure naturally associated to a decomposable complex structure $\mathcal{I}$. If $\Omega$  is a weak symplectic form  on $\mathbb{E}$, then $\Omega$ and $\mathcal{J}$ are compatible if and only if $\Omega$ and $\mathcal{I}$ are compatible. If $g$ is a pre-Hilbert product on $\mathbb{E}$, we denote by $g_\mathcal{S}$ the neutral inner product defined by $g_{\mathcal{S}}(u,v)=-g(\mathcal{S}u, v)$. Then $g_\mathcal{S}$ and $ \mathcal{J}$ are compatible if and only if $g$ and $\mathcal{I}$ are compatible.
\end{remark}

According to this Remark, from Proposition \ref{P_CompatibilitySymplecticStructureInnerProductDecomposableComplexStructures}, we obtain:
\begin{proposition}
\label{P_CompatibilitySymplecticStructureInnerProductParaComplexStructures}
Let $\mathbb{E}$ be a Banach space.

\begin{enumerate}
\item[1.]  
Assume  that  there exists a para-complex structure $\mathcal{J}$ on $\mathbb{E}$ and a neutral product $g$ on 
$\mathbb{E}$ which are compatible.  Then there exists a weak symplectic structure $\Omega$ on $\mathbb{E}$ given by 
{$\Omega( u,v)={g}(\mathcal{J}u, v)$ for all $u,v\in \mathbb{E}$}
which  is compatible with $\mathcal{J}$. Moreover if  $\mathbb{E}=\mathbb{E}_1\oplus \mathbb{E}_2$ is a decomposition  associated with $\mathcal{J}$ such that $\mathbb{E}_1$ and $\mathbb{E}_2$ are $g$-orthogonal, then $\Omega$ is a Darboux form and $\mathbb{E}_1$ and $\mathbb{E}_2$ are supplemented Lagrangian subspaces.
\item[2.]  
Assume that there exists a neutral inner product $g$ on $\mathbb{E}$. Consider an inner product $<\;,\;>$ on $\mathbb{E}$ canonically associated to a decomposition $\mathbb{E}=\mathbb{E}^+\oplus \mathbb{E}^-$
 associated to  $g$ and let $\widehat{\mathbb{E}}$ the Hilbert space associate do $g$. For any weak symplectic structure $\Omega$ (not necessarily compatible with $g$) such that  $\Omega^{\flat}(\mathbb{E})\subset g^{\flat}(\mathbb{E})$. 
 If $\Omega^{\flat}(\mathbb{E})$ is dense in  $\widehat{\mathbb{E}}$, then there exists  a neutral product $\bar{g}$  and  a para-complex structure $\mathcal{J}$ on $\mathbb{E}$ 
 which are  compatible with $\Omega$. Moreover there exists a decomposition $\mathbb{E}=\mathbb{E}_1\oplus \mathbb{E}_2$  associated to $\mathcal{J}$ such that  the 
 restriction of $\bar{g}$ to $\mathbb{E}_1$ (resp. $\mathbb{E}_2$) is positive definite (resp. negative definite), $\Omega$ is a Darboux form and $\mathbb{E}_1$ and $\mathbb{E}_2$  are supplemented Lagrangian subspaces.
\item[3.] 
Assume that  there exists a  para-complex structure $\mathcal{J}$ on $\mathbb{E}$ compatible with a weak symplectic structure $\Omega$. Then $g(u,v)=\Omega(u, \mathcal{I}v)$ is a neutral product on $\mathbb{E}$ which is compatible with $\mathcal{J}$. Moreover, there exists a decomposition $\mathbb{E}=\mathbb{E}_1\oplus \mathbb{E}_2$ which is associated to $\mathcal{J}$ and to $\Omega$ such that  $\mathbb{E}_1$ and $\mathbb{E}_2$ are Lagrangian and  the restriction of ${g}$ to $\mathbb{E}_1$ (resp. $\mathbb{E}_2$) is positive definite (resp. negative definite).
\end{enumerate}
\end{proposition}

Consider the direct sum $\mathbb{H}_{0}=\mathbb{L}\oplus\mathbb{L}$, where
$\mathbb{L}$ is a Banach space, provided with the Hilbert  product
\[
<u_{1}+u_{2},v_{1}+v_{2}>_{\mathbb{H}_{0}}=<u_{1},v_{1}>_{\mathbb{L}}%
+<u_{2},v_{2}>_{\mathbb{L}}.
\]
We denote by $\Omega_{D}$ the canonical Darboux form on $\mathbb{L}%
\oplus\mathbb{L}$,  $g_{\mathrm{can}}$ is the canonical neutral  inner product on $\mathbb{H}_0$ defined by

 $g_{\mathrm{can}}=<\;,\;>$ on the first factor $\mathbb{L}$ and  $g_{\mathrm{can}}=<\;,\;>=-<\;,\;>$ on the second factor $\mathbb{L}$ and $\mathcal{J}_{\mathrm{can}}$ is the unique para-complex structure on $\mathbb{H}_{0}$ compatible with $g_{\mathrm{can}}$ and $\Omega_{D}$.

\begin{corollary}
\label{C_DarbouxFormNeutralInnerProductParaComplexStructure} 
Consider a weak Darboux form $\Omega$, a neutral inner product $g$ and a  para-complex structure $\mathcal{J}$ on a Banach space $\mathbb{E}$. Assume that any pair of such a triple is
compatible. Then the third one is compatible with anyone of the given
pair.\\
There exist a Hilbert space $\mathbb{H}$ such that $\mathbb{E}$
is continuously and densely embedded in $\mathbb{H}$ and an isomorphism
$A_{\mathbb{H}}:\mathbb{H}\rightarrow\mathbb{H}_{0}=\mathbb{L}\oplus
\mathbb{L}$ such that the restriction of $A$ of $A_{\mathbb{H}}$ to
$\mathbb{E}$ has the following properties:
\[
A^{\ast}\Omega_{D}=\Omega,A^{\ast}g_{\mathrm{can}}=g 
\textrm{ and } 
A^{\ast}\mathcal{I}_{\mathrm{can}}=\mathcal{J}.
\]
\end{corollary}

\begin{remark}
\label{R_DecompositionInLagrangianOrthogonalspacesParaComplex} 
As in the previous section, under the assumption of Corollary \ref{C_DarbouxFormNeutralInnerProductParaComplexStructure}, we have a decomposition $\mathbb{E}=\mathbb{E}_1\oplus \mathbb{E}_2$ such that $\mathbb{E}_1$ and $ \mathbb{E}_2$ are isomorphic, Lagrangian and the restriction of ${g}$ to $\mathbb{E}_1$ (resp. $\mathbb{E}_2$) is positive definite (resp. negative definite)
\end{remark}

\begin{definition}
\label{D_WeakPara-KahlerBanachSpace} 
A Banach space $\mathbb{E}$ is called a  weak  (resp. strong) para-K\"ahler Banach space if there exists on $\mathbb{E}$ a weak (resp. strong) Darboux  form $\Omega$, a weak (resp. strong)  neutral inner product $g$ and a  para-complex structure such that any pair of such a triple data is
compatible
\end{definition}

\bigskip

\noindent{\small {\textsc{Unit\'e Mixte de Recherche 5127 CNRS, Universit\'e
de Savoie Mont Blanc\newline Laboratoire de Math\'ematiques (LAMA) }%
\newline\textit{Campus Scientifique \newline73370 Le Bourget-du-Lac, France}%
}\newline
\noindent{\textsc{E-mail address: }
\textit{patrickcabau@gmail.com}, \textit{fernand.pelletier@univ-smb.fr}} }


\begin{thebibliography}{9999}                                                                                            

\bibitem[AbMa]{AbMa}M. Abbati, A. Mani\`{a}, \textit{On Differential Structure for Projective Limits of Manifolds,} J. Geom. Phys. \textbf{29} 1-2 (1999) 35--63.

\bibitem[ADGS]{ADGS}M. Aghasi, C.T. Dodson, G.N. Galanis, A. Suri,
\textit{Conjugate connections and differential equations on infinite
dimensional manifolds,} J. Geom. Phys. (2008).

\bibitem[AgSu]{AgSu}M. Aghasi, A. Suri, \textit{Splitting theorems for the
double tangent bundles of Fr\'{e}chet manifolds}, Balkan Journal of Geometry
and Its Applications \textbf{15} 2 (2010) 1--13.

\bibitem[Bam]{Bam}D. Bambusi, \textit{On the Darboux theorem for weak
symplectic manifolds}, Proccedings AMS Vol 127-11(1999) 3383--3391.

\bibitem[Bel1]{Bel1}D. Belti\c{t}\u{a}, \textit{Integrability of Analytic Almost Complex 
Structures on Banach Manifolds}, Annals of Global Analysis and Geometry 28(1):59--73 2005.

\bibitem[Bel2]{Bel2}D. Belti\c{t}\u{a}, \textit{Smooth Homogeneous Structures in
Operator Theory}, Chapman \& Hall/CRC Monographs and Surveys in Pure and
Applied Mathematics, vol. 137, Chapman \& Hall/CRC Press, Boca Raton-London-New York-Singapore, 2006.

\bibitem[Bog]{Bog}J. Bognar,\textit{\ Indefinite Inner Product Spaces}, Springer-Verlag, Berlin Heidelberg New York (1974).

\bibitem[Bom]{Bom}J. Boman,\textit{\ Differentiability of a function and of
its compositions with functions of one variable}, Mathematica Scandinavica 20
(1967) 249--268.

\bibitem[Bou]{Bou}N. Bourbaki, \textit{Vari\'{e}t\'{e}s diff\'{e}rentiables et
analytiques. Fascicule de r\'{e}sultats}, \S \S 1--7, Hermann, Paris, 1967.

\bibitem[Cab]{Cab}P. Cabau, \textit{Strong projective limits of Banach Lie
algebroids}, Portugaliae Mathematica, Volume 69, Issue 1 (2012).

\bibitem[CaPe]{CaPe}P. Cabau, F. Pelletier, \textit{Integrability of direct limits of Banach manifolds}, accepted for publication in Annales de la Facult\'{e} des Sciences de Toulouse.

\bibitem[Can]{Can}Ana Cannas da Silva, \textit{Lectures on Symplectic Geometry} Springer Lecture Notes
in Mathematics 1764, second printing (2008).

\bibitem[Car]{Car}W. C. Carrera, \textit{A Banach space with a countable infinite number of complex structures}, Journal of Functional Analysis Vol 267, Issue 5, 1 September 2014.
%arXiv:1401.1781v1  [math.FA]  8 Jan 2014 A BANACH SPACE WITH A COUNTABLE INFINITE NUMBER OF COMPLEX STRUCTURES W. CUELLAR CARRERA.

\bibitem[ChMa]{ChMa}P. Chernoff, J. Marsden, \textit{Properties of Infinite Dimensional Hamiltonian Systems}, Springer-Verlag Lecture Notes in Mathematics 425 (1974).

\bibitem[ChSt]{ChSt}D. Chillingworth, P. Stefan, \textit{Integrability of
singular distributions on Banach manifolds}, Math. Proc. Cambridge Phil. Soc. 79 (1976) 117--128.

\bibitem[ClGo]{ClGo}R. S. Clark, D. S. Goel, \textit{Almost cotangent manifolds,} J. Differential Geom. 9 (1974) 109--122.

\bibitem[CrTh]{CrTh}M. Crampin, G. Thompson, \textit{Affine bundles and integrable almost tangent structures}, Math. Proc. Cambridge Philos. Soc. 98 (1985).

\bibitem[Dah]{Dah}R. Dahmen, \textit{Direct limit constructions in infinite
dimensional Lie theory}, Dissertation Universit\"{a}t Paderborn. http://nbn-resolving.de/urn:nbn:de:hbz:466:2-239.

\bibitem[Die]{Die}J. Dieudonn\'e, \textit{Complex structures on real Banach spaces}, Proc. Amer. Math. Soc. (3) 1 (1952).

\bibitem[DGV]{DGV}C.T.J. Dodson, G. Galanis, E. Vassiliou, \textit{Geometry in
a Fr\'{e}chet Context: A projective Limit Approach}, London Mathematical Society, Lecture Notes Series 428. Cambridge University Press, 2015.

\bibitem[EgWu]{EgWu}M. Egeileh, T. Wurzbacher, \textit{Infinite-Dimensional Manifolds as Ring Spaces}, Publications of the Research Institute for Mathematical Sciences, Volume 53, Issue 1, 2017, 187--209.

\bibitem[Eli]{Eli}H.I. Eliasson, \textit{Geometry of manifolds of maps}, J. Differential Geometry
1 (1967) 169--194.  

\bibitem[Fer]{Fer} V. Ferenczi, \textit{Infinite-Dimensional Manifolds as Ring Spaces}, Advances in Mathematics 213 (2007).

\bibitem[FeGa]{FeGa} V. Ferenczi, E. Galego, \textit{Countable groups of isometries on Banach spaces}, Trans. Amer. Math. Soc. 362-8 (2010).

\bibitem[FiTe]{FiTe}D. Filipovi\'{c}, J. Teichmann, \textit{Existence of invariant manifolds for stochastic equation in infinite dimension}, Journal of Functional Analysis 197 (2003) 398--432.

\bibitem[FrKr]{FrKr}A. Fr\"{o}licher, A. Kriegl, \textit{Linear Spaces and Differentiation Theory}, Pure and Applied Mathematics, J. Wiley, Chichester 1988.

\bibitem[Gal1]{Gal1}G.N. Galanis, \textit{Projective Limits of Banach-Lie groups}, Periodica Mathematica Hungarica \textbf{32} (1996) 179--191.

\bibitem[Gal2]{Gal2}G.N. Galanis, \textit{Projective Limits of Banach Vector Bundles}, Portugaliae Mathematica \textbf{55} 1 (1998) 11--24.

\bibitem[Gal3]{Gal3}G.N. Galanis, \textit{Differential and geometric structure for the tangent bundle of a projective limit manifold}, Rend. Sem. Mat. Univ. Padova \textbf{112} (2004) 103--115.

\bibitem[Glo]{Glo}H. Gl\"{o}ckner, \textit{Direct limits of infinite-dimensional Lie groups compared to direct limits in related categories}, Journal of Functional Analysis 245 (2007) 19--61.

\bibitem[Gra]{Gra}J. Grabowski, \textit{Derivative of the Exponential Mapping for Infinite Dimensional Lie Groups}, Annals of Global Analysis and Geometry \textbf{11} issue 3 (1993) 213--220.

\bibitem[Ham]{Ham}R.S. Hamilton, \textit{The Inverse Function Theorem of Nash and Moser}, Bulletin of the American Mathematical Society \textbf{7} 1 (1982) 65--222.

\bibitem[Jar]{Jar}H. Jarchow, \textit{Locally convex spaces}, Teubner, Stutgart, 1981.

\bibitem[KaSw]{KaSw}N. J. Kalton R. C. Swanson, \textit{A symplectic Banach space with no Lagrangian subspaces},  Transactions of the American Mathematical Society, Vol. 273, No. 1 (Sep., 1982), 385--392. 

\bibitem[Klo]{Klo}M. Klotz, \textit{The Automorphism Group of a Banach Principal Bundle with }$\left\{  1\right\} -$\textit{structure}, Geometriae Dedicata, October 2011, Volume 154, Issue 1, 161--182.

\bibitem[Kob]{Kob}S. Kobayashi, \textit{Transformation groups in differential
geometry}, Ergebnisse der Mathematik und ihrer Grenzgebiete, Band 70, Berlin: Springer, 1972.

\bibitem[KrMi]{KrMi}A. Kriegel, P.W. Michor, \textit{The convenient Setting of
Global Analysis }(AMS Mathematical Surveys and Monographs) \textbf{53}, 1997.

\bibitem[Kum1]{Kum1} P. Kumar, \textit{Almost complex structure on Path space},
 International Journal of Geometric Methods in Modern Physics Vol. 10, No. 3 (2013).

\bibitem[Kum2]{Kum2}P. Kumar, \textit{ Darboux chart on projective limit of weak
symplectic Banach manifold}, Int. J. Geom. Methods Mod. Phys., 12, (2015).

\bibitem[Lan]{Lan}S. Lang, \textit{Fundamentals of Differential Geometry},
Springer, New York, 1999.

\bibitem[Les]{Les}J.A. Leslie, \textit{On a differential structure for the
group of diffeomorphisms}, Topology \textbf{46} (1967) 263--271.

\bibitem[Lib]{Lib}P. Libermann, \textit{Sur les structures presque complexes et autres structures infinit\'{e}simales réguli\`{e}res}, Bulletin de la S. M. F., tome 83 (1955) 195--224.

\bibitem[Mar]{Mar}J. E. Marsden \textit{ Lectures on Geometric Methods in
Mathematical Physics}, SIAM, (1981).

\bibitem[MaRa]{MaRa}J.E. Marsden, T.S. Ratiu, \textit{Introduction to Mechanics and
Symmetry}, Texts in Applied Mathematics 17, Springer-Verlag New-York, 1999.

\bibitem[Mol]{Mol}P. Molino, \textit{Sur quelques propri\'{e}t\'{e}s des }
$G$\textit{-structures}, Journal of Differential Geometry 7 (1972) 489--518.

\bibitem[Nee]{Nee}K-H Neeb,\textit{Towards a Lie theory for locally convex
groups}, arXiv:1501.06269 [math.RT].

\bibitem[Omo]{Omo}H. Omori, \textit{Infinite-dimensional Lie groups},
Translations of Mathematical Monographs vol 158 (American Mathematical
Society) 1997.

\bibitem[OnVi]{OnVi}A.L. Onishchik, E.B. Vinberg, \textit{Lie Groups and Lie
Algebras III: Structure of Lie Groups and Lie Algebras}, Encyclopaedia of
Mathematical Sciences, vol. 41. Springer, Berlin 1994.

\bibitem[Pat]{Pat}I. Patyi, \textit{On the }$\overline{\partial}$\textit{-equation on a Banach space}, Bull. Soc. math. France, 128 (2000) 391--406.

\bibitem[Pay]{Pay}S. Paycha, \textit{Basic Prerequisites in Differential
Geometry end Operator Theory in View of Applications to Quantum Field Theory},
Preprint Universit\'{e} Blaise Pascal, Clermont, France, 2001.

\bibitem[Pel1]{Pel1}F. Pelletier, \textit{Integrability of weak distributions
on Banach manifolds}, Indagationes Mathematicae 23 (2012) 214--242.

\bibitem[Pel2]{Pel2}F. Pelletier, \textit{On a Darboux theorem for symplectic
forms on direct limits of ascending sequences of Banach manifolds}, arXiv:1807.04501.

\bibitem[PiTa]{PiTa}P. Piccione, D. V. Tausk, \textit{The theory of
connections and G-structures. Applications to affine and isometric
immersions}, XIV escola de geometria diferencial -- XIV school of differential
geometry, Salvador, Brazil, July 17--21, 2006.

\bibitem[RoRo]{RoRo}A. P. Robertson, W. Robertson, \textit{Topological Vector
Spaces}, Cambridge Tracts in Mathematics. 53. Cambridge University Press 1964.

\bibitem[Sau]{Sau}D. J. Saunders, \textit{The Geometry of Jet Bundles},
\textbf{142}, London Mathematical Society Lectures Notes Series, Cambridge
University Press, 1989.

\bibitem[ThSc]{ThSc}G. Thompsonandu, U. Schwardmann, \textit{Almost tangent
and cotangent structure in the large}, Trans AMS vol 327 N1 (1991) 313--328.

\bibitem[Tum1]{Tum1}B. Tumpach, \textit{Hyperk\"{a}hler structures and infinite-dimensional Grassmannians}, Journal of Functional Analysis 243 (2007) 158--206.

\bibitem[Tum2]{Tum2}B. Tumpach, \textit{Structures k\"{a}hlériennes et hyperk\"{a}hlériennes en dimension infinie}, Th\`{e}se, Ecole Polytechnique, 2005.  

\bibitem[Vai]{Vai}I. Vaisman, \textit{Lagrange geometry on tangent manifolds},
Intern. J. of Math. and Math. Sci., 2003(51) (2003) 3241--3266.

\bibitem[Wei]{Wei}A. Weinstein, \textit{Symplectic Manifolds and Their
Lagrangian Submanifolds}, Advances in Math. 6 (1971) 329--346.

\bibitem[Wen]{Wen}J. Wenzel, \textit{Real and complex operator Ideals}, Quaestiones Mathematicae, 18 (1995).

\end{thebibliography}
\end{document}